\documentclass{article}
\usepackage[english]{babel}
\usepackage[T1]{fontenc}
\usepackage[utf8]{inputenc}
\usepackage{amsmath,amssymb,amsfonts,amsthm}
\usepackage[initials,nobysame]{amsrefs}
\usepackage{bm}
\usepackage{mathtools}
\usepackage{ascmac}
\usepackage{framed}
\usepackage{ulem,cancel}
\usepackage{graphicx}
\usepackage{bbm}
\usepackage{tikz}
\usepackage{pgfplots}
\usetikzlibrary{intersections, calc, arrows.meta, patterns}
\pgfplotsset{compat=1.18}
\allowdisplaybreaks

\usepackage{color}								% You may remove in the future manuscript.
\definecolor{lightgray}{rgb}{0.75,0.75,0.75}

\usepackage{titlesec}
\usepackage{authblk}

 {\endMakeFramed}

\newcommand{\eps}{\varepsilon}
\newcommand{\C}{\mathbb{C}}
\newcommand{\bbH}{\mathbb{H}}

\newcommand{\R}{\mathbb{R}}
\renewcommand{\S}{\mathbb{S}}

\newcommand{\Z}{\mathbb{Z}}
\newcommand{\boA}{\mathcal{A}}

\newcommand{\boD}{\mathcal{D}}
\newcommand{\boE}{\mathcal{E}}

\newcommand{\boH}{\mathcal{H}}
\newcommand{\boI}{\mathcal{I}}
\newcommand{\boJ}{\mathcal{J}}
\newcommand{\boK}{\mathcal{K}}

\newcommand{\boM}{\mathcal{M}}
\newcommand{\boN}{\mathcal{N}}

\newcommand{\boR}{\mathcal{R}}

\newcommand{\la}{\lambda}

\newcommand{\al}{\alpha}
\newcommand{\be}{\beta}
\newcommand{\gam}{\gamma}

\newcommand{\si}{\sigma}

\newcommand{\vph}{\varphi}

\DeclareMathOperator{\Ran}{{\rm Ran}}
\renewcommand{\Im}{\mathop{{\rm Im}}\nolimits}
\DeclareMathOperator{\Ker}{{\rm Ker}}
\renewcommand{\Re}{\mathop{{\rm Re}}\nolimits}

\DeclareMathOperator{\Span}{{\rm Span}}

\DeclareMathOperator{\diag}{{\rm diag}}

\DeclareMathOperator{\ind}{{\rm ind}}
\DeclareMathOperator{\sech}{{\rm sech}}

\newcommand{\ess}{\mathrm{ess}}

\newtheorem*{claim*}{Claim}

\newtheorem*{cor*}{Corollary}
\newtheorem{lem}{Lemma}[section]
\newtheorem{prop}[lem]{Proposition}

\newtheorem{thm}{Theorem}[section]

\theoremstyle{definition}
\newtheorem{definition}{Definition}[section]

\newtheorem*{remark}{Remark}

\theoremstyle{remark}

\makeatletter												%
  \@addtoreset{equation}{section}						%
\makeatother

\newcommand{\nor}[2]{\left\| {#1} \right\|_{#2}}		%	Norm
\newcommand{\jb}[1]{\langle {#1} \rangle{} }			%	Japanese blacket
\newcommand{\ovl}[1]{\overline{#1}}					%	Complex conjugate
		
\newcommand{\inp}[2]{\left\langle {#1} , {#2} \right\rangle }	%	inner product

\newcommand{\del}{{\delta}}								%	Small delta
\newcommand{\rd}{{\partial}}								%	Round d
\newcommand{\nab}{{\nabla}}							%	Nabla

\newcommand{\sym}{\mathrm{sym}}

\newcommand{\bh}{{\mathbf{h}}}

\newcommand{\boe}{{\mathbf{e}}}
\newcommand{\bom}{{\mathbf{m}}}

\newcommand{\bow}{{\bm{w}}}

\newcommand{\bxi}{{\bm{\xi}}}
\newcommand{\bF}{{\mathbf{F}}}

\begin{document}

\title{Asymptotic stability of Bloch wall for dispersive Landau--Lifshitz equation
}
\author{Slim Ibrahim, Christof Melcher, and Ikkei Shimizu}

% \affil[1]{Department of Mathematics and Statistics, University of Victoria, BC-Canada}
% \affil[2]{}
% \affil[3]{Department of Mathematics, Graduate School of Science, Kyoto University, Kyoto, Japan}

% \author{Slim Ibrahim
% , Christof Melcher, and Ikkei Shimizu}
\date{}
\maketitle

%\section{Asymptotic stability under odd perturbation: $\eps>0$}
\begin{abstract}
% We consider the Bloch wall stationary solutions to the dispersive Landau--Lifshitz equation with Dzyaloshinskii--Moriya interaction
% with coupling \sout{parameter} {\color{blue}constant} $\eps$
% . 
% %{\color{red}I would request the improvement of this expression, as the sign of $\eps$ is also important.}
% First, we establish a dichotomy for the orbital stability and instability depending on the value of $\eps$. 
% %We prove that when $\eps>0$, the Bloch wall is orbitally stable. 
% Then, under symmetric perturbations, we prove that if $\eps>0$, the Bloch wall is asymptotically stable.
%
%
%{\color{blue}
We consider the Bloch wall stationary solutions to the dispersive Landau--Lifshitz equation with Dzyaloshinskii--Moriya interaction (DMI). 
First, we establish a dichotomy for the orbital stability and instability depending on the value of a coupling constant.
%the relative strength of DMI
Then, under symmetric perturbations, we prove that the Bloch wall is asymptotically stable in the energetically favorable regime.
%}
% \textit{2020 MSC code}: 35B40, 82D40, 35Q55, 47A40, 47A75\\
% \textit{Keywords}: Landau--Lifshitz equation, Bloch wall, asymptotic stability, Virial-type estimate
\end{abstract}

\begingroup
\renewcommand{\thefootnote}{}
\footnotetext{%

}
\endgroup

\section{Introduction and Main Results}

\subsection{Physical background and the Landau--Lifshitz model}
%{\color{teal} %Alternatively start: 
We consider domain wall configurations $\bom:\R \to \mathbb{S}^2$ in an anisotropic chiral ferromagnet governed by the energy 
\begin{align}\label{eq:Energy}
E_\eps(\bom) = \int_\R \frac 12 |\rd_x \bom|^2 +\eps (\rd_x \bom \times \bom)\cdot \hat\boe_1 + \frac 12 (1-m_3^2)\,dx.
\end{align}
Here, $\eps$ is a coupling constant related to the relative strength of  Dzyaloshinkii--Moriya interaction (DMI), a spin-orbit effect arising in magnetic systems without inversion symmetry. The specific variant considered here comes from the three-dimensional helicity density $\bom \cdot (\nabla \times \bom)$. 
This chiral coupling breaks the mirror symmetry of the energy, selecting domain walls of a definite handedness and, for sufficiently large $|\eps|$, destabilizing the uniform ferromagnetic ground state in favor of helically modulated phases in which the magnetization winds periodically through space.\par
The effective field, i.e., minus the variational gradient of $E_\eps$,
given by
\[
\bh_\eps (\bom) = \rd_{xx} \bom -2\eps \hat\boe_1\times \rd_x\bom + m_3 \hat\boe_3, 
\]
%}
%\textcolor{blue}{projection to tangent plane?}
%{\color{teal}
exerts a torque on the magnetization vector,
%}
%{\color{blue}
driving a continuum precessional dynamics. 
The corresponding Hamiltonian flow is governed by the Landau--Lifshitz equation (LL): 
%formulated by
% The Hamiltonian flow for  nonequilibrium 
% magnetization 
% states, called the Landau--Lifshitz equation (LL), are governed by the continuum precessional dynamics driven by the torque effect due to the interaction with the effective field; namely,
\begin{align}\label{E5.1Mar}
%\label{E5.1Mar}
\rd_t \bom = \bom \times \bh_\eps(\bom)
\end{align}
%governed by the
where the magnetization fields $\bom=\bom(t,x):\R\times\R\to\S^2$ is the unknown function for \eqref{E5.1Mar}. 
%Formal computation shows that 
The flow preserves the energy functional \eqref{eq:Energy}, 
%. Here, 
and the states for which the torque $\bom \times \bh_\eps$ vanishes correspond to the stationary solutions to \eqref{E5.1Mar}. 
For any $\eps\in\R$, \eqref{E5.1Mar} admits an $\eps$-independent exact stationary solution, known as the \textit{Bloch wall}, 
%. This is a topological soliton, 
explicitly given by
\[
\bom_0(x) = 
\begin{pmatrix}
    0 \\ b(x) \\ p(x)
\end{pmatrix}
,\qquad 
b(x) =\sech x,\qquad p(x) =\tanh x.
\]
%
%}
%
% {\color{gray}
% The Landau--Lifshitz (LL) equation is a fundamental model in condensed matter physics, describing the continuum precessional dynamics of magnetization in ferromagnetic materials. We consider the one-dimensional LL equation:
% \begin{align}\label{E5.1Mar}
% \rd_t \bom = \bom \times \bh_\eps(\bom).
% \end{align}
% Here, the unknown $\bom=\bom(t,x)\colon \R\times \R\to\mathbb{S}^2$ represents the magnetization vector field on the unit sphere $\mathbb{S}^2\subset \R^3$. The parameter $\eps\in\R$ characterizes specific magnetic anisotropies or chiral interactions. 

% The effective magnetic field $\bh_\eps$ is the Fr\'echet derivative of the magnetic energy functional $E_\eps$, given by:
% \begin{align*}
% \bh_\eps (\bom) &= \rd_{xx} \bom -2\eps \hat\boe_1\times \rd_x\bom + m_3 \hat\boe_3 = E_\eps'(\bom),\\
% E_\eps(\bom) &= \int_\R \frac 12 |\nab \bom|^2 +\eps (\rd_x \bom \times \bom)\cdot \hat\boe_1 + \frac 12 (1-m_3^2)\,dx.
% \end{align*}
% For any $\eps\in\R$, \eqref{E5.1Mar} admits an $\eps$-independent exact stationary solution, known as the \textit{Bloch wall}. This topological soliton is explicitly given by:
% \[
% \bom_0(x) = 
% \begin{pmatrix}
%     0 \\ b(x) \\ p(x)
% \end{pmatrix}
% ,\qquad 
% b(x) =\sech x,\qquad p(x) =\tanh x.
% \]
% }
%{\color{teal}
The Bloch wall is a paradigmatic example of a magnetic kink structure. It originally 
emerges from the competition between exchange and magnetic anisotropy and can be characterized 
as an absolute energy minimizer subject to boundary conditions connecting two antipodal vacua 
(see e.g. \cites{KMM2006,MR4616656}). The variational problem features a double degeneracy due to translation 
$\bom_0 \mapsto \bom_0(\cdot + x_0)$ and gauge invariance $\bom_0 \mapsto e^{i \varphi} \bom_0$ 
with respect to rotations about the anisotropy axis. The additional effect of the helicity term 
is to lift rotational degeneracy. Indeed $
    E_\eps(e^{i \varphi} \bom_0) = 2 - \eps \pi \cos \varphi$,
so the characteristic Bloch wall $\varphi=0$ is energetically preferable for $\eps>0$ (while reflections 
of $\bom_0$ are preferred in the opposite regime $\eps<0$). In fact it can be deduced from 
our analysis that, for $\eps>0$, the Hessian is strictly positive modulo the translation mode 
$\bom_0'$ and in turn $\bom_0$ is a local energy minimum (adapting Section 5 in \cite{Li_Melcher2018}). The energy expansion 
also marks the transition point $\eps=2/\pi$ between the ferromagnetic regime where 
$E_\eps \ge 0$ and the absolute energy minimum is attained for $m_3 \equiv \pm 1$ and the 
helical regime where periodically modulating configurations become the ground state 
\cite{Bogdanov_Hubert}. \par
%
%\textcolor{teal}{
%Very recently, 
A different magnetic system, a magnetic nanowire, was considered in \cite{MR4616656}, where the anisotropy axis is aligned with the propagation direction of the wall. In this setting, the interplay with helicity preserves rotational invariance about the propagation axis, giving rise to modulated static domain walls and precessing dynamic domain walls, both available in closed form. These works address the stability of such solitons in the so-called ferromagnetic regime of small relative DMI strength, where the ground state is the uniform state aligned with the wire. In the present work, we consider a domain wall model arising in the multidimensional context of a uniaxially anisotropic chiral magnet. Here, modulation and precession are suppressed, but in contrast, the traveling wall profiles are not available in closed form \cite{MR3957215}. %\cite{Komineas_Venakides_Melcher_2019}
Moreover, we address the entire range of DMI strength, including the helical regime of large coupling.\par
%}
%
% {\color{blue}Nevertheless of the above energetics, since \eqref{E5.1Mar} is an energy-conserving flow, the dynamics in the vicinity of the Bloch wall state is quite a nontrivial question. In particular, it is worth examining %the asymptotic stability which stems from
% the relaxation to the Bloch wall caused by 
% the underlying dispersive effect in \eqref{E5.1Mar}. 
% In light of this, the primary objective of this paper is to rigorously examine the stability of the Bloch wall $\bom_0$ under the dispersive flow of \eqref{E5.1Mar}. 
% % the dynamical stability of $\bom_0$ under \eqref{\eqref{E5.1Mar}} is a mathematically nontrivial question. 
% % In fact, \eqref{\eqref{E5.1Mar}}
% }
% {\color{gray}
% The primary objective of this paper is to rigorously examine the stability of the Bloch wall $\bom_0$ under the dispersive flow of \eqref{E5.1Mar}.}
% \bigskip
%{\color{teal} 
%Reformulated paragraph:
Despite the energetic considerations discussed above, the dynamics in the vicinity of the Bloch wall state remain a subtle problem, owing to the energy-conserving nature of the flow governed by \eqref{E5.1Mar}. In particular, understanding how the dispersive effects inherent in \eqref{E5.1Mar} influence the long-time behavior of perturbations and potentially drive the system toward the Bloch wall is a fundamental question. The central objective of this paper is to establish a rigorous stability theory for the Bloch wall $\bom_0$ under the dispersive dynamics of \eqref{E5.1Mar}, thereby elucidating the interplay between energy conservation, chiral symmetry breaking, dispersion, and the persistence of this coherent structure.
%}

\subsection{Mathematical Backgrounds and Challenges}
%
%
% {\color{blue}
% When the energy is \textcolor{teal}{reduced to} the Dirichlet integral $\int_\R \frac 12 |\bom'|^2 dx$, then the corresponding equation is particularly called the Schr\"odinger maps, the name of which arises from the feature that the skew-symmetric operation $\bom\times$ is exactly the complex structure of the sphere. 
% In fact, the classical trick known as the Hasimoto transform reduces the equation to cubic nonlinear Schr\"odinger equations (NLS), which reveals the nature of an integrable system for the Schr\"odinger maps.\par
% Not specified to our case, the dispersive Landau--Lifshitz flow has been studied extensively in many decades [xxxxxx]. As confronted in these works, the mathematical analysis of \eqref{E5.1Mar} still presents formidable challenges, situating this problem at the forefront of research in nonlinear dispersive PDEs. Moreover, \eqref{E5.1Mar} poses its specific difficulties arising from the DMI term. Let us describe the detailed challenges as follow.
% }
%{\color{teal} 
%Alternative text:
When the energy is reduced to the Dirichlet functional
$\int_\R \frac 12 |\rd_x \bom|^2 dx$,
the associated geometric evolution is known as the Schrödinger map equation. Its name reflects the fact that the skew-symmetric operator $\bom\times$, acting on the tangent space of the sphere, induces the natural complex structure on $\mathbb S^2$. This geometric structure places Schrödinger maps in close connection with dispersive Schrödinger dynamics. Indeed, through the classical Hasimoto transform, the equation can be reformulated as a cubic nonlinear Schrödinger equation, revealing its integrable structure in one spatial dimension.\par
A closely related geometric interpretation is provided by the binormal flow, which describes the evolution of a space curve in the direction of its binormal vector, with velocity proportional to its curvature. Through the identification of the tangent vector of the curve with a map into the sphere, the binormal flow is equivalent to the one-dimensional Schrödinger map equation. In this sense, Schrödinger maps may be viewed as a geometric counterpart of the free Schrödinger equation, in which the linear dispersive mechanism is encoded by the curvature of the evolving map, while the nonlinear effects arise from the underlying geometric constraint.\par
Beyond this integrable setting, dispersive Landau--Lifshitz flows have been extensively investigated over the past decades, see \cites{MR866199,MR2304153,MR1877231,MR1737504,MR2800718,MR3090782,MR3090180,MR2725187} and references therein. Despite substantial progress, the rigorous analysis of \eqref{E5.1Mar} remains challenging, owing to the intricate interplay between nonlinear geometric effects and dispersive dynamics. In addition to the difficulties inherent in the classical Landau--Lifshitz equation, the presence of the Dzyaloshinskii--Moriya interaction introduces further structural complications that require new analytical ideas. We now describe the principal difficulties specific to the problem under consideration.

%%%%%%%%SLIM START

\paragraph{1. Asymptotic stability in dispersive PDEs.} 
Establishing asymptotic stability for solitons in purely dispersive PDEs relies on global-in-time smoothing estimates for linearized operators in matrix form, which require stringent spectral assumptions \cites{MR2219305, MR4701407}, as well as an intricate analysis of nonlinear interactions between bounded modes and radiation \cites{MR4864876,MR4430297,MR2843104}. 
For two-dimensional Schr\"odinger maps, in fact, 
the Hasimoto transform reduces the system to a standard NLS, allowing for global bounds via Strichartz estimates \cites{MR2462113, MR2725187, MR1737504}. However, for \eqref{E5.1Mar}, the derivative nonlinearity induced by the DMI presents a severe obstruction and might be a source of ill posedness
in Sobolev spaces \cites{MR5076732, Christ}. 
An approach as in \cite{MR4482034} to handle this term as magnetic potential also confronts the difficulty, due to the potential presence of long-range factors.

% For two-dimensional Schr\"odinger maps, stability under an equivariant ansatz can be studied via the Hasimoto transform \cites{MR2462113, MR2725187, MR1737504}, 
% {\color{blue}by which the global-in-time bound can be obtained via standard Strichartz regime. However, for \eqref{E5.1Mar}, 

% the derivative nonlinearity produced by DMI causes serious obstruction. 
% Sorry for overlapping the editing. Let me undo my revision, by ikkei 
%}
% {\color{gray}
% which reduces the system to a standard NLS. 
% However, for \eqref{E5.1Mar}, the DMI term fundamentally obstructs this approach. Under the Hasimoto transform, DMI manifests as a quadratic derivative nonlinearity, which is known to cause ill-posedness in standard Sobolev spaces \cites{MR5076732, Christ}. While local-in-time well-posedness can be achieved by treating the DMI as a magnetic potential \cite{MR4482034}, extending this to a long-time analysis requires global dispersive estimates for magnetic Schr\"odinger operators. This remains highly nontrivial due to the potential presence of long-range factors.
% }

\paragraph{2. Quasilinear derivative nonlinearity.} 
Because the Hasimoto transform is not feasible, we must analyze \eqref{E5.1Mar} directly. The interplay between the complex structure $\bom\times$ and the effective field $\bh_\eps$ renders the problem quasilinear. The analysis of quasilinear dispersive systems is notoriously delicate \cites{MR2096797,MR2955206,MR4953177}, as the lack of maximal regularity prevents perturbative treatments of second-order derivative nonlinearities. In the long-time regime, these nonlinearities also disrupt the virial-type arguments recently developed in \cite{MR3630087}. Furthermore, establishing Strichartz estimates for Schr\"odinger operators with solution-dependent coefficients, following the approach in \cite{MR2955206}, is inherently obstructed by long-range factors that preclude global-in-time estimates.

\paragraph{3. Absence of dissipation.} 
Equation \eqref{E5.1Mar} is purely dispersive and conservative, preserving both the $L^2$-norm and the Hamiltonian $E_\eps$. This contrasts sharply with the dissipative Landau--Lifshitz--Gilbert (LLG) equation, where damping mechanisms facilitate proofs of asymptotic stability for magnetic solitons \cites{MR4616656,MR2895286,IS1}. Without dissipation, we cannot adapt these methods and must rely strictly on much weaker dispersive decay mechanisms. Handling spatial derivative nonlinearities in this setting requires highly sophisticated functional and frequency analysis, as well as modified energy methods \cites{MR1230709,MR2278007,MR4881654, ABIN25}.

%%%%%%%%%SLIM END

\subsection{Main Results}

%{\color{gray}
%\sout{To circumvent some of these obstacles}}, 
%We begin by introducing solution classes. Define
We define the complete metric space
\begin{align*}
\boM &\coloneqq \{ \bom:\R\to\mathbb{S}^2\ |\ \bom-\bom_0\in H^2(\R)\},
%\\
% \boM_\sym &\coloneqq \{\bom\in\boM\ |\ m_1,m_3\colon\text{odd}, m_2\colon\text{even}\},
\end{align*}
%In the present manuscript, we
and consider the solutions in $C_t\boM$ that are continuous in time into $\boM$. 
%and $C_t\boM_\sym$, where  we denote $C_tX$ by the space of functions continuous in time into complete metric space $X$. 
An interesting subclass $\boM_\sym$ that 
features the symmetry of \eqref{E5.1Mar} and $\bom_0$ is given by
\[
\boM_\sym \coloneqq \{\bom\in\boM\ |\ m_1,m_3\colon\text{odd}, m_2\colon\text{even}\}.
\]
 %$C_t\boM_\sym$, where

%and the subclass, 

% {\color{gray}
% First, we denote $C_tX$ by the space of functions continuous in time into complete metric space $X$. 
% %
% We first study 
% the solutions in the class $C_t\boM$, where
% %continuous in time into the function space $\boM$, where
% %natural finite-energy class:
% \[
% \boM \coloneqq \{ \bom:\R\to\mathbb{S}^2\ |\ \bom-\bom_0\in H^2(\R)\}.
% \]
% Second, observing that $\bom_0$ naturally resides in a symmetric class, we also investigate refined long-time stability under symmetric perturbations:
% \[
% \boM_\sym \coloneqq \{\bom\in\boM\ |\ m_1,m_3\colon\text{odd}, m_2\colon\text{even}\}.
% \]
% }
We first show the existence and uniqueness of solutions. We have %local well-posedness of 
%\subsubsection{Local Well-Posedness}
%establishes 
%the local-in-time theory for 
%\eqref{E5.1Mar} in %the energy space 
%C_t\boM$. 

\begin{thm}
System \eqref{E5.1Mar} is locally well-posed in $C_t\boM$. Moreover, if 
\[
\sup_{t\in I_{\max}}\nor{\bom(t)-\bom_0}{H^1(\R)} <\infty,
\]
where $I_{\max}$ is the maximal time interval of existence in $C_t\boM$, then $I_{\max}=\R$.
\end{thm}

%\subsubsection{Orbital Stability}
Next, we examine the orbital stability of $\bom_0$. Recall
%, measured up to spatial shifts.

\begin{definition}[Orbital stability]
$\bom_0$ is \textit{orbitally stable} in $\boM$ (resp. \textit{stable} in $\boM_\sym$) if there exist $C,\del_1>0$ such that if $\del=\inf_{y\in\R} \nor{\bom(0,\cdot) - \bom_0 (\cdot-y)}{H^1} <\del_1$ (resp. $\del=\nor{\bom(0,\cdot) - \bom_0}{H^1} <\del_1$), then the corresponding solution in $C_t\boM$ is global in time, and satisfies
\[
\inf_{y\in\R} \nor{\bom(t,\cdot) - \bom_0 (\cdot-y)}{H^1} <  C\del\qquad
(\text{resp. }\nor{\bom(t,\cdot) - \bom_0}{H^1} <  C\del)
\]
for all $t\in\R$. 
Otherwise, we say that $\bom_0$ is orbitally unstable (resp. \textit{unstable} in $\boM_\sym$).
\end{definition}

Then, our second result reveals a threshold value of the coupling constant $\eps$ for the orbital stability and instability of $\bom_0$.
%the dependence of orbital stability and instability on the coupling constant $\eps$.
% \dlt{Stability heavily depends on the sign of the anisotropy parameter $\eps$ and the symmetry of the perturbation.}

\begin{thm}[Orbital stability]\label{Tm1}
There exists a critical threshold $\eps_*<0$ such that the following hold:
\begin{itemize}
\item If $\eps>0$, then $\bom_0$ is orbitally stable in $\boM$.
\item If $\eps <\eps_*$, then $\bom_0$ is orbitally unstable in $\boM$.
\item If $\eps>\eps_*$, then $\bom_0$ is stable in the symmetric class $\boM_\sym$.
\end{itemize}
\end{thm}

\begin{remark}
%\dlt{(i)} 
The critical value $\eps_*$ corresponds to the exact threshold at which the spectral stability of the linearized operator changes (see Proposition \ref{Tx1} for details). 
A numerical evaluation yields $\eps_* \approx -1.3$ (see Section \ref{SA.6}).
This illustrates a nontrivial gap from the threshold of the energetic stability at $\eps=0$. It reveals the parameter regime $\eps\in[\eps_*,0]$ where $\bom_0$ is energetically unstable despite remaining linearly stable.
% \\
% \dlt{(ii) Through numerical evaluations of the linear operator's kernel, $\eps_* \approx -1.3$.} \\
% \dlt{(iii) We also note that traveling wave solutions for the 1D Landau--Lifshitz equation have been constructed and analyzed in related contexts; see Komineas, Melcher, and Venakides} \cite{KMV2019}.\\
% (iii) \textcolor{red}{As shown in Proposition \ref{Tx1}, the threshold of spectral stability for the ``linearlized operator'' 
% One can observe the nontrivial gap of two thresholds: that of the energetic stability at $\eps=0$ and that of spectral stability  
% of the energetic stability $\eps=0$ and that of the 
% Threshold of spectral stability is at $\eps=\eps_*$: BW remains spectrally stable in the energetically unstable parameter regime on $[\eps_*,\eps]$.}
\end{remark}

%\subsubsection{Asymptotic Stability}
Finally, we address the asymptotic stability of $\bom_0$ under symmetric perturbations via local energy decay. 
Our third and main result is:
\begin{thm}[Asymptotic stability]\label{Tx1.2Jun}
Let $\eps>0$, and let $\bom (t,x)$ be the solution to \eqref{E5.1Mar} with initial data $\bom(0,\cdot) \in \boM_\sym$. 
There exist $0<\del_0\ll 1$ and $A\gg 1$ such that if $\nor{\bom (0,\cdot)-\bom_0}{H^1}<\del_0$, then the solution to \eqref{E5.1Mar} is global in time, and satisfies the local decay estimate:
\begin{equation}\label{Eq:Jun4}
\int_0^\infty \left( \nor{b_A \rd_x(\bom(t) -\bom_0)}{L^2_x}^2 + \nor{b_A(\bom(t)-\bom_0)}{L^2_x}^2 \right) dt <\infty,
\end{equation}
where the weight function is given by $b_A(\cdot) \coloneqq \sech\left(\frac{\cdot}{A}\right)$.
\end{thm}

The integrability condition \eqref{Eq:Jun4} implies asymptotic stability in a weak sense. It guarantees the existence of 
%it guarantees 
a sequence of times $t_j \to \infty$ along which the local $H^1$-norm of the perturbation converges to zero, reflecting the radiation of dispersive energy away from the soliton.
\subsection{Key ideas and technical novelties}
\subsubsection*{Moving frame reduction (Section 2)}
Throughout the paper, we write $\bom$ in the specific moving frame $\{\hat\boe_1,\bom_0^\perp,\bom_0\}$: %in \eqref{eq:Sep2}
%, that is,
\begin{align}\label{eq:Sep3}
\bom = w_1 \hat\boe_1 +w_2 \bom_0^\perp + w_3 \bom_0,
\end{align}
which transforms $\bom$ into a new $\S^2$-valued unknown $\bow ={}^t (w_1,w_2,w_3)$. 
% Then the resulting equation of $\bow$ is
% \begin{gather}
% \rd_t \bm{w} =  \bm{w}\times \bbH \bm{w} -2(b-\eps) w_1\bm{w}' -2bpw_2\, \hat\boe_3\times \bm{w},\label{E2.4Aug}\\
% \bm{w} = \begin{pmatrix}
%     u_1\\u_2\\1+u_3
% \end{pmatrix}
% ,\qquad \bbH = \diag (L_\eps , L_0,-\rd_{xx}).\notag
% \end{gather}
Alternatively, we write the difference from $\bom_0$ as
% we reduce $\bom-\bom_0$ to the corresponding coordinates with respect to the specific moving frame $\{\hat\boe_1,\bom_0^\perp,\bom_0\}$, that is,
\begin{equation}\label{eq:Sep1}
\bom-\bom_0 = u_1 \hat\boe_1+ u_2 \bom_0^\perp+ u_3\bom_0.
\end{equation}
As long as $\nor{\bom-\bom_0}{H^1}\ll 1$, $\bom-\bom_0$ is equivalent to $u=(u_1,u_2)^t$ part. 
%the constraint $|\bom|\sim 1$
Then, the corresponding equation of $u$ becomes
%can be derived from \eqref{E5.1Mar} as 
\begin{align}\label{eq:Aug91}
\rd_t 
\begin{pmatrix}
    u_1\\ u_2
\end{pmatrix}
= 
\boA_\eps 
\begin{pmatrix}
u_1\\ u_2
\end{pmatrix}
+ 
\begin{pmatrix}
u_3u_2'' - u_2u_3''\\
u_1u_3''-u_3u_1''
\end{pmatrix}
+ \text{other nonlinear terms}, 
\end{align}
where
\[
\boA_\eps \coloneqq
\begin{pmatrix}
    0 & -L_0\\
    L_\eps & 0
\end{pmatrix}
,\qquad
L_\eps \coloneqq -\rd_x^2 +1 -2b^2+2\eps b,\qquad 
L_0\coloneqq L_\eps|_{\eps=0}.
\]
\subsubsection*{Local well-posedness (Section \ref{SD})}

We establish local well-posedness by transforming $\bom$ to $\bow$ via \eqref{eq:Sep3}, where the assumption $\bom-\bom_0\in H^s$ translates to $\bow-\hat\boe_3 \in H^s(\R)$ for any integer $s\ge 1$. Although \eqref{E5.1Mar} is a subcritical perturbation of the 1D Schr\"odinger map, proving local well-posedness remains highly nontrivial due to the derivative nonlinearity induced by the DMI term. 

While classical energy methods \cites{MR1877231,MR2304153,MR4482034,MR4099214} yield local well-posedness in $H^s$ for $s\ge 2$, our goal is to establish a blow-up criterion with respect to the $H^1$-norm. For standard Schr\"odinger maps, such criteria are obtained using the Hasimoto transform and Strichartz estimates \cites{MR1737504,MR2547134}. In our setting, however, applying the Hasimoto transform converts the DMI term into a delicate quadratic derivative nonlinearity. Equations of this type are severely problematic; indeed, Kondo and Okamoto \cite{MR5076732} have shown that certain quadratic derivative NLS equations are ill-posed in any Sobolev space. 

In light of this obstacle, we utilize a modified energy method inspired by the work of Hayashi, Ozawa, and Visciglia \cite{MR4881654} on the derivative NLS. Specifically, we construct an explicit $H^2$-energy, $\boE_2^*(\bow)$, satisfying:
\begin{align}\label{eq:Sep4}
\boE_2^*(\bow) \gtrsim \nor{\rd_x^2 \bow}{L^2},\qquad 
|\rd_t \boE_2^*(\bow) | \le C \big(1+P(\nor{\bow-\hat\boe_3}{H^1})\big) \nor{\bow-\hat\boe_3}{H^2}^2,
\end{align}
where $P$ is a polynomial.
\subsubsection*{Linear analysis (Sections \ref{S3},\ref{SA},\ref{SB},\ref{SC})}

%{\color{red}
We investigate the spectral properties of the linearized operator $\boA_\eps$ in \eqref{eq:Aug91} and the corresponding bounds for the propagator $e^{t\boA_\eps}$. Because $\boA_\eps$ structurally resembles linearized operators for NLS standing waves, we employ techniques developed in that context. The main results are summarized in Section \ref{S3}, with proofs detailed in Sections \ref{SA}, \ref{SB}, and \ref{SC}.

Proposition \ref{Tx1} establishes a threshold $\eps_*<0$ given in Theorem \ref{Tm1} such that $\boA_\eps$ is spectrally stable for $\eps>\eps_*$, but admits real eigenvalues $\pm\la_*$ ($\la_*>0$) for $\eps<\eps_*$. The proof in Section \ref{SA} relies on Grillakis--Shatah--Strauss theory \cite{MR1081647} and additional refined techniques from \cites{MR2220744,MR2368894}.

For $\eps>0$, the propagator $e^{t\boA_\eps}$ satisfies dispersive smoothing estimates on the complement of the generalized kernel. To prove this, we shall rule out internal modes, threshold resonances, and embedded eigenvalues (see Proposition \ref{P3.2Jun} and Section \ref{SB}). Excluding internal modes and threshold resonances follows the approach of Krieger and Schlag \cite{MR2219305}. However, 
% we need to guarantee the sign-definiteness of the function which  
% scaler resolvent
% because the spectral structure of $L_\eps$ differs from that in \cite{MR2219305}, 
% \textcolor{blue}{we need to rule out the zeros of the scalar resolvent 
% % supply an additional argument
% (Lemma \ref{Lop2.2}) to verify their framework's assumptions}.
because the spectral structure of $L_\eps$ differs from that in \cite{MR2219305}, we need to proceed differently. 
Indeed, the internal modes and threshold resonances are characterized by zeros of a scaler resolvent function. We rule them out by Lemma \ref{Lop2.2}.
\par
To exclude embedded eigenvalues, we apply a recent method by Collot, Germain, and Pacherie \cite{MR5093257}, which derives a representation formula for the inverses of $L_\eps$ and $L_0$ based on the asymptotic behavior of their eigenfunctions as $|x|\to\infty$. While their original argument relies on a non-trivial $L^2$-kernel of $L_\eps$ with a known sign, our setting lacks such an $L^2$-kernel. Instead, we utilize the Jost solution (a generalized kernel), whose sign is fixed by the positivity of $L_\eps$. This adaptation highlights the approach's robustness.

Consequently, Proposition \ref{P3.2Jun} yields various smoothing estimates, including Strichartz estimates from \cite{MR2219305} and Kato smoothing estimates akin to Cuccagna and Maeda \cites{MR4864876,MR4910476}. A technical distinction arises: whereas \cite{MR4864876} derives resolvent estimates by exploiting perturbative proximity to an integrable problem, our setting is not perturbatively integrable. Therefore, we develop a novel, broadly applicable proof for the corresponding resolvent estimate in Section \ref{SD}.

\subsubsection*{Orbital stability/instability (Sections \ref{S4}, \ref{S5})}

The proof of orbital stability relies on a classical coercivity argument combining energy conservation with modulation \cites{MR783974,MR820338}. While our overall strategy follows C\^ote and Ignat \cite{MR4616656}, a crucial distinction in our setting is the absence of rotational symmetry, which restricts us solely to translational modulation. Decomposing the perturbation as $u={}^t(u_1,u_2)$, an appropriate translational modulation yields coercivity for the $u_2$-component of the Hessian. However, coercivity for the $u_1$-component requires either $\eps>0$, or $\eps>\eps_*$ with an odd $u_1$. This latter constraint fundamentally necessitates the symmetry assumptions in Theorem \ref{Tm1}.

Orbital instability for $\eps<\eps_*$ follows from the presence of an unstable mode of the linearized operator. Specifically, perturbing $\bom_0$ in the unstable direction drives exponential norm-growth. Adapting the framework of Georgiev and Ohta \cite{MR2916078}, we control the nonlinearity via Strichartz-type norms to ensure the linear dynamics dominate. However, applying \cite{MR2916078} directly to \eqref{eq:Aug91} introduces second-order derivative nonlinearities, which would require bounding the $H^2$-norm growth. To circumvent this, we instead apply their framework to the equation obtained via the stereographic projection of $\bow$. This preserves the linear part while eliminating second-order derivative nonlinearities, allowing the estimates to close in the $H^1$-norm, which is inherently controlled by energy conservation. While this semilinearization trick via stereographic projection is well-known for Schr\"odinger maps \cites{MR1737504,MR2278007}, applying it to the \textit{gauge-transformed variable} $\bow$ for stability analysis appears to be a novel strategy.

\subsubsection*{Asymptotic stability (Sections \ref{S6}, \ref{S7})}

Our proof of asymptotic stability under symmetric perturbations follows the bootstrap framework of Cuccagna and Maeda \cite{MR4864876}. We control the high-frequency component and spatial tail via a virial-type estimate \cite{MR3630087}, while bounding the low-frequency, spatially localized component using dispersive smoothing estimates for the linearized equation.

The virial-type estimate is detailed in Proposition \ref{P4:4.1} and proved in Section \ref{S6}. In standard NLS settings, one introduces a virial functional whose time derivative controls localized energy, treating the remainder perturbatively. Here, quasilinearity obstructs this approach: the second-order derivative nonlinearities in \eqref{eq:Aug91} cannot be treated perturbatively. Furthermore, the stereographic projection used for orbital instability fails here, as it destroys the nonlinear structure required for integration by parts against the spatially growing virial multiplier.

To overcome this, we construct a novel virial functional \eqref{eq:Sep9} tailored to \eqref{eq:Aug91}. Inspired by Haldane's momentum functional for Schr\"odinger maps \cite{PhysRevLett.57.1488}, this formulation resolves both structural issues, enabling us to adapt the strategies of \cites{MR4404795,MR4552229}.

Finally, Section \ref{S7} establishes the localized, low-frequency control required to close the bootstrap. To circumvent quasilinearity in this regime, we again utilize equation \eqref{E5.5Apr}, derived via stereographic projection on $\bow$. Following \cite{MR4864876}, we introduce a spatial cutoff in \eqref{E5.5Apr} and apply smoothing estimates for $e^{t\boA_\eps}$. The symmetry assumption in Theorem \ref{Tx1.2Jun} automatically projects the equation onto the continuous spectral subspace of $\boA_\eps$. Notably, unlike in \cite{MR4864876}, bounding the derivative nonlinearity here requires the endpoint Strichartz estimate in conjunction with the Kato smoothing estimate (see \eqref{eq:Sep10}).

\section{Moving frame reduction}
\label{S2}

For a solution $\bom$ to \eqref{E2.1Apr}, we write
\begin{align}\label{eq:Aug71}
\bxi(t,x) = \bom(t,x)-\bom_0(x).
\end{align}
Set
\begin{align}\label{eq:Sep2}
\hat\boe_1 =
\begin{pmatrix}
1 \\ 0 \\ 0
\end{pmatrix}
,\qquad 
\bom_0^\perp =
\begin{pmatrix}
0 \\ -p(x) \\ b(x)
\end{pmatrix}
,
\end{align}
and note that $\{\hat\boe_1,\bom_0^\perp, \bom_0\}$ is an orthonormal basis in $\R^3$ for each $x\in\R$. 
Then we decompose
\begin{align}\label{eq:Aug6}
\bxi = u_1 \hat\boe_1 + u_2 \bom^\perp_0 + u_3 \bom_0.
\end{align}
That is
\[
\bom = u_1 \hat\boe_1 + u_2 \bom_0^\perp + (1+u_3) \bom_0,
\]
with the constraint
\begin{equation}\label{E2.1Apr}
u_1^2 + u_2^2 + (1+u_3)^2 =1.
\end{equation}
Letting $u={}^t(u_1,u_2)$, we can write 
$u_3 = -\frac{|u|^2}{1+\sqrt{1-|u|^2}}$, and thus $u_3$ is of higher order then $|u|$. 
% $u_3$ is of higher order than $|u|$ because}
% $
% u_3 = -\frac{|u|^2}{1+\sqrt{1-|u|^2}}.
% $
Therefore, as long as $|u|\le \frac 12$, we have
\begin{align}\label{eq:Aug5}
u_3 = O(|u|^2),&& \rd_x u_3 = O(|u||\rd_x u|).
\end{align}
Throughout this paper, we frequently reformulate the problem in terms of the new variable $u={}^t(u_1,u_2)$. The governing equation for $u$ is obtained by substituting \eqref{eq:Aug6} into \eqref{E5.1Mar}. 
First, we compute the derivatives as
\begin{align*}
&\rd_x \bom_0 = b \bom_0^\perp,\qquad \rd_x \bom_0^\perp = -b \bom_0,\\
& \rd_x \bom = \rd_x u_1 \hat\boe_1 + (\rd_x u_2 + b(1+u_3)) \bom_0^\perp 
+ (\rd_x u_3 - bu_2) \bom_0,\\
& \rd_{xx} \bom= u_1''\hat\boe_1 +
(u_2'' - bp (1+u_3) + 2bu_3' - b^2u_2) \bom_0^\perp \\
&\hspace{60pt}
+(u_3'' + bpu_2 -2bu_2' -b^2 (1+u_3)) \bom_0.
\end{align*}
%Hence we have
Also, since
$%\begin{align*}
\rd_t \bom
=
\rd_t u_1 \boe_1 + \rd_t u_2 \bom_0^\perp + \rd_t u_3 \bom_0,
$ %\end{align*}
then we have,
\begin{align*}
\bom \times \bom''
&= (u_1 \hat\boe_1 + u_2 \bom_0^\perp + (1+u_3)\bom_0) \\
&\qquad \times 
[u_1''\hat\boe_1 +
(u_2'' - bp (1+u_3) + 2bu_3' - b^2u_2) \bom_0^\perp \\
&\hspace{60pt}
+(u_3'' + bpu_2 -2bu_2' -b^2 (1+u_3)) \bom_0
]\\
&= 
[(1+u_3)u_2'' -bp(1+u_3)^2 + 2b(1+u_3)u_3' - u_2u_3'' -bpu_2^2 +2bu_2u_2']\hat\boe_1 \\
&\qquad
+[u_1u_3''+bpu_1u_2-2bu_1u_2'-b^2u_1(1+u_3) -(1+u_3)u_1''] \bom_0^\perp\\
&\qquad 
+[u_2u_1''-u_1u_2''+bpu_1(1+u_3)-2bu_1u_3'+b^2u_1u_2] \bom_0.
\end{align*}
Next,
\begin{align*}
m_3 &= (u_1\hat\boe_1 +u_2\bom_0^\perp +(1+u_3)\bom_0) \cdot \hat\boe_3 =
bu_2 + p(1+u_3),\\
\bom\times \hat\boe_3
&= 
(-pu_2 + b(1+u_3)) \hat\boe_1 - u_1 \hat\boe_2\\
&= (-pu_2 + b(1+u_3)) \hat\boe_1 
+ pu_1 \bom_0^\perp -bu_1 \bom_0,
\end{align*}
and hence
\begin{align*}
m_3\bom\times \hat\boe_3 &= 
(bu_2 + p(1+u_3))
[(-pu_2 + b(1+u_3)) \hat\boe_1 
+ pu_1 \bom_0^\perp -bu_1 \bom_0]\\
&=
[-bpu_2^2 +(b^2-p^2)u_2(1+u_3) +bp(1+u_3)^2] \hat\boe_1\\
&\qquad +[bpu_1u_2 +p^2u_1(1+u_3)] \bom_0^\perp +
[-b^2u_1u_2 -bpu_1(1+u_3)]\bom_0.
\end{align*}
Finally,
\begin{align*}
\hat\boe_1\times \bom' 
&= \hat\boe_1 \times 
(u_1' \hat\boe_1 + (u_2'+b(1+u_3))\bom_0^\perp +(u_3'-bu_2)\bom_0)\\
&= (u_3'-bu_2) \bom_0^\perp 
-(u_2'+b(1+u_3)) \bom_0
\end{align*}
and hence
\begin{align*}
\bom\times (\hat\boe_1\times \bom')
&= 
\left((1+u_3)u_3' -bu_2(1+u_3) +u_2u_2' +bu_2(1+u_3)\right) \hat\boe_1 \\
&\qquad 
+(-u_1u_2'-bu_1(1+u_3)) \bom_0^\perp + (-u_1u_3' + bu_1u_2) \bom_0.
\end{align*}
Combining these, we obtain
%
%
%\begin{screen}
\begin{align}\label{E2.2Apr}
\begin{aligned}
\rd_t u_1 
&=
-L_0 u_2 + u_3 u_2'' - u_2 u_3''
-2(b-\eps) u_1u_1'
-2bpu_2^2 + (2b^2-1)u_2u_3 \\
\rd_t u_2 
&=
L_\eps u_1 + u_1u_3'' -u_3u_1'' 
- 2(b-\eps)u_1u_2' +2bp u_1u_2  +(1-2b^2+2\eps b) u_1u_3 \\
\rd_t u_3 
&= u_2 u_1'' - u_1 u_2'' - 2(b-\eps) u_1 u_3' - 2\eps bu_1u_2,
\end{aligned}
\end{align}
where 
\begin{equation}\label{eq:Aug7}
L_\eps \coloneqq -\rd_x^2 +1 -2b^2 +2\eps b.
\end{equation}
%\end{screen}
Alternatively, 
we may reorganize \eqref{E2.2Apr} into
%\begin{screen}
\begin{gather}
\rd_t \bm{w} =  \bm{w}\times \bbH \bm{w} -2(b-\eps) w_1\bm{w}' -2bpw_2\, \hat\boe_3\times \bm{w},\label{E2.4Aug}\\
\bm{w} = \begin{pmatrix}
    u_1\\u_2\\1+u_3
\end{pmatrix}
,\qquad \bbH = \diag (L_\eps , L_0,-\rd_{xx}).\notag
\end{gather}
%\end{screen}
%
%
%

\section{Linear properties}
\label{S3}
Since the constraint \eqref{E2.1Apr} determines $u_3$ in terms of $u_1$ and $u_2$, \eqref{E5.1Mar} reduces to the following system for $(u_1,u_2)$:
\begin{align}\label{eq:Aug81}
\rd_t 
\begin{pmatrix}
    u_1 \\u_2
\end{pmatrix}
=
\boA_\eps 
\begin{pmatrix}
    u_1 \\u_2
\end{pmatrix}
+ \text{h.o.t},\qquad
\boA_\eps 
\coloneqq
\begin{pmatrix}
0 & -L_0\\
L_\eps & 0
\end{pmatrix}.
\end{align}
Crucially, the structure of $\boA_\eps$ parallels that of linearized operators around NLS standing waves \cites{MR2368894,MR2219305}. This structural similarity allows us to derive spectral properties and linear estimates by adapting established NLS frameworks.
\subsection{Spectral properties}

We first characterize the spectral stability of $\boA_\eps$ by identifying a threshold parameter.

\begin{prop}\label{Tx1}
There exists $\eps_* <0$ such that
\begin{itemize}
\item If $\eps \ge \eps_*$, then $\si (\boA_\eps) \subset i\R$. 
\item If $\eps < \eps_*$, then $\si (\boA_\eps) \not\subset i\R$. 
More precisely, $\boA_\eps$ has a symmetric pair of real eigenvalues $\pm\la_*$, with $\la_*>0$ (arising from the Hamiltonian symmetry). Moreover, the corresponding eigenfunctions are smooth and exponentially decaying as $|x|\to\infty$.
% In detail, there exists a pair of point spectrum $\pm \la_*$ with $\la_* >0$. Moreover, the corresponding eigenfunction is smooth, exponentially decaying as $|x|\to\infty$.
\end{itemize}
\end{prop}

% \textcolor{teal}{I'd prefer....
%  More precisely, $\boA_\eps$ has a symmetric pair of real eigenvalues $\pm\la_*$, with $\la_*>0$ (arising from the Hamiltonian symmetry). Moreover, the corresponding eigenfunctions are smooth and exponentially decaying as $|x|\to\infty$.
% }

%\end{screen}
We postpone the proof to Section \ref{SA}. 
To establish smoothing estimates for $e^{t\boA_\eps}$, we need further information on the spectrum of $\boA_\eps$. In particular, we establish the following spectral properties:

\begin{prop}\label{P3.2Jun}
For $\eps>0$, the following hold:
\begin{itemize}
    \item $\si(\boA_\eps) = \{0\} \cup \pm i[1,\infty)$.
    \item $\si_p(\boA_\eps)= \{0\}$, and the generalized kernel $\boN_g$ of $\boA_\eps$ is spanned by 
    $\begin{bmatrix}
    0 \\ b
    \end{bmatrix}$ and     $\begin{bmatrix}
    L_\eps^{-1} b \\ 0
    \end{bmatrix}$.
    \item $\boA_\eps$ does not have threshold resonance nor embedded eigenvalues; 
    see the appendix for the definitions of these notions.
\end{itemize}   
\end{prop}
%\end{itembox}
%Since this part is technical, 
%To keep the momentum, 
%\dlt{We postpone the proof into} 
%\textcolor{blue}{
The proof is deferred to Section \ref{SB}.
%}  %Proposition \ref{P3.2Jun} and its proof into the appendix. 
% The thing here is that once Proposition \ref{P3.2Jun} has verified, we may apply the known smoothing estimates.
%
%
\subsection{Smoothing estimates for $e^{t\boA_\eps}$}

% \dlt{Thanks to} \textcolor{blue}{By virtue of} Proposition \ref{P3.2Jun}, as well as the property of potential (\textcolor{blue}{namely, smoothness, exponential decay, and evenness} \dlt{smooth, exponential decay, even}),
% %the exponential decay of the potential, 
% a series of known smoothing estimates for $e^{t\boA_\eps}$ \dlt{are} \textcolor{blue}{is} available. 
% Let us set
% \[
% P_d \coloneqq \frac{1}{2\pi i} \oint_\gam (zI -\boA_\eps)^{-1} dz,\qquad P_c\coloneqq  I-P_d,
% \]
% where $\gam$ is a contour enclosing $z=0$.
% Note that $P_d$ corresponds to the projection of $L^2(\R)^2$ onto $\boN_g$.

%\textcolor{teal}{
Given the smoothness, exponential decay, and even parity of the potential, Proposition \ref{P3.2Jun} allows us to invoke several established smoothing estimates for $e^{t\boA_\eps}$. We define the Riesz projection $P_d$ onto the generalized kernel $\boN_g$, and the complementary spectral projection $P_c$, as
\[
P_d \coloneqq \frac{1}{2\pi i} \oint_\gam (zI-\boA_\eps)^{-1}\,dz,
\qquad
P_c \coloneqq I-P_d,
\]
where $\gam$ is a contour enclosing $z=0$. 

The first result is the Kato smoothing estimate established in \cite{MR4864876}. To formulate it, we introduce the weighted $L^2$-norm:
\[
\nor{f}{L^{2,s}(\R)} \coloneqq \nor{\jb{x}^sf}{L^2(\R)},\qquad s\in\R.
\]

\begin{prop}[Homogeneous Kato smoothing]\label{P4:3.3}
If $\eps>0$, then for $s>1$, we have
\[
\nor{P_ce^{t\boA_\eps}f}{L^2_tL^{2,-s}_x} \le C \nor{f}{L^{2}}.
\]
\end{prop}
\begin{prop}[Inhomogeneous smoothing estimate]\label{P4:3.2}
If $\eps>0$, then for $s>\frac 32$ and $\tau>\frac 12$, we have
\[
\nor{\int_0^t P_c e^{(t-t')\boA_\eps} g(t') dt'}{L^2_tL^{2,-s}_x} \le C \nor{g}{L^2_tL^{2,\tau}_x}.
\]
\end{prop}
While analogous estimates for NLS standing waves appear in \cite{MR4864876}*{Propositions 2.3 and 2.4}, their derivation relies on the proximity of the underlying model to an integrable system (specifically for the resolvent estimate in \cite{MR4864876}*{Proposition 8.6}). To bypass this restriction, Section \ref{SC} provides an alternative, integrability-free proof of the corresponding resolvent estimate, making it applicable to a broader class of problems.

The second result is the Strichartz estimate for $e^{t\boA_\eps}$:
\begin{prop}[Strichartz estimate \cite{MR2219305}*{Proposition 7.1}]
\label{P4:3.4}
\begin{gather*}
\nor{\int_0^t P_c e^{(t-t')\boA_\eps}g(t',x)dt'}{L^4_tL^\infty_x} \le C \nor{g}{L^{\frac 43}_t L^1_x},   
\\
\nor{\int_{-\infty}^\infty P_c e^{-t'\boA_\eps}g(t',x)dt'}{L^2_x} \le C \nor{g}{L^{\frac 43}_t L^1_x}.
\end{gather*}
\end{prop}
The proof is given in \cite{MR2219305}*{Proposition 7.1}; by Proposition \ref{P3.2Jun}, $\boA_\eps$ satisfies all required hypotheses of \cite{MR2219305}. For the subsequent nonlinear analysis, we combine the Kato smoothing and Strichartz estimates as follows:
%}

\begin{prop}\label{P4:3.5}
For $s>\frac 32$, we have
\begin{align}\label{E3.1Apr}
\nor{\int_0^t P_c e^{(t-t')\boA_\eps}g(t',x)dt'}{L^2_tL^{2,-s}_x} \le C \nor{g}{L^{\frac 43}_t L^1_x}.  
\end{align}
\end{prop}
\begin{proof}
By Propositions \ref{P4:3.3} and \ref{P4:3.4}, we obtain
\begin{align*}
\nor{\int_{-\infty}^\infty P_c e^{(t-t')\boA_\eps}g(t',x)dt'}{L^2_tL^{2,-s}_x}
&=
\nor{e^{t \boA_\eps}P_c\int_{-\infty}^\infty P_c e^{-t'\boA_\eps}g(t',x)dt'}{L^2_tL^{2,-s}_x}\\
&\le C\nor{\int_{-\infty}^\infty P_c e^{-t'\boA_\eps}g(t',x)dt'}{L^2_x}\le C\nor{g}{L^{\frac 43}_tL^1_x}.
\end{align*}
By the Christ--Kiselev lemma \cite{MR1809116}, we conclude the retarded version as desired.
%claimed in \eqref{E3.1Apr}. 
\end{proof}

\section{Orbital stability}
\label{S4}
%
%\subsection{Stability}
%
The main claim of this section is the following:
\begin{thm}\label{T2.1}
The following hold.\\
(i) If $\eps>0$, then $\bom_0$ is orbitally stable in $\boM$.\\ 
(ii) If $\eps>\eps_*$, then $\bom_0$ is stable in $\boM_\sym$. 
\end{thm}
% \begin{thm}[Orbital stability]\label{T2.1}
% The following hold.\\
% (i) 
% %(Orbital stability in $\boM$ for $\eps>0$.)
% If $\eps>0$, then $\bom_0$ is orbitally stable in $\boM$. More precisely,  %In more detail, the following holds:
% there exist $\del_1>0$, $C>0$ such that if $\del \coloneqq \inf_{y\in\R} \nor{\bom(0,\cdot) - \bom_0 (\cdot-y)}{H^1} < \del_1$, then the solution $\bom\in C_t\boM$ is global in time, and satisfies
% \begin{align}\label{eq:Aug74}
% \inf_{y\in\R} \nor{\bom(t,\cdot) - \bom_0 (\cdot-y)}{H^1} \le C\del,\qquad t\in \R.%(T_{\max}^-,T_{\max}^+).
% \end{align}
% %{\color{red}Remark: We are not able to verify global solvability without symmetry. See the argument below.}\\
% (ii) If $\eps>\eps_*$, then $\bom_0$ is orbitally stable in $\boM_\sym$. More precisely, there exist $\del_1>0$, $C>0$ such that if $\del \coloneqq \nor{\bom(0,\cdot) - \bom_0}{H^1} < \del_1$, then 
% the solution $\bom\in C_t\boM_{\sym}$ is global in time, and satisfies
% \begin{align}\label{eq:Aug75} 
% \nor{\bom(t,\cdot) - \bom_0}{H^1} \le C\del,\qquad t\in\R.
% \end{align}
% \end{thm}
%
%
%{\color{blue}
\begin{remark}
% We provide some remarks on the difference between (i) and (ii). 
% First, in (i), it remains open whether the solution globally exists. In fact, the blow-up criterion in Proposition \ref{P:wp2} is only characterized by $\nor{\bom-\bom_0}{H^1}$, and thus even by \eqref{eq:Aug74}, it is a priori possible the translational modulation $y(t)$ would blow up in finite time. Second, in the case 
When $\eps \in [\eps_*,0]$, the proof of (i) fails due to the lack of coercivity of the $u_1$-part in \eqref{eq:Aug77}, see Lemma \ref{C2.1}. 
%If one would be able to make $u_1$ orthogonal to $b$, which is the case when $\bom\in\boM_{\sym}$, then ....
%In general, however, 
Since \eqref{E5.1Mar} has no rotational symmetry, there seems to be no possible modulation framework to make $u_1$ orthogonal to $b$.  
\end{remark}
%}
%
%
\subsection{Energy expansion}\label{S2.3}
The proof is based on %\dlt{the} 
a classical argument using translational modulation and the coercivity of $\nab^2E$ (see \cites{MR783974,MR820338}). 
Given $y\in \R$, we write
\begin{align}\label{eq:Aug73}
\bom =\tau_y (\bom_0+\bxi),
\end{align}
instead of \eqref{eq:Aug71}, where $\tau_y f= f(\cdot -y)$. 
We begin by deriving the formula of the energy Hessian in terms of $\bxi$. Let $e_\eps(\bom)$ be the energy density, namely $E_\eps (\bom) = \int_{\R}e_\eps(\bom)dx$. Then, we have
\begin{align*}
e_\eps (\bom_0+\bxi) - e_\eps (\bom_0)
&=(-\rd_x^2 \bom_0 - m_{0,3} \hat\boe_3 +2\eps \hat\boe_1\times \rd_x \bom_0) \cdot \bxi\\
&\qquad + \frac 12 |\rd_x \bxi|^2 - \frac 12 \xi_3^2 + \eps \rd_x \bxi \times \bxi \cdot \boe_1.
\end{align*}
% %
% \begin{align*}
% E_\eps (\bom_0 +\bxi)
% &= 
% \int_\R \frac 12 |\rd_x \bom_0 + \rd_x \bxi|^2 
% + \frac 12 \left( 1- (m_{0,3} + \xi_3)^2 \right) \\
% &\hspace{40pt} +\eps \left(
% \rd_x (\bom_0 +\bxi) \times  (\bom_0 +\bxi)
% \right) \cdot \hat\boe_1 dx\\
% &=
% E_\eps (\bom_0) 
% + \int_\R (-\rd_x^2 \bom_0 - m_{0,3} \hat\boe_3 +2\eps \hat\boe_1\times \rd_x \bom_0) \cdot \bxi \\
% &\hspace{40pt} 
% +\int_\R \frac 12 |\rd_x \bxi|^2 - \frac 12 \xi_3^2 + \eps \rd_x \bxi \times \bxi \cdot \boe_1 dx
% \end{align*}
Here, $\bom_0$ satisfies
\[
\rd_x^2 \bom_0 + m_{0,3} \hat\boe_3 
= (1-2b^2) \bom_0,\qquad 
\hat\boe_1 \times \rd_x \bom_0 = -b \bom_0.
\]
Thus
\begin{align*}
(-\rd_x^2 \bom_0 - m_{0,3} \hat\boe_3 +2\eps \hat\boe_1\times \rd_x \bom_0) \cdot \bxi 
&=
\left( -1+2b^2  -2\eps b\right) \bom_0 \cdot \bxi \\
&=
\frac 12 \left(1-2b^2  +2\eps b\right) |\bxi|^2,
\end{align*}
where we used $\bom_0 \cdot \bxi = -\frac 12 |\bxi|^2$ given by the constraint. 
% \[
% 1 =|\bom|^2 =
% |\bom_0 +\bxi|^2 \quad \Longrightarrow \quad 
% \bom_0 \cdot \bxi = -\frac 12 |\bxi|^2.
% \]
Hence we obtain
\begin{align}\label{eq:Sep22}
\begin{aligned}
H(\bxi) &\coloneqq 2\int_{\R} (e_\eps (\bom_0+\bxi) - e_\eps (\bom_0)) dt\\
&= \int_\R |\rd_x \bxi|^2 - \xi_3^2 + 2\eps \rd_x\bxi \times \bxi \cdot \hat\boe_1 
+ \left(1-2b^2  +2\eps b\right) |\bxi|^2 dx.
\end{aligned}
\end{align}
We remark that although a map $\bom\in\boM$ does not define its energy due to the lack of integrability, $H(\bxi)$ is well-defined as $\bxi\in H^1$.
% \dlt{We remark that although the energy of $\bom_0$ or $\bom_0+\bxi$ is not defined due to the non-integrability, its difference $H(\bxi)$ is well-defined provided $\bxi\in H^1$.}
% \begin{screen}
% \begin{gather*}
% E_\eps (\bom_0 +\bxi)
% = 
% E_\eps (\bom_0) + \frac 12 H(\bxi),
% \\
% H(\bxi) := 
% \int_\R |\rd_x \bxi|^2 - \xi_3^2 + 2\eps \rd_x\bxi \times \bxi \cdot \hat\boe_1 
% + \left(1-2b^2  +2\eps b\right) |\bxi|^2 dx.
% \end{gather*}
% \end{screen}

\subsection{Coercivity of $H$}\label{S2.4}
Now we apply the moving frame reduction \eqref{eq:Aug6} for $\bxi$. Note that \eqref{E2.1Apr} and \eqref{eq:Aug5} still hold regardless of $y$. 
% \dlt{Here we shall rewrite $H(\bxi)$ with respect to $u$.} 
% Let
% \[
% \bom_0^\perp (x) = 
% \begin{pmatrix}
% 0 \\ -\tanh x \\ \sech x
% \end{pmatrix}
% .
% \]
% Then $\hat\boe_1, \bom_0^\perp, \bom_0$ forms an orthonormal frame of $\R^3$. 
% Using this, we write
% \[
% \bxi = u_1 \hat\boe_1 + u_2 \bom_0^\perp + u_3 \bom_0.
% \]
% Note that if $|\bxi|\ll 1$, then $\bxi$ and $u=(u_1,u_2)$ is one to one, as
% \begin{align*}
% |\bom_0 + \bxi|^2 =1 \quad \Longrightarrow \quad 
% |u|^2 + (u_3+1)^2 =1 
% \quad 
% &\Longrightarrow \quad  u_3 = -1 + \sqrt{1-|u|^2}\\
% &\Longrightarrow \quad u_3 = -\frac{|u|^2}{1+ \sqrt{1-|u|^2}}.
% \end{align*}
Since 
\[
\rd_x \bxi = \rd_x u_1 \hat\boe_1 + (\rd_x u_2 + bu_3) \bom_0^\perp 
+ (\rd_x u_3 - bu_2) \bom_0,
\]
we compute each term in the integrand of $H(\bxi)$ as
\begin{align*}
|\rd_x \bxi|^2 &= |\nab u|^2 + b^2 u_2^2 + O(|u|^2|\rd_xu|+ |u||\rd_xu|^2 + |u|^3),\\
\hat\boe_1 \times \rd_x \bxi \cdot \bxi &= -bu_2^2 +u_2 \rd_x u_3 - u_3 \rd_x u_2 - bu_3^2\\
&=
-bu_2^2 + O(|u|^2|\rd_xu| + |u|^3),\\
|\bxi|^2 &= |u|^2 + O(|u|^3).
\end{align*}
% Now we compute as
% \begin{align*}
% & \rd_x \bom_0 = b \bom_0^\perp,\qquad \rd_x \bom_0^\perp = -b \bom_0,\\
% & \rd_x \bxi = \rd_x u_1 \hat\boe_1 + (\rd_x u_2 + bu_3) \bom_0^\perp 
% + (\rd_x u_3 - bu_2) \bom_0, \\
% & |\rd_x \bxi|^2 = |\nab u|^2 + b^2 u_2^2 + O(|u|^2|\rd_xu|+ |u||\rd_xu|^2 + |u|^3), \\
% & \hat\boe_1 \times \rd_x \bxi \cdot \bxi = -bu_2^2 + O(|u|^2|\rd_xu| + |u|^3), \\
% &  |\bxi|^2 = |u|^2 + O(|u|^3).
% \end{align*}
Therefore,
\begin{align}\label{eq:Aug77}
H(\bxi) 
% &= 
% \int_\R (-\rd_x^2 u_1 + (1-2b^2 +2\eps b) u_1) \cdot u_1 + 
% (-\rd_x^2 u_2 + (1-2b^2 ) u_2) \cdot u_2 dx \\
&= \inp{L_\eps u_1}{u_1}_{L^2} + \inp{L_0 u_2}{u_2}_{L^2} + O(\nor{u}{H^1}^3),
\end{align}
where $L_\eps$ is defined in \eqref{eq:Aug7}. 
We now examine the coercivity of $H(\bxi)$.
% \dlt{Here, let us observe the coercivity of $H(\bxi)$.}
%
%
\begin{lem}\label{L2.1}
%If $\mu_\eps >0$,
(i) If $\eps>0$, then there exists $C=C_\eps>0$ such that
\[
\inp{L_\eps f}{f} \ge C \nor{f}{H^1(\R)}^2,\qquad f\in H^1.
\]
(ii) If $\eps>\eps_*$, then there exists $C=C_\eps>0$ such that
\[
\inp{L_\eps f}{f} \ge C \nor{f}{H^1(\R)}^2,\qquad \text{for } f\in H^1 \text{ with } f\perp b \text{ in } L^2.
\]
%In particular, this is true 
%if $\eps > \eps_*$.
\end{lem}

\begin{proof}
The proof follows by a classical contradiction argument using the spectral gap of $L_\eps$ and Proposition \ref{P1.3}. For the details, see \cite{MR4899905}*{Lemma 2.5} for example.

%%%
%%% Proof of Lemma 4.1
%%%
%\begin{myleftbar}
% We only give a proof of (ii); The proof of (i) can be conducted similarly. 
% By Proposition \ref{P1.3}, if $\eps>\eps_*$ it holds that
% \[
% \inf_{\substack{f\perp b \text{ in } L^2 \\ f\neq 0}} 
% \inp{L_\eps f}{f}_{L^2} > 0.
% \]
% If we suppose the contrary to Lemma \ref{L2.1}, then we have
% \[
% \inf_{\substack{f\perp b \text{ in } L^2 \\ f\neq 0}} 
% \frac{\inp{L_\eps f}{f}_{L^2} }{\nor{f}{H^1}^2} =0. 
% \]
% Then there exists $\{f_n \}$ with $f_n\perp b$ such that
% \[
% \nor{f_n}{H^1} =1,\qquad \inp{L_\eps f_n}{f_n} \xrightarrow{n\to\infty} 0.
% \]
% By weak compactness, we have $f_n \wto {}^\exists f_*$ in $H^1(\R)$ with $f_*\perp b$. 
% Using the Rellich-Kondrashov theorem and the decay of the potential, we have
% \[
% \inp{(-2b^2+2\eps b) f_n}{f_n} \xrightarrow{n\to\infty} 
% \inp{(-2b^2+2\eps b) f_*}{f_*}.
% \]
% Then
% \begin{align*}
% 1= \nor{f_n}{H^1}^2 &= 
% \inp{L_\eps f_n}{f_n} - \inp{(-2b^2+2\eps b) f_n}{f_n}\\
% &\xrightarrow{n\to\infty} \inp{(2b^2-2\eps b) f_*}{f_*},
% \end{align*}
% which especially implies $f_*\neq 0$. Hence by the assumption $\mu_\eps>0$, 
% \begin{align*}
% 0 < \inp{L_\eps f_*}{f_*} &= \nor{f_*}{H^1}^2 + \inp{(-2b^2+2\eps b) f_*}{f_*} \\
% &\le \liminf_{n\to\infty} \nor{f_n}{H^1}^2 + \lim_{n\to\infty} \inp{(-2b^2+2\eps b) f_n}{f_n} \\
% &= \liminf_{n\to\infty} \inp{L_\eps f_n}{f_n} = 0
% \end{align*}
% which is a contradiction. 
% % \dlt{The last assertion follows from \dlt{by} Proposition \ref{Tx1}.}
% %\end{myleftbar}
%%%
%%% Proof of Lemma 4.1 end
%%%
\end{proof}

\begin{lem}\label{C2.1}
There exist $C>0$, $\del_{11}>0$ such that if $\eps>0$ and $\nor{\bxi}{H^1}<\del_{11}$, then 
\[
H(\bxi) \ge C\nor{\bxi}{H^1}^2,\qquad \text{if } u_2\perp b \text{ in } L^2.
\]
\end{lem}
\begin{proof}
If $\nor{\bxi}{H^1}<\del_{11} \ll1$, then \eqref{eq:Aug5} and the Sobolev embedding imply $\nor{u}{H^1} \sim \nor{\bxi}{H^1} \ll 1$. Therefore, Lemma \ref{L2.1} yields
\begin{align*}
H(\bxi) &= \inp{L_\eps u_1}{u_1}_{L^2} + \inp{L_0 u_2}{u_2}_{L^2} + O(\nor{u}{H^1}^3) \\
&\gtrsim \nor{u}{H^1}^2 + O(\nor{u}{H^1}^3) \gtrsim \nor{u}{H^1}^2\sim \nor{\bxi}{H^1}^2.
\end{align*}
\end{proof}

\subsection{Modulation}
%
%
%{\color{blue}
Now we show the existence of modulation which achieves the orthogonality as in Lemma \ref{C2.1}. Indeed, it is realized by the minimizer of $y\mapsto \nor{\bom-\tau_y\bom_0}{L^2}$. 
More precisely, we have the following:
% \[
% d_{L^2} (\bom) \coloneqq \inf_{y\in\R} \nor{\bom-\tau_y \bom_0}{L^2}.
% \]
%
\begin{lem}\label{L2.2}
There exists $\del_{12}>0$ such that if $\inf_{y\in\R} \nor{\bom-\tau_y\bom_0}{L^2}<\del_{12}$, then there exists unique $y=y(\bom)$ which attains the minimum of \eqref{eq:Aug105}. Moreover, the mapping $\bom\mapsto y(\bom)$ is $C^1$. Furthermore, if we write $\bom$ as in \eqref{eq:Aug73} with $y=y(\bom)$, we have $\inp{u_2}{b}_{L^2}=0$.
\end{lem}
\begin{proof}
Since
\[
\limsup_{|y|\to+\infty} \nor{\bom-\tau_y\bom_0}{L^2}^2 \ge  \limsup_{|y|\to+\infty} \nor{b-\tau_yb}{L^2}^2 = 2 \nor{b}{L^2}^2, 
\]
the existence of the minimizer is guaranteed provided $\del_{12}\ll1$. 
The uniqueness follows by exploiting the convexity of the map $y\mapsto \nor{\bom-\tau_y \bom_0}{L^2}$, which proceeds along lines similar to \cite{MR2290708}*{Theorem 2.1, Lemma 4.6}, and thus we leave the details to the reader. If $y$ is the minimizer, then $y$ should satisfy
\begin{align}\label{eq:Aug107}
F(\bom,y)\coloneqq \int_\R (\bom-\tau_y\bom_0)\cdot \tau_y\bom_0' dx =0\quad \implies\quad \inp{u_2}{b}_{L^2}=0.
\end{align}
We note that $F$ is $C^1$, and satisfies 
%\[
$F( \bom_0 ,0) = 0$, 
$\frac{\rd F}{\rd y} ( \bom_0 ,0) = \nor{b}{L^2}^2 >0$.
%\]
Therefore, the implicit function theorem implies that the mapping $y\mapsto y(\bom)$ is also $C^1$. Hence the proof is complete.
\end{proof}

\subsection{Proof of Theorem \ref{T2.1}}
%{\color{blue}
To simplify the notation, we denote
\begin{align}\label{eq:Aug105}
d (\bom) := \inf_{y\in \R} \nor{\bom - \tau_y \bom_0}{H^1}.
\end{align}
By the same argument as in Lemma \ref{L2.2}, \eqref{eq:Aug105} also attains its minimum at a unique point $y_*=y_*(\bom)$. %We observe that %$\nor{\bom-\tau_y\bom_0}{L^2}\le $
Letting $y=y(\bom)$, we have
\begin{align*}
\nor{\bom-\tau_{y}\bom_0}{L^2} = \inf_{y'\in\R} \nor{\bom-\tau_{y'}\bom_0}{L^2} \le \nor{\bom-\tau_{y^*}\bom_0}{L^2} \le d(\bom),
\end{align*}
and thus
\begin{align}\label{eq:Sep21}
\begin{aligned}
|y(\bom)-y_*(\bom)|\lesssim
\nor{\tau_y\bom_0 -\tau_{y_*}\bom_0}{L^2} 
&\le \nor{\bom-\tau_y\bom_0}{L^2} + 
\nor{\bom-\tau_{y^*}\bom_0}{L^2}\\
&\lesssim d(\bom).
\end{aligned}
\end{align}
To prove Theorem \ref{T2.1}, 
we only consider the forward time direction. 
% , while the backward case follows similarly. 
Let $\bom\in C((T_{\max}^-,T_{\max}^+)\colon \boM)$ be a solution as in Proposition \ref{P:wp1}. 
Let $0<\del_1<\del_{13} 
%<\min\{\del_1,\del_2\}
$ be numbers determined later, and 
suppose that $d(\bom(0)) <\del_1$. 
Define
\[
T_1 := \sup \left\{
T\in (0,T_{\max}^+) \,\left|\, \sup_{t\in [0,T)} d(\bom(t)) < \del_{13}
\right.
\right\}. 
\]
We claim that $T_1=T_{\max}^+$ by choosing $\del_1,\del_{13}$ sufficiently small. 
Suppose $\del_{13}<\del_{12}$, then for $t\in [0,T_1)$, 
Lemma \ref{L2.2} defines a function $y(t)=y(\bom(t))$. 
% which minimizes \eqref{eq:Aug105} for each $t$
Set $\bxi(t),u_1(t),u_2(t)$ as in \eqref{eq:Aug73} and \eqref{eq:Aug6}. 
By virtue of \eqref{eq:Sep21}, we have
\[
\nor{\bxi}{H^1} = \nor{\bom-\tau_y\bom_0}{H^1} \le d(\bom(t)) + C|y-y_*|
\lesssim \del_{13}.
\]
By choosing $\del_{13}$ sufficiently small, it follows from Lemma \ref{C2.1} that
\[
d(\bom(t))^2 \le 
\nor{\bxi(t)}{H^1}^2 \lesssim H(\bxi(t)) =  H(\bxi)|_{t=0} \lesssim \del_1^2.
\]
Hence, by choosing $\del_1$ sufficiently small, we have $d(\bom(t))<\del_{13}$ for all $t\in [0,T_1)$, and hence $T_1=T_{\max}^+$ follows.\par
%}
% {\color{gray}
% We only discuss forward in time, while the backward case follows similar. 
% Let $\bom\in C((T_{\max}^-,T_{\max}^+)\colon \boM)$ be a solution as in Proposition \ref{P:wp1}. 
% Let $0<\del_0<\del_3 <\min\{\del_1,\del_2\}$ be a number determined later, and 
% suppose that $d(\bom(0)) <\del_0$. 
% Define
% \[
% T_1 := \sup \left\{
% T\in (0,T_{\max}) \,\left|\, \sup_{t\in [0,T)} d(\bom(t)) < \del_3
% \right.
% \right\}. 
% \]
% We claim that if $\del_0,\del_3$ is sufficiently small, then $T_1=T_{\max}$. 
% Since $\del_3<\del_2$, for $t\in [0,T_1)$, 
% Lemma \ref{L2.2} defines a function $y(t)=y(\bom(t))$ which minimizes \eqref{eq:Aug105} for each $t$. Let $\bxi(t),u_1(t),u_2(t)$ be defined as in \eqref{eq:Aug73} and \eqref{eq:Aug6}. Then we have
% \[
% \nor{\bxi}{H^1} = \nor{\bom-\tau_y\bom_0}{H^1} < \del_3.
% \]
% Since $\del_3<\del_1$, it follows from Lemma \ref{C2.1} that
% \[
% d(\bom(t))^2 =\nor{\bxi(t)}{H^1}^2 \le CH(\bxi(t)) = C H(\bxi)|_{t=0} \le C \del_0^2.
% \]
% Hence, by choosing $\del_0$ sufficiently small, we have $d(\bom(t))<\del_3$ for all $t\in (T_{\max}^-,T_{\max}^+)$, and thus $T_1=T_{\max}$.}\par 
%{\color{blue}
Next, let us show $T_{\max}^+=\infty$. Since $\bom(t)$ remains bounded, the only possible scenario to be ruled out is the blow-up of the translational modulation parameter. 
Letting $y(t)=y(\bom(t))$ as in Lemma \ref{L2.2}, we claim that
\begin{align}\label{eq:Aug103}
|\dot{y}(t)|\le C \nor{\bom(t)-\bom_0}{H^1_x},\qquad t\in (0,T_{\max}^+),
\end{align}
where $C>0$ is independent of $\bom$ and $t$. Once \eqref{eq:Aug103} is shown, we have
\begin{align*}%\label{eq:Aug104}
|\dot{y}(t)| \le C \nor{\bom (t)-\bom_0}{H^1} 
&\le C \nor{\bom(t)-\tau_{y}\bom_0}{H^1} + C\nor{\tau_y\bom_0-\bom_0}{H^1}\\
&\le C\del_{1} + C|y(t)|
\end{align*}
by \eqref{eq:Sep21}, and hence the Gronwall inequality yields $\sup_{t\in (0,T_{\max}^+)}|y(t)|<\infty$. Therefore, we have $\sup_{t\in (0,T_{\max}^+)}\nor{\bom(t)-\bom_0}{H^1}<\infty$, and hence Proposition \ref{P:wp2} implies $T_{\max}^+ =  \infty$.\par
%}\par
%
%
%, which concludes (i). 
% {\color{gray}
% Next, let us show $T_{\max}^\pm=\pm\infty$. Since $\bom(t)$ remains bounded, the only possible scenario is the blow-up of the modulation parameter $y(t)$. For this, we claim that
% \begin{align}\label{eq:Aug103}
% |\dot{y}(t)|\le C \nor{\bom(t)-\bom_0}{H^1_x},\qquad t\in (T_{\max}^-,T_{\max}^+),
% \end{align}
% where $C>0$ is independent of $\bom$ and $t$. Once \eqref{eq:Aug103} is shown, we have
% \begin{align*}%\label{eq:Aug104}
% |\dot{y}(t)| \le C \nor{\bom (t)-\bom_0}{H^1} 
% &\le C \nor{\bom(t)-\tau_{y}\bom_0}{H^1} + C\nor{\tau_y\bom_0-\bom_0}{H^1}\\
% &\le C\del + C|y(t)|,
% \end{align*}
% and hence the Gronwall inequality yields $\sup_{t\in (T_{\max}^-,T_{\max}^+)}|y(t)|<\infty$. Therefore, we have $\sup_{t\in (T_{\max}^-,T_{\max}^+)}\nor{\bom(t)-\bom_0}{H^1}<\infty$, and hence Proposition \ref{P:wp2} implies $T_{\max}^\pm = \pm \infty$.
% }\par 
%
%
To show \eqref{eq:Aug103}, 
differentiating 
%$F(\bom(t),y(t))=0$ 
\eqref{eq:Aug107} in $t$ yields
\begin{gather*}
\int_\R \rd_t\bom \cdot \tau_y\bom_0'dx +\int_\R \dot{y} \tau_y\bom_0'\cdot \tau_y\bom_0'dx +\int_\R (\bom-\tau_y\bom_0 )\cdot \dot{y}\tau_y\bom_0''dx=0\\
\Longleftrightarrow \quad  
\left(\nor{b}{L^2}^2 + \int_\R (\bom-\tau_y\bom_0 )\cdot \dot{y}\tau_y\bom_0''dx\right) \dot{y} = - \int_\R \rd_t\bom \cdot \tau_y\bom_0'dx.
\end{gather*}
Here, we have
\begin{align*}
\nor{b}{L^2}^2=2,
&&\int_\R (\bom-\tau_y\bom_0 )\cdot \dot{y}\tau_y\bom_0''dx \le C \nor{\bom-\tau_y\bom_0}{H^1} |\dot{y}| \le C\del_1 |\dot{y}|. 
\end{align*}
Therefore, if $\del_1$ is sufficiently small, we have
\[
\dot{y} = \left(-\frac 12 + O(\del_1)\right) \int_\R \rd_t\bom\cdot \tau_y\bom_0'dx. 
\]
Now we substitute \eqref{E5.1Mar}. Writing $\bom=\bom_0+\bxi$, we have
\begin{align*}
&\int_\R \bom\times(\bom''+m_3\hat\boe_3) \cdot \tau_y\bom_0'dx 
%&=  \int_\R (\bom_0+\bxi)\times (\bom_0''+\bxi''+(p+\xi_3)\hat\boe_3)\cdot \tau_y\bom_0'dx
\\
&=
%\int_\R \bom_0\times \bom_0''\cdot \tau_y\bom_0'dx
%+
\int_\R (\bom_0+\bxi)\times \bxi''\cdot \tau_y\bom_0'dx
+\int_\R (\bom_0+\bxi)\times \xi_3\hat\boe_3\cdot \tau_y\bom_0dx\\
&\qquad +\int_\R \bxi\times (\bom_0''+p\hat\boe_3)\cdot \tau_y\bom_0'dx\\
&\eqqcolon X_1+X_2+X_3.
\end{align*}
Here, integration by part yields
\begin{align*}
|X_1| &= \left|\int_\R -(\bom_0'+\bxi')\times \bxi'\cdot \tau_y\bom_0' - (\bom_0+\bxi)\times \bxi'\cdot \tau_y \bom_0'' dx \right|\\
&= \left|\int_\R -\bom_0'\times \bxi'\cdot \tau_y\bom_0' - (\bom_0+\bxi)\times \bxi'\cdot \tau_y \bom_0'' dx \right| \le C \nor{\bxi'}{L^2_x}.
\end{align*}
On the other hand, we have 
$|X_2|+|X_3|\le C \nor{\bxi}{L^2_x}$. 
% \begin{align*}
% |X_2|+|X_3|\le C \nor{\bxi}{L^2_x}.
% \end{align*}
Similarly,
\begin{align*}
&\left|\int_\R \bom\times (\hat\boe_1\times\bom') \cdot \tau_y\bom_0'dx\right|\\
% &= \left|\int_\R (\bom_0+\bxi)\times (\hat\boe_1\times (\bom_0'+\bxi'))\cdot \tau_y \bom_0'dx\right|\\
&=
\left|\int_\R \bxi\times (\hat\boe_1\times \bom_0')\cdot \tau_y \bom_0'dx + \int_\R (\bom_0+\bxi)\times (\hat\boe_1\times \bxi')\cdot \tau_y \bom_0'dx\right|\le C\nor{\bxi}{H^1_x}.
\end{align*}
Therefore, we obtain \eqref{eq:Aug103}, and hence Theorem \ref{T2.1} (i) is concluded.
% \begin{claim*}
% $T_{\max}= +\infty$. 
% \end{claim*}
% The proof will be given later.
%
% {\color{red}
% \[
% \nor{\bom(t)-\bom_0}{H^1} \le \nor{\bom(t)-\tau_y\bom_0}{H^1}+ \nor{\tau_y\bom_0 - \bom_0}{H^1} 
% \le C\del_0 + C|y|
% \]
% {\bf We are not able to bound $y(t)$.}}
\par
Finally, we shall conclude (ii). Suppose $\eps>\eps_*$, and let $\bom\in \boM_{\sym}$. Then $u_1 =m_1 \in H^1_{odd} (\R)$, and $u_2= -pm_2 + bm_3 
\in H^1_{odd} (\R)$. 
%If $\bom \in \boM_\sym$, then
% \begin{align*}
% u_1 =%\bom\cdot \hat\boe_1 =
% m_1 \in H^1_{odd} (\R), 
% &&
% u_2 %= \bom\cdot \bom_0^\perp 
% = -pm_2 + bm_3 
% \in H^1_{odd} (\R).
% \end{align*}
In particular, we have $u_1,u_2\perp b$ in $L^2(\R)$ since $b$ is even function. Therefore, by fixing $y=0$, the same coercivity estimate as in Lemma \ref{C2.1} holds. The remaining part of the proof is the same as above.

% Hence the same argument as above shows
% \[
% \nor{\bom(t)-\bom_0}{H^1}=  \nor{\bxi}{H^1} \le C \nor{\bom|_{t=0}-\bom_0}{H^1},\quad t\in (0,T_{\max}^+),
% \]
% %Moreover, by Proposition \ref{P:wp2}, this implies $T_{\max}^+=+\infty$, 
% which completes the proof.

\section{Orbital instability
%for $\eps<\eps_*$
}
\label{S5}
%
%
%\subsection{Main Theorem II: Orbital instability for $\eps<\eps_*$}
%\subsubsection{Outline}
%
Suppose $\eps<\eps_*$. 
%Then $\boA_\eps$ has an unstable eigenvalue. 
We claim the following:
%
%\begin{screen}
\begin{thm}\label{T3.1}
(Orbital instability for $\eps<\eps_*$) 
If $\eps<\eps_*$, then $\bom_0$ is orbitally unstable.
\end{thm}
%\end{screen}
%
%
Our proof follows the framework of Georgiev and Ohta \cite{MR2916078}, exploiting the unstable mode of $\boA_\eps$ established in Proposition \ref{Tx1}. However, this approach cannot be applied directly to \eqref{E2.2Apr}, as second-order derivative nonlinearities obstruct standard nonlinear estimates. To circumvent this, we apply a stereographic projection to the $\S^2$-valued map $\bow$:
%introduce a nonlinear transformation of the form
\begin{equation}\label{eq:Aug3}
v_1 = \frac{u_1}{2+u_3},\qquad v_2= \frac{u_2}{2+u_3}.
\end{equation}
% This coincides with the stereographic projection for the $\S^2$-valued map $\bm{w}=(u_1,u_2,1+u_3)$. 
Then we get
\begin{align}\label{E5.5Apr}
\begin{cases}
\rd_t v_1 = - L_0 v_2 + \tilde{N}_1(u),\\
\rd_tv_2 = L_\eps v_1 + \tilde{N}_2(u),
\end{cases}
\end{align}
where
\begin{align*}
\tilde{N}_1(u) &= 
\frac{2u_2'u_3'}{(2+u_3)^2} - \frac{2u_2(u_3')^2}{(2+u_3)^3}
+\frac{((u_1')^2+(u_2')^2+(u_3')^2)u_2}{(2+u_3)^2}\\
&\qquad 
+\frac{-2(b-\eps) u_1u_1'
-2bpu_2^2 + (2b^2-1)u_2u_3 }{2+u_3}
-\frac{u_1(-2(b-\eps)u_1u_3'-2\eps b u_1u_2)}{(2+u_3)^2},
\\
\tilde{N}_2(u) &= 
-\frac{2u_1'u_3'}{(2+u_3)^2} + \frac{2u_1(u_3')^2}{(2+u_3)^3}
-\frac{((u_1')^2+(u_2')^2+(u_3')^2)u_1}{(2+u_3)^2}\\
&\qquad - \frac{- 2(b-\eps)u_1u_2' +2bp u_1u_2  +(1-2b^2+2\eps b) u_1u_3}{2+u_3} + \frac{u_2(2(b-\eps)u_1u_3'+2\eps b u_1u_2)}{(2+u_3)^2}.
\end{align*}
This formulation is advantageous because it eliminates second-order derivative nonlinearities. Using \eqref{E5.5Apr}, we bound the low-frequency component by the $H^1$-norm (see \eqref{4.7}), while the high-frequency part is controlled by energy conservation (Proposition \ref{P4.6}). Note that \eqref{eq:Aug3} and \eqref{E5.5Apr} will also be reused in the proof of asymptotic stability.
\subsection{
%Preliminaries
Key estimate}
\label{S4.2}
Let $\la_*$ be an eigenvalue with maximal real value; namely, $\Re \la_* = \max \{ \Re \la \ |\ \la\in \si(\boA_\eps) \}$.
% : 
% \[
% \Re \la_* = \max \{ \Re \la \ |\ \la\in \si(\boA_\eps) \}.
% \]
Set 
\begin{align*}
k= 1\qquad \text{ if } \Im \la_* =0,&& 
k= \exp\left(
\frac{2\pi \Re \la_*}{|\Im \la_*|} 
\right) \qquad \text{ if } \Im \la_* \neq 0.
\end{align*}
Let $\phi = 
\begin{bmatrix}
\phi_1 \\ \phi_2
\end{bmatrix}
$ be a corresponding eigenfunction. 
We may assume $\Re \phi\neq 0$ by phase correction if necessary.
% \dlt{$\la_*$ and $\phi$ might be complex-valued, \dlt{but} we may assume $\Re \phi\neq 0$ by phase correction.} 
\par
Now, let $\eps_0,\del>0$ be two numbers satisfying $0<k\del<\eps_0 (\ll 1)$. 
Then, let $\bom_\del$ be the solution to \eqref{E5.1Mar} with initial data
\[
\bom_\del (0) = \frac{1-\del^2 |\Re\phi|^2}{1+\del^2|\Re\phi|^2}\bom_0 + \frac{2\del \Re\phi_1}{1+\del^2|\Re\phi|^2} \hat\boe_1 + 
\frac{2\del \Re\phi_2}{1+\del^2|\Re\phi|^2} \hat\bom_0^\perp.
\]
%\dlt{Note that} 
By virtue of Proposition \ref{P:wp1}, 
$\bom$ is a smooth solution since $\phi\in H^\infty$. We write it as
\[
\bom_\del = \bom_0 + \bxi ,\qquad 
\bxi = u_1 \hat\boe_1 + u_2 \bom_0^\perp + u_3 \bom_0.
\]
Suppose $\del$ is sufficiently small such that $|u(0,x)| < \frac 12$. 
Let $v$ be defined as in \eqref{eq:Aug3}. Then $v$ solves \eqref{E5.5Apr} with initial data $v|_{t=0} = \del \Re \phi$.

%\subsection{
%Proof part 2: Choice of $T_\del$
%}\label{S4.3}
%
\begin{lem}\label{L4.1}
For $\eps_0, \del>0$ as above, there exists $T_\del>0$ such that
\[
\log \frac{\eps_0}{k\del} \le \Re \la_* T_\del \le \log \frac{\eps_0}{\del},\qquad 
\Im \la_* T_\del \in 2\pi \Z.
\]
\end{lem}
\begin{proof}
If $\Im \la_*=0$, 
%\dlt{this is clear} 
then the conclusion follow immediately. If $\Im \la_*\neq 0$, the first condition 
% \dlt{for the real part} 
is translated into
\[
\frac{1}{\Re \la_*} \log \frac{\eps_0}{\del} - \frac{2\pi}{|\Im \la_*|} \le T_\del \le 
\frac{1}{\Re \la_*} \log \frac{\eps_0}{\del}.
\]
Since the length of this interval is $\frac{2\pi}{|\Im \la_*|}$, there must be a integer $n$ such that 
\[
\frac{1}{\Re \la_*} \log \frac{\eps_0}{\del} - \frac{2\pi}{|\Im \la_*|} \le 
\frac{2\pi}{|\Im \la_*|} n \le
\frac{1}{\Re \la_*} \log \frac{\eps_0}{\del}.
\]
Then, $T_\del = \frac{2\pi}{|\Im \la_*|} n$ satisfies the desired properties.
\end{proof}

%\subsection{Key estimate}\label{S4.4}
%
The key estimate in this framework is the following:
%\begin{screen}
\begin{prop}\label{P4.2}
Let $T_\del$ be as in Lemma \ref{L4.1}. 
There exist $c,C>0$ such that if $\eps_0$ is sufficiently small, independently of $\del$, we have
\begin{equation}\label{4.3}
\nor{u(T_\del)}{H^1} \le c\eps_0,\qquad 
\nor{u(T_\del) - 2\del \Re ( e^{\la_*T_\del} \phi) }{L^2} \le C \eps_0^{2}.
\end{equation}
\end{prop}
%\end{screen}
Deferring the proof to the next subsection, we first demonstrate the orbital instability. 
% Postponing the proof into later subsections, we first observe \dlt{the} \textcolor{blue}{its} implications. 
Denoting 
\[
(\Re \phi)^\perp := \Re \phi -
\nor{b}{L^2}^{-2}
\inp{\Re\phi}{ %{}^t[0,b]
\begin{bmatrix}
0 \\ b
\end{bmatrix} 
 }_{L^2} 
\begin{bmatrix}
0 \\ b
\end{bmatrix} .
\]
we conclude, by \eqref{4.3}, 
\begin{align*}
&\left| \inp{u(T_\del)}{(\Re \phi)^\perp}_{L^2}  - 2\del e^{\Re \la_* T_\del} \nor{(\Re \phi)^\perp}{L^2}^2 \right| \\
&= \left| \inp{u(T_\del) - 2\del \Re (e^{\la_*T_\del })\phi }{ (\Re \phi)^\perp}_{L^2} \right|
\le C \nor{(\Re \phi)^\perp}{L^2} \eps_0^{2}.
\end{align*}
Therefore,
\begin{equation}\label{4.4}
\begin{aligned}
\inp{u(T_\del)}{(\Re \phi)^\perp}_{L^2} 
&\ge 
2\del e^{\Re \la_* T_\del} \nor{(\Re \phi)^\perp}{L^2}^2 - C \nor{(\Re \phi)^\perp}{L^2}\eps_0^{2}\\
&\ge \frac{2\eps_0}{k} \nor{(\Re \phi)^\perp}{L^2}^2 - C \nor{(\Re \phi)^\perp}{L^2}\eps_0^{2}\\
&\ge \frac{\eps_0}{k} \nor{(\Re \phi)^\perp}{L^2}^2,
\end{aligned}
\end{equation}
provided $\eps_0 \le (kC)^{-1} \nor{(\Re \phi)^\perp}{L^2}$.
%
%
%\subsection{Proof part 4: Conclusion}\label{S4.5}
%
%\dlt{We} \textcolor{blue}{Now, let us} conclude the orbital instability from \eqref{4.3} and \eqref{4.4}. 
Setting
\[
\Theta_\del := \inf_{y\in\R} \nor{\bom_\del (T_\del) - \bom_0(\cdot +y)}{L^2},
\]
we first observe that this infimum is attained at some $y_\del\in \R$. 
Indeed, we have
\begin{align*}
\liminf_{y\to\pm\infty} \nor{\bom_\del (T_\del) - \bom_0(\cdot +y)}{L^2}^2 
&\ge 
\liminf_{y\to\pm\infty} \nor{m_{\del,2} (T_\del) - b(\cdot +y)}{L^2}^2 \\
&= \nor{m_{\del,2} (T_\del)}{L^2}^2 + \nor{b}{L^2}^2 
\ge \nor{b}{L^2}^2 ,
\end{align*}
whereas \eqref{4.3} implies
\[
\Theta_\del \le 
\nor{\bom_\del (T_\del) - \bom_0}{L^2} 
= \nor{\bxi(T_\del)}{L^2} 
\lesssim \nor{u(T_\del)}{L^2} + \nor{u_3(T_\del)}{L^2} 
\lesssim \eps_0.
\]
Thus the minimum exists if $\eps_0$ is sufficiently small. 
Then we have
\[
\Theta_\del = 
\nor{\bxi (T_\del) +\bom_0 - \bom_0(\cdot +y_\del)}{L^2}\le C\eps_0,
\qquad 
\nor{\bxi(T_\del)}{L^2} \le c\eps_0.
\]
Thus
\[
\nor{b - b(\cdot+y_\del)}{L^2} \le 
\nor{\bom_0 - \bom_0(\cdot +y_\del)}{L^2} \le C \eps_0,
\]
and hence we have $y_\del = O(\eps_0)$. 
Therefore, writing $(\Re \phi)^\perp = {}^t[\tilde{\phi}_1,\tilde{\phi}_2]$, we obtain
%Therefore, by \eqref{4.4}, we obtain
\begin{align*}
&\inp{\bom_\del(T_\del) - \bom_0(\cdot +y_\del)}{ \tilde\phi_1\hat\boe_1+ \tilde\phi_2 \bom_0^\perp
%(\Re \phi)^\perp \cdot [\hat\boe_1 ,\bom_0]
}_{L^2} \\
&= \inp{\bxi(T_\del) +\bom_0 - \bom_0(\cdot +y_\del)}{ \tilde\phi_1\hat\boe_1+ \tilde\phi_2 \bom_0^\perp}_{L^2} \\
&= \inp{\bxi(T_\del)}{ \tilde\phi_1\hat\boe_1+ \tilde\phi_2 \bom_0^\perp}_{L^2}
- y_\del \inp{\rd_x \bom_0}{ \tilde\phi_1\hat\boe_1+ \tilde\phi_2 \bom_0^\perp}_{L^2} +
o(\eps_0)\\
&= \inp{u(T_\del)}{(\Re \phi)^\perp}_{L^2} - y_\del 
\inp{
\begin{bmatrix}
0 \\ b
\end{bmatrix}
 }{(\Re \phi)^\perp}_{L^2} + o(\eps_0) \\
&= \inp{u(T_\del)}{(\Re \phi)^\perp}_{L^2}+ o(\eps_0)
\ge \frac{\eps_0}{2k} \nor{(\Re \phi)^\perp}{L^2}^2
\end{align*}
by \eqref{4.4} and by choosing $\eps_0$ sufficiently small. 
This implies that
\[
\Theta_\del = \nor{\bom_\del(T_\del) - \bom_0(\cdot +y_\del)}{L^2} \ge C,
\]
where $C>0$ is independent of $\del$. 
Since $\nor{\bom_\del (0)-\bom_0}{H^1} \lesssim \del$ and since $\del$ can be taken arbitrarily small, 
% \dlt{this concludes the orbital instability} 
we obtain the desired conclusion.

\subsection{Proof of Proposition \ref{P4.2}}\label{S4.6}

Since $v(0)= \del \Re \phi$, we can rewrite \eqref{E5.5Apr} as
\begin{align*}
v(t) 
% &= \del e^{t\boA_\eps} \Re \phi + \int_0^t e^{(t-s)\boA_\eps} \tilde{N}(u(s)) ds\\
= \del \Re (e^{\la_* t} \phi) + \int_0^t e^{(t-s)\boA_\eps} \tilde{N}(u(s)) ds.
\end{align*}
We first set three constants $\la,\nu,\mu$ satisfying
\[
\frac \mu2<\la <\Re \la_*<\nu < \mu < 2\Re\la_*.
\]
Also, for $\la>0$, we write $e^{\pm}_\la (t) = e^{\pm \la t}$. 
Then we have the following family of estimates:
\begin{lem}
\[
\nor{e^{t\boA_\eps}}{L^2\to L^2} \le C e^{\nu t},\qquad t>0.
\]
\end{lem}
\begin{lem}
If $u=u(t,x)$ satisfies
\[
u(t) = e^{t\boA_\eps} \psi + \int_0^t e^{(t-s)\boA_\eps} f(s) ds,
\]
with $\psi(x)\in L^2$ and $f(t,x)\in L^1_tL^2_x+ L^{\frac 43}_tL^1_x$, then 
\begin{equation}\label{4.5}
\nor{u(t)}{L^2} \le Ce^{\nu t} \nor{\psi}{L^2} + C e^{\mu t} \nor{e_\mu^- f}{L^1_tL^2_x + L^{\frac 43}_tL^1_x}
\end{equation}
\end{lem}
The proof of these lemmas can be found in \cite{MR2916078}. 
The necessity of $\Re\la_*$ being the maximum lies in the conditions for these lemmas. 
We remark that 
%the dual Strichartz norm in the bound 
\eqref{4.5} follows simply from the Strichartz estimate for the free Schr\"odinger operator, since the potential term in $A_\eps$ can be handled perturbatively %\dlt{thanks to the help of} 
with the aid of the exponential factor.\par
Now set
\[
T_* := \sup \left\{T\in (0,T_{\max}) \mid\sup_{0\le t\le T}\nor{u(t)}{L^\infty} \le \frac 12  \right\}.
\]
Note that $T_*>0$ 
%\dlt{, since} 
by $|u(0,x)|< \frac 12$ and the continuity of $t\mapsto  u(t) \in H^1\subset L^\infty$. 
%\dlt{and $\nor{u(t)}{H^1}$ is continuous with respect to $t$.} 
%\dlt{Let $T<T_*$.} 
Now, fix an arbitrary $T\in (0,T_{\max})$. 
For $p\in[1,\infty]$ and for
the Lebesgue and Sobolev spaces $Y$, we write 
$\nor{f}{L_T^pY} := \left[\int_0^T \nor{f(t,\cdot)}{Y}^p  dt \right]^{\frac 1p}$. 
% \[
% \nor{f}{L_T^pY} := \left[\int_0^T \nor{f(t,\cdot)}{Y}^p  dt \right]^{\frac 1p}.
% \]
%Also, 
Set
\[
\nor{u}{X_T} := e^{\la T} \nor{e_\la^- u}{L^\infty_T H^1_x}.
\]
Here, by noting \eqref{eq:Aug5},  $\tilde{N}$ consinsts of the terms of the form
\[
u^2,\qquad uu',\qquad (u')^2,
\]
multiplied by functions bounded uniformly in $[0,T_*]\times \R$. 
With this, 
\eqref{4.5} and \eqref{eq:Aug5} implies
\begin{align*}
&e^{-\la t} \nor{u(t)}{L^2} \le Ce^{-\la t} \nor{v(t)}{L^2}\\
&\lesssim 
\del e^{(\Re \la_* -\la) t} \nor{\phi}{L^2} 
+  e^{(\mu -\la)t} \nor{e_\mu^- \tilde{N}(u)}{L^1_tL^2_x+ L^{\frac 43}_tL^1_x} \\
&\lesssim \del e^{(\Re \la_* -\la) t} \nor{\phi}{L^2} + e^{(\mu-\la)t} \nor{e_\mu^- (u^2+uu')}{L^1_tL^2_x} 
+ e^{(\mu-\la)t} \nor{e_\mu^- (u')^2}{L^{\frac 43}_tL^1_x}\\
&\lesssim 
\del e^{(\Re \la_* -\la) t} \nor{\phi}{L^2} 
+ e^{(\mu -\la)t} \int_0^t e^{-\mu s} \nor{u(s)}{H^1}^2 ds + 
e^{(\mu -\la)t} \left[\int_0^t e^{-\frac 43\mu s} \nor{u(s)}{H^1}^{\frac 83} ds\right]^{\frac 34}
\\
&\lesssim  
\del e^{(\Re \la_* -\la) t} \nor{\phi}{L^2} 
+  e^{(\mu -\la)t} \int_0^t e^{(2\la-\mu) s} ds \nor{e_\la^- u}{L^\infty_TH^1_x}^2 
+ e^{(\mu -\la)t} \left[\int_0^t e^{\frac 43 (2\la-\mu) s} ds 
\right]^{\frac 34} \nor{e_\la^- u}{L^\infty_TH^1_x}^2 \\
&\lesssim 
\del e^{(\Re \la_* -\la) T} \nor{\phi}{L^2} 
+ 
e^{\la T} \nor{e_\la^- u}{L^\infty_TH^1_x}^2,
\end{align*}
which yields
\begin{equation}\label{4.7}
e^{\la T} \nor{e_\la^- u}{L^\infty_TL^2_x} 
\lesssim 
\del e^{\Re \la_* T} \nor{\phi}{L^2} 
+ \nor{u}{X_T}^2.
\end{equation}
For the high frequency part, we claim the following:
\begin{prop}\label{P4.6}
For $T\in [0,T_*]$, we have
%There exists $K\ge \frac 32$ such that
\begin{equation}\label{4.8}
e^{\la T} \nor{e_\la^- u'}{L^\infty_TL^2_x} 
\le 
C e^{\Re \la_* T} \del 
+ C \left( \nor{u}{X_T}^{\frac 32} + \nor{u}{X_T}^2 \right).
\end{equation}
\end{prop}
%We \dlt{will} prove this \textcolor{blue}{claim} in Section \ref{S4.7}. 
The proof is deferred to the end of this subsection.
Combining \eqref{4.7} and \eqref{4.8}, we have
\[
\nor{u}{X_T}
\le 
C \del e^{\Re \la_* T} 
+ C \left( \nor{u}{X_T}^{\frac 32} + \nor{u}{X_T}^2 \right).
\]
Especially if $T< \min\{T_\del,T_*\}$, we have
\[
\nor{u}{X_T}
\le 
C \eps_0 
+ C \left( \nor{u}{X_T}^{\frac 32} + \nor{u}{X_T}^2 \right).
\]
Here, note that
\begin{equation}\label{eq:Aug1}
\nor{u}{L^\infty_TH^1} \le e^{\la T} \nor{e^-_\la u}{L^\infty H^1} \le \nor{u}{X_T}.
\end{equation}
Therefore, if $\eps_0$ is sufficiently small, 
the continuity argument, combined with \eqref{eq:Aug1}, implies $T_\del< T_*$, and 
%\[
$\nor{u}{X_{T_\del}}
\le 
C \eps_0$. 
%\]
In particular, \eqref{eq:Aug1} implies 
$\nor{u}{L^\infty_{T_\del}H^1} \le C\eps_0$. 
On the other hand, 
the same estimate as \eqref{4.7} yields
%\eqref{4.5} and \eqref{4.6} yield
\begin{align*}
e^{-\la t}\nor{v(t) - \del \Re (e^{\la_*t} \phi)}{L^2} 
&\le 
C e^{(\mu-\la) t} \nor{e_\mu^- \tilde{N}(u)}{L^1_t L^2_x + L^{\frac 43}_tL^1_x} \le C e^{\la T_\del} \nor{e_\la^- u}{L^\infty_{T_\del}H^1_x}^2 ,
\end{align*}
and thus
\begin{align*}
\nor{u(T_\del) - 2\del \Re (e^{\la_* T_\del} \phi)}{L^2} &
\le \nor{u-2v}{L^\infty_{T_\del}L^2} +
2\nor{v(T_\del) - \del \Re (\la_* T_\del \phi)}{L^2}
\\
&\le C\nor{u}{L^\infty_{T_\del}H^1}^3+
Ce^{ \la T_\del}
\nor{e^{-\la t} (v(t) - \del \Re (\la_*t \phi))}{L^\infty_{T_\del} L^2} \\
&\le C \left( \nor{u}{X_{T_\del}}^2 + \nor{u}{X_{T_\del}}^3 \right) \le C \eps_0^2.
\end{align*}
Hence the assertion in Proposition \ref{P4.2} follows.

\begin{proof}[Proof of Proposition \ref{P4.6}]
%\label{S4.7}
%
%
Let us recall the Hessian $H(\bxi)$, defined by the energy expansion as in \eqref{eq:Sep22}. 
% \dlt{Let us recall the energy expansion:} %{\color{blue}\eqref{eq:Sep22}}.
% {\color{gray}
% \begin{align*}
% H(\bxi) &= 2(E(\bom_0+\bxi) -E(\bom_0))\\
% &=
% \int_\R |\rd_x \bxi|^2 - \xi_3^2 + 2\eps \rd_x\bxi \times \bxi \cdot \hat\boe_1 
% + \left(1-2b^2  +2\eps b\right) |\bxi|^2 dx.
% \end{align*}
% }
Note that $H(\bxi)$ remains constant in $t$ by the conservation of energy for \eqref{E5.1Mar}. 
% \dlt{
% By the energy conservation for \eqref{E5.1Mar}, 
% $H(\bxi(t))$ is independent of $t$.} 
Recalling \eqref{eq:Aug77}, we can write
% \dlt{
% By recalling the moving frame reduction, we have}
\begin{align}\label{eq:Aug4}
H(\bxi) = \int_\R |\rd_x u|^2 + (1-2b^2)|u|^2 + 2\eps bu_1^2 + O(|u||u'|^2+|u|^2|u'| + |u|^3)dx.
\end{align}
Therefore, for $t\in [0,T]$, we have
\begin{align*}
e^{-2\la t}\nor{\rd_x u(t)}{L^2}^2
&\le e^{-2\la t} |H(\bxi(t))| + C \nor{e_\la^- u(t)}{L^2}^2 + C e^{\la t}\nor{e_\la^- u(t)}{H^1}^3 \\
&\le |H(\bxi (0))| + C 
\nor{e_\la^-u}{L^\infty_TL^2_x}^2 + C e^{\la T} \nor{e_\la^- u}{L^\infty_TH^1_x}^3.
\end{align*}
Here, \eqref{eq:Aug4} yields 
%\[
$|H(\bxi(0))| \le C \del^2$ 
%\]
by using $\del\lesssim 1$. Therefore, 
%multiplying $e^{2\la T}$ and 
using \eqref{4.7}, we obtain
\begin{align*}
e^{2\la T} \nor{e_\la^-u'}{L^\infty_TL^2}^2 
&\le Ce^{2\la T} \del^2 \nor{\phi}{H^1}^2 +
C e^{\la T}\nor{e_\la^-u}{L^\infty_TL^2_x}^2 
+ C \nor{u}{X_T}^3\\
&\le C (\del e^{\Re \la_* T})^2 + C (\nor{u}{X_T}^3+ \nor{u}{X_T}^4),
\end{align*}
which implies \eqref{4.8}.
\end{proof}

\section{Virial type estimate}
\label{S6}
Toward the asymptotic stability, we introduce a key estimate, which is a new extension of a virial-type estimate in \cite{MR3630087} to the Landau--Lifshitz equation. The main %\dlt{statement} 
result in this section is as follows:
\begin{prop}[Virial type estimate]\label{P4:4.1}
There exist $\mu_0, A_0, C>0$ such that the following holds: 
Let $T,A,\del>0$ satisfy $A>A_0$ and $A\del<\mu_0$. Let $(u_1,u_2,u_3)$ be a solution to \eqref{E2.2Apr} with $\nor{u}{L^\infty_TH^1_x}\le \del.$ Then, we have
% {\color{gray}
% Let $T,A,\del_0>0$ be constants satisfying 
% Let $(u_1,u_2,u_3)$ be the solution to \eqref{E2.2Apr} with $\nor{u}{L^\infty_tH^1_x}\le C\del.$ 
% Let $b_A = b(\cdot/A)$. Then for $A,\del>0$ with 
% $A\gg 1$ and $A\del \ll 1$, we have}
\begin{align}\label{E4.1Apr}
\int_0^T \left(
\nor{b_A u'}{L^2}^2
+\frac{1}{A^2}\nor{b_A u}{L^2}^2
\right) dt
\le 
CA\del^2 + C \int_0^T \nor{b_{10} u}{L^2}^2 dt,
\end{align}
where $b_A = b(\frac \cdot A)$.
\end{prop}
Since $b_{10}$ is more localized than $b_A$, \eqref{E4.1Apr} is interpreted as 
\begin{center}
\text{Weighted $L^2_tH^1_x$ Bounds} $\lesssim$ 
\text{$L^\infty_tH^1_x$-bound} $+$ 
\text{$L^2_tL^2_x$-bounds near the origin}.
\end{center}

%{\color{red}
Consequently, this estimate bounds the high-frequency component and spatial tails using the localized information of the low-frequency part. This virial-type method was originally developed for nonlinear Klein--Gordon equations \cite{MR3630087}, and later extended to standing wave perturbations for nonlinear Klein--Gordon \cite{MR4404795} and Schr\"odinger (NLS) equations \cites{MR4552229,MR4864876,MR4685914,MR4910476,MR4777089}. However, directly applying this framework to our setting presents significant obstructions. 

To clarify these challenges, recall the standard approach for the NLS perturbation equation \cites{MR4404795,MR4552229}:
\begin{equation}\label{eq:m5}
i \rd_tu = (- \rd_x^2 +1) u + \big(f(Q+u)-f(Q)\big), \quad u=u_1+iu_2,
\end{equation}
where $f$ is the nonlinearity and $Q$ is the standing wave profile. The key idea relies on the virial functional
\[
\boI(u) = \int_\R (u_1S_A u_2 - u_2S_Au_1) \,dx = \int_\R 2\vph_A (u_1u_2'-u_2u_1') \,dx,
\]
where $S_A = 2\vph_A\rd_x + \vph_A'$ and $\vph_A$ is a weight function satisfying $\vph_A' \sim b_A^2$. The time derivative of $\boI$ yields the lower bound
\begin{align}\label{E4.2Apr}
-\frac 12 \frac{d}{dt} \boI(u) \ge \nor{b_A u'}{L^2}^2 + \frac{1}{A^2}\nor{b_Au}{L^2}^2 + \inp{S_A\big(f(Q+u)-f(Q)\big)}{u}.
\end{align}
Using integration by parts, the last term in \eqref{E4.2Apr} can be rewritten as
\begin{align}\label{eq:Aug82}
\int_\R \big(\vph_A' O(u^2) + \vph_A Q O(u^2)\big) \,dx,
\end{align}
which is then treated perturbatively to establish \eqref{E4.1Apr}.

Applying this standard functional directly to \eqref{eq:Aug81} fails because the second-order derivative nonlinearities (e.g., $u_2u_3'' - u_3u_2''$) cannot be bounded perturbatively. Alternatively, evaluating $\boI$ on the semilinearized equation \eqref{E5.5Apr} also fails, as the transformed nonlinearity lacks the algebraic structure required to produce a remainder of the form \eqref{eq:Aug82}. 

To overcome these structural obstacles, we introduce a \textit{novel virial functional}:
\begin{align}\label{eq:Aug83}
\boI(u) = \int_\R \frac{u_1 S_Au_2 - u_2S_Au_1}{2+u_3} \,dx = \int_\R 2\vph_A \frac{u_1u_2' - u_2u_1'}{2+u_3} \,dx.
\end{align}
We will show that this modified functional recovers a coercivity estimate analogous to \eqref{E4.2Apr} while effectively canceling the problematic nonlinearities. Namely, we obtain
%}

\begin{equation}\label{E7.1Apr}
-\frac 12\frac{d}{dt} \boI(u) \ge \nor{b_A u'}{L^2}^2 + \frac{1}{A^2}\nor{b_Au}{L^2}^2
+ \text{Harmless Remainders.}
\end{equation}
%\end{screen}
%\end{itembox}
It is worth remarking that 
\eqref{eq:Aug83} 
%this definition of virial is 
is inspired by the Haldene's momentum functional for Schr\"odinger maps \cite{PhysRevLett.57.1488}:
\[
\int_\R \frac{m_1\rd_x m_2 - m_2\rd_xm_1}{1+m_3} dx,\qquad 
\bom = {}^t (m_1,m_2,m_3)\colon \R\to\S^2.
\]
\subsection{%Preliminaries
Virial functional}\label{S6.1}
%
%
% In the sequel, we give a proof of Proposition \ref{P4:4.1}. 
% For this, \textcolor{blue}{we introduce some notations and recall some properties from \cite{MR4404795}.}
% \dlt{we prepare some notations and recall their properties.} 
%the virial-type estimate \eqref{E4.1Apr}. 
%
%

Let $\chi:\R\to\R$ be smooth even cut-off function with
\[
\chi(x)=1 \text{ if } |x|\le 1,\qquad \chi(x)=0 \text{ if } |x|\ge 2,\qquad \chi'\le 0 \text{ on } [0,\infty).
\]
For given $A\gg 1$, set
\begin{align*}
\zeta_A(x) 
\coloneqq \exp (-\frac{|x|}{2A} (1-\chi(x))),\qquad 
\vph_A(x) \coloneqq \int_0^x \zeta_A^2(y) dy.
\end{align*}
We define the virial functional
\begin{align}\label{eq:Sep9}
\boI(u) 
= \int_\R 2\vph_A \frac{u_1u_2' - u_2u_1'}{2+u_3}dx
= \int_\R \frac{u_1 S_Au_2 - u_2S_Au_1}{2+u_3} dx
,\qquad S_A= 2\vph_A\rd_x +\vph_A'.
\end{align}

% In the sequel, we give a proof of Proposition \ref{P4:4.1}. 
% For this, 

% \textcolor{blue}{we introduce some notations and recall some properties from \cite{MR4404795}.}
% \dlt{we prepare some notations and recall their properties.} 

We first recall the properties of $\zeta_A$ and $\vph_A$.
\begin{lem}
There exists $C>0$, independent of $x\in\R$ and $A\ge 1$, such that the following hold:
\begin{gather}
0< \vph_A'(x) = \zeta^2_A(x)\le 1,\qquad |\vph_A(x)| \le |x|,\label{eq:Aug15}\\
|\zeta_A'(x)|\le \frac{C}{A} \zeta_A(x),
\label{eq:Aug12}\\
|\vph_A(x)| + |x\vph_A'(x)| \le CA,\label{eq:Aug16}\\
C^{-1} b_A (x) \le \zeta_A^2(x) \le C b_A(x) \label{eq:Aug14}.
\end{gather}
\end{lem}
For the proof, we refer to \cite{MR4404795}.

%%%
%%% Proof of Lemma 6.2
%%%
% \begin{proof}
% \eqref{eq:Aug15} is immediate. \eqref{eq:Aug12} follows from $\zeta_A'/\zeta_A = \frac{1}{2A} |\rd_x (|x| (1-\chi))|$. If $x\ge 2$,
% \begin{align}
% |\vph_A(x)| 
% &\le \int_0^2 |\vph_A'(t)|dt dx + \int_2^x |\vph_A'(t)| dt\\
% &\le 2 + \int_2^x \exp(-\frac{t}{A})dt \le 2+ A\exp (-\frac 2A) \le CA.
% \end{align}
% If we write
% \begin{gather*}
% |x\vph_A'(x)| = A \rho_A(\frac{|x|}{A}),\qquad \rho_A (t) \coloneqq t\exp (-t(1-\chi(At))),
% \end{gather*}
% then we have
% \begin{gather}
% |\rho_A(t)|\le t \le \frac 2A
% \quad\text{if } t\le \frac 2A,\qquad |\rho_A(t)|=|te^{-t}| \le C
% \quad\text{if } t\ge \frac 2A,
% \end{gather}
% and hence \eqref{eq:Aug16} follows. Next, we observe that
% \[
% b_A(x) \zeta_A^{-2}(x) \le b(\frac xA) e^{\frac{|x|}{A}} \le C
% \]
% and thus the first inequality in \eqref{eq:Aug14} follows. If $|x|\le 2$, then
% \[
% \zeta_A^2(x) \le 1 \le Cb(x) \le Cb_A(x)
% \]
% for $A\ge 1$. If $|x|\ge 2$, then 
% \begin{align*}
% \zeta_A^2(x) = e^{-\frac{|x|}{A}} \le Cb_A(x).
% \end{align*}
% Hence the second inequality in \eqref{eq:Aug14} follows.
% \end{proof}
%%%
%%% Proof of Lemma 6.2 end
%%%

% Here, we define the virial functional
% \begin{align}\label{eq:Sep9}
% \boI(u) 
% = \int_\R 2\vph_A \frac{u_1u_2' - u_2u_1'}{2+u_3}dx
% = \int_\R \frac{u_1 S_Au_2 - u_2S_Au_1}{2+u_3} dx
% ,\qquad S_A= 2\vph_A\rd_x +\vph_A'.
% \end{align}
%
\subsection{Time derivative of $\boI(u)$}
We compute the time derivative $\boI(u)$. 
For simplicity, 
we abbreviate \eqref{E2.2Apr} as
\begin{align*}
\rd_t u_1
&=
(1+u_3) u_2'' - u_2u_3''  
+N_1(u),\\
\rd_t u_2
&=
u_1u_3'' - (1+u_3)u_1''   +N_2(u), \\
\rd_t u_3 &= u_2 u_1'' - u_1u_2'' + N_3(u),
\end{align*}
where the terms $N_j(u)$, $j=1,2,3$ involve derivatives of $u$ up to order one.
%\dlt{consist of up to first-order derivatives of $u$}.
% with
% \begin{align*}
% N_1(u) &= - (1-2b^2) u_2-u_2u_3 - 2(b-\eps) u_1u_1' + b(-2p u_2^2 + 2b u_2u_3) ,\\
% N_2(u) &= (1-2b^2+2\eps b) u_1+u_1u_3 - 2(b-\eps) u_1u_2' + b(2p u_1u_2 + (-2b+2\eps ) u_1u_3) ,\\
% N_3(u) &= - 2(b-\eps) u_1u_3' + -2\eps b u_1u_2 .
% \end{align*}
Then,
\[
\frac{d}{dt} \boI(u) = \int_\R 2\vph_A \left(
\frac{
\dot{u}_1 u_2' + u_1 \dot{u}_2' - \dot{u}_2u_1' - u_2 \dot{u}_1'
}{2+u_3}
- \frac{\dot{u}_3 (u_1 u_2' -u_2u_1')}{(2+u_3)^2}
\right) dx.
\]
Integration by part yields
\begin{align*}
\int_\R 2\vph_A \frac{u_1 \dot{u}_2' -u_2 \dot{u}_1'}{2+u_3} dx
&= \left[ 2\vph_A \frac{u_1\dot{u}_2 - u_2 \dot{u}_1}{2+u_3}\right]_{-\infty}^\infty + \int_\R \dot{u}_1 \left(\frac{2\vph_A u_2}{2+u_3}\right)'dx - \int_\R \dot{u}_2\left(\frac{2\vph_A u_1}{2+u_3}\right)'dx\\
&= \int_\R 2\vph_A' \frac{u_2\dot{u}_1- u_1 \dot{u}_2}{2+u_3}dx +\int_\R 2\vph_A \left[ \left(\frac{ u_2}{2+u_3}\right)' \dot{u}_1 - \left(\frac{ u_1}{2+u_3}\right)' \dot{u}_2\right]dx.
\end{align*}
Therefore, we have
\begin{align*}
&\frac{d}{dt} \boI(u)\\
&= 
\int_\R 2\vph_A' \frac{u_2\dot{u}_1- u_1 \dot{u}_2}{2+u_3}dx\\
&\qquad+ \int_\R2\vph_A \left(
\left(\frac{ u_2}{2+u_3}\right)' \dot{u}_1 - \left(\frac{ u_1}{2+u_3}\right)' \dot{u}_2+ \frac{\dot{u}_1u_2' - \dot{u}_2 u_1'}{2+u_3} -
\frac{\dot{u}_3 (u_1 u_2' -u_2u_1')}{(2+u_3)^2}
\right)dx\\
&=\int_\R 2\vph_A' \frac{u_2\dot{u}_1- u_1 \dot{u}_2}{2+u_3}dx\\
&\qquad + \int_\R 2\vph_A\left(
\left(\frac{2u_2'}{2+u_3} - \frac{u_2u_3'}{(2+u_3)^2}
\right) \dot{u}_1 
-\left(\frac{2u_1'}{2+u_3} - \frac{u_1u_3'}{(2+u_3)^2}
\right) \dot{u}_2 
- \frac{(u_1 u_2' -u_2u_1')}{(2+u_3)^2}\dot{u}_3 
\right)dx\\
&=\int_\R 2\vph_A' \frac{u_2\dot{u}_1- u_1 \dot{u}_2}{2+u_3}dx\\
&\qquad + \int_\R 2\vph_A
\left(\frac{2u_2'}{2+u_3} - \frac{u_2u_3'}{(2+u_3)^2}
\right) ((1+u_3)u_2'' - u_2u_3'') dx\\
&\qquad + \int_\R 2\vph_A
\left(\frac{2u_1'}{2+u_3} - \frac{u_1u_3'}{(2+u_3)^2}
\right) ((1+u_3)u_1'' - u_1u_3'')dx\\
&\qquad + \int_\R 2\vph_A
\left( \frac{u_2u_1' - u_1u_2'}{(2+u_3)^2}
\right) (u_2u_1'' - u_1u_2'')dx\\
&\qquad + \int_\R 2\vph_A \left(
\left(\frac{2u_2'}{2+u_3} - \frac{u_2u_3'}{(2+u_3)^2}
\right) N_1(u) 
- \left(\frac{2u_1'}{2+u_3} - \frac{u_1u_3'}{(2+u_3)^2}
\right) N_2(u) +
\left( \frac{u_2u_1' - u_1u_2'}{(2+u_3)^2}
\right) N_3(u)
\right)dx
\end{align*}
Now, recall that
\[
u_1^2+ u_2^2+(1+u_3)^2=1,\qquad u_1u_1' +u_2u_2' + (1+u_3)u_3' =0.
\]
By using this, direct computation yields
\begin{align*}
&\left(\frac{2u_2'}{2+u_3} - \frac{u_2u_3'}{(2+u_3)^2}
\right) ((1+u_3)u_2'' - u_2u_3'')+
\left(\frac{2u_1'}{2+u_3} - \frac{u_1u_3'}{(2+u_3)^2}
\right) ((1+u_3)u_1'' - u_1u_3'')\\
&\qquad +\left(\frac{u_2u_1'-u_1u_2'}{(2+u_3)^2}
\right) (u_2u_1''-u_1u_2'')
= u_1'u_1'' +u_2'u_2''+u_3'u_3'' . 
%\frac 12 \left( (u_1')^2 + (u_2')^2 + (u_3')^2 \right)'.
\end{align*}
We also have
\begin{align*}
\int_\R 2\vph_A' \frac{u_2\dot{u}_1- u_1 \dot{u}_2}{2+u_3}dx
&=\int_\R 2\vph_A' \frac{u_2 ((1+u_3)u_2'' - u_2u_3'' +N_1) 
- u_1 (u_1u_3'' - (1+u_3)u_1'' +N_2) }{2+u_3}\\
&= \int_\R 2\vph_A' \frac{
(1+u_3) (u_1u_1'' +u_2u_2'') - (u_1^2+u_2^2) u_3''
}{2+u_3}
+ \int_\R 2\vph_A' \frac{u_2N_1 -u_1 N_2}{2+u_3}
\\
&= \int_\R \vph_A' (u_1u_1''+u_2u_2'') dx+ \int_\R \vph_A' 
\frac{u_3(u_1u_1''+u_2u_2'') -2 (u_1^2+u_2^2)u_3''}{2+u_3}dx\\
&\qquad + \int_\R 2\vph_A' \frac{u_2N_1 -u_1 N_2}{2+u_3}dx.
\end{align*}
Therefore, we obtain
\begin{equation}\label{eq:Aug11}
\begin{aligned}
\frac{d}{dt}\boI(u) 
&= 
\int_\R (S_Au_1)u_1'' + (S_Au_2)u_2'' dx + \int_\R 2\vph_A' \frac{u_2N_1 -u_1 N_2}{2+u_3}dx\\
&\qquad 
+ \int_\R \vph_A' 
\frac{u_3(u_1u_1''+u_2u_2'') -2 (u_1^2+u_2^2)u_3''}{2+u_3}dx
+ \int_\R 2\vph_A u_3' u_3''dx\\
&\qquad 
+ \int_\R 2\vph_A \left(
\left(\frac{2u_2'}{2+u_3} - \frac{u_2u_3'}{(2+u_3)^2}
\right) N_1(u) 
- \left(\frac{2u_1'}{2+u_3} - \frac{u_1u_3'}{(2+u_3)^2}
\right) N_2(u) \right.\\
&\hspace{230pt} \left.+
\left( \frac{u_2u_1' - u_1u_2'}{(2+u_3)^2}
\right) N_3(u)
\right)dx\\
&\eqqcolon \int_\R (S_Au_1)u_1'' + (S_Au_2)u_2'' dx  + \boR(u).
\end{aligned}
\end{equation}
The first term is the essential factor to produce the desired coercivity.
%\dlt{LHS of \eqref{E7.1Apr}.} 
Indeed, 
by \cite{MR4404795}*{Lemma 1}, we have
\begin{align}\label{eq:Aug32}
\begin{aligned}
\nor{(\zeta_Au_1)'}{L^2}^2 + \nor{(\zeta_Au_2)'}{L^2}^2 
&\le
-\frac 12 \int_\R (S_Au_1)u_1'' + (S_Au_2)u_2'' dx 
+\frac{C}{A} \int_\R \mathbbm{1}_{1\le|x|\le 2} |\zeta_A u|^2 dx\\
&\le -\frac 12\int_\R (S_Au_1)u_1'' + (S_Au_2)u_2'' dx 
+ C \int_\R b_{10}^2 |u|^2 dx
\end{aligned}
\end{align}
since $A\gg 1$, $\mathbbm{1}_{1\le|x|\le 2}\lesssim b_{10}^2$, and since $S_A$ is skew-symmetric.
\subsection{Estimate of the remainder}
Now we claim the following:
\begin{lem}
There exists $C>0$, independent of $t,\del,A$, such that  
\begin{align}\label{eq:Aug31}
|\boR(u)| \le C\nor{b_{10}u}{L^2}^2 + C(A\nor{u}{L^\infty} + A^2\nor{u}{L^\infty}^2) \nor{(\zeta_Au)'}{L^2}^2.
\end{align}
\end{lem}
For the terms related to $N_j$, $j=1,2,3$, we separate them into% the following groups:
\[
N_j = \sum_{k=1}^4 N_{kj},\qquad j=1,2,3,
\]
where
\begin{itemize}
    \item $N_{11}= -(1-2b^2)u_2$, $N_{12}=(1-2b^2+2\eps b)u_1$, $N_{13}=0$;
    \item $N_{21}= -u_2u_3$, $N_{22}=u_1u_3$, $N_{23}=0$;
    \item $N_{31}=- 2(b-\eps) u_1u_1'$, $N_{32}=- 2(b-\eps)u_1u_2'$, $N_{33}=- 2(b-\eps) u_1u_3'$;
    \item $N_{41}=bF_1$, $N_{42}=bF_2$, $N_{43}=bF_3$.
    % \item $N_{51}=\dot{y}u_1'$, $N_{52}= \dot{y} (u_2'+b(1+u_3))$, $N_{53}=\dot{y} (u_3'-bu_2)$.
\end{itemize}
For each group, we denote
\begin{align*}
I_k \coloneqq &\int_\R 2\vph_A' \frac{u_2N_{k1} -u_1 N_{k2}}{2+u_3}dx +\int_\R 2\vph_A \left(
\left(\frac{2u_2'}{2+u_3} - \frac{u_2u_3'}{(2+u_3)^2}
\right) N_{k1}(u) \right.\\
&\qquad \left.
- \left(\frac{2u_1'}{2+u_3} - \frac{u_1u_3'}{(2+u_3)^2}
\right) N_{k2}(u) 
+
\left( \frac{u_2u_1' - u_1u_2'}{(2+u_3)^2}
\right) N_{k3}(u)
\right)dx.
\end{align*}
%
%\subsubsection*{Group 1}
%
Then, direct computation yields
\begin{align*}
I_1 &= \int_\R 2\vph_A' \frac{-u_2 (1-2b^2)u_2 -u_1 (1-2b^2+2\eps b) u_1}{2+u_3} dx\\
&\qquad +\int_\R 2\vph_A \left(
- (1-2b^2)\left(\frac{2u_2'}{2+u_3} - \frac{u_2u_3'}{(2+u_3)^2} \right) u_2 \right.\\
&\hspace{70pt} \left.- (1-2b^2+2\eps b)
\left(\frac{2u_1'}{2+u_3} - \frac{u_1u_3'}{(2+u_3)^2} \right) u_1
\right)dx\\
&=\int_\R 2\vph_A' \frac{- (1-2b^2)u_2^2 - (1-2b^2+2\eps b) u_1^2}{2+u_3} dx\\
&\qquad +\int_\R 2\vph_A \left(
- (1-2b^2)\left(\frac{u_2^2}{2+u_3} \right)'  - (1-2b^2+2\eps b)
\left(\frac{u_1^2}{2+u_3} \right)' 
\right)dx\\
&= \int_\R -2(1-2b^2) \left(\vph_A \frac{u_2^2}{2+u_3}\right)' -2 (1-2b^2+2\eps b) \left(\vph_A\frac{u_1^2}{2+u_3}\right)' dx.
\end{align*}
Therefore, integration by parts yields
\begin{align*}
I_1&=\int_\R 8\vph_A b^2p \frac{u_2^2}{2+u_3} + 4\vph_A bp (2b- \eps) 
\frac{u_1^2}{2+u_3}dx\le C \nor{b_{10}^2 |u|}{L^2}^2
\end{align*}
by noting $|\vph_Ab| \le |x|b \lesssim b_{10}^4$. For $I_2$, we have
\begin{align*}
I_2 &= \int_\R 2\vph_A' \frac{-u_2 (u_2u_3) -u_1 (u_1u_3))}{2+u_3} dx\\
&\qquad +\int_\R 2\vph_A \left(
 -\left(\frac{2u_2'}{2+u_3} - \frac{u_2u_3'}{(2+u_3)^2} \right) u_2u_3 
 %\right.\\
% &\hspace{70pt} \left.
-
\left(\frac{2u_1'}{2+u_3} - \frac{u_1u_3'}{(2+u_3)^2} \right) u_1u_3
\right)dx\\
&=\int_\R 2\vph_A' \frac{- u_3 (u_1^2+ u_2^2)}{2+u_3} dx
% \\
% &\qquad 
+\int_\R 2\vph_A \left(
- \left(\frac{u_2^2}{2+u_3} \right)'u_3  - 
\left(\frac{u_1^2}{2+u_3} \right)' u_3
\right)dx\\
&=\int_\R 2\vph_A' u_3^2 dx
% \\
% &\qquad 
+\int_\R 2\vph_A u_3u_3'dx\\
&=\int_\R \vph_A' u_3^2 dx \le C \int_\R \zeta_A^2 |u|^4 dx \le C A^2 \nor{u}{L^\infty}^{2} \nor{(\zeta_Au)'}{L^2}^2,
\end{align*}
where we used the following inequality:
%By \cite{MR4404795}*{Claim 1}, for $\al>0$ we have
\begin{equation}\label{eq:Aug13}
\int_\R \zeta_A^2 |f|^{2\al+2}
\le C_\al 
A^2 \nor{f}{L^\infty}^{2\al} \int_\R |(\zeta_A f)'|^2 dx,\qquad \al>0.
\end{equation}
(See \cite{MR4404795}*{Claim 1}.) 
For $I_3$, noting that the second integral vanishes by cancellation, we have
% \par 
% Next, we observe $I_3$. Indeed, the second integral vanishes by cancellation. Thus,
\begin{align*}
I_3 &= \int_\R 2\vph_A' \frac{- 2(b-\eps) u_2(u_1u_1') + 2(b-\eps) u_1 (u_1u_2')}{2+u_3} dx\\
&= \int_\R -4(b-\eps) \zeta_A \frac{ u_1u_2 (\zeta_A u_1)' -  u_1^2 (\zeta_A u_2)'}{2+u_3} dx\\
&= C \nor{(\zeta_A u)'}{L^2} \left[\int \zeta_A^2 |u|^4 dx\right]^{\frac 12} 
%+\frac{C}{A} \int_\R \zeta_A^2 |u|^3dx \le  %\frac{1}{100} \nor{(\zeta_A u)'}{L^2}^2 
\le CA \nor{u}{L^\infty} 
\nor{(\zeta_Au)'}{L^2}^2.
\end{align*}
Next, we observe $I_4$. 
The first term of $I_4$ is bounded by 
%\[
$C\int_\R \zeta_A^2 b |u|^3 dx \le 
C \nor{u}{L^\infty} \nor{b_{10}^2u}{L^2}^2$. 
%CA^2\nor{u}{L^\infty} \nor{(\zeta_Au)'}
%{L^2}^2
%\]
The second integral of $I_4$ can be bound as
\begin{align*}
C\int_\R |\vph_A| b (|u|^2 |u'| + |u|^4 |u'|)dx
&\le C \int_\R |\vph_A| b \zeta_A^{-1} 
(|u|^2 + |u|^4) |(\zeta_A u)' -\zeta_A' u| dx\\
&\le C \int_R (|u|+|u|^3) |b_{10}^2u| |(\zeta_A u)'| dx
+ C\int_\R \frac{1}{A} (|u|+|u|^3) |b_{10}^2 u|^2\\
&\le C (\nor{u}{L^\infty} + \nor{u}{L^\infty}^3 ) \left( \nor{(\zeta_A u)'}{L^2}^2 + \nor{b_{10}^2 u}{L^2}^2 \right),
\end{align*}
where we used \eqref{eq:Aug12} and 
%\[
$|\vph_A(x) b(x)\zeta_A^{-1}(x)| \le |x|b(x) e^{\frac{|x|}{2A}} \lesssim b_{10}^2(x)$. 
%\]
%
%
Finally, we observe the terms in $\boR(u)$ which the second-derivative of $u$ involves. 
We first have
\begin{multline*}
-\int_\R \vph_A' \frac{2(u_1^2+u_2^2)u_3''}{2+u_3} + \int_\R 2\vph_Au_3'u_3'' 
= \int_\R (2\vph_A' u_3u_3'' + 2\vph_Au_3'u_3'') \\
= \int_\R -2\vph_A''u_3u_3' -3\vph_A'(u_3')^2 = \int_\R -4\zeta_A \zeta_A' u_3u_3' - 3\zeta_A^2 (u_3')^2.
%\\
% &=\int_\R \vph_A''' u_3^2 - 3\vph_A'(u_3')^2\\
% &= \int_\R (\zeta_A^2)'' u_3^2 - 3\zeta_A^2(u_3')^2\\
% &\le \frac{C}{A^2}\int_\R \zeta_A^2|u|^4dx + C\int_\R \zeta_A^2|u|^2|u'|^2dx\\
% &\le C\nor{u}{L^\infty}^2 \nor{(\zeta_Au)'}{L^2}^2 .
%= \int_\R (S_Au_3)u_3''dx.
\end{multline*}
%Each term can be estimated in the same manner as above; indeed, 
Using \eqref{eq:Aug12} and \eqref{eq:Aug13}, we obtain
\begin{align*}
\left| \int_\R \zeta_A \zeta_A' u_3u_3'dx\right| 
&\lesssim \int_\R \frac 1A \zeta_A^2 |u|^3|u'| dx 
= \int_\R \frac 1A \zeta_A |u|^3 |(\zeta_A u)'-\zeta_A'u|dx\\
&\lesssim \int_\R \frac 1A \zeta_A|u|^3 |(\zeta_A u)'|dx + \int_\R \frac 1{A^2} \zeta_A^2|u|^4dx\\
&\lesssim C \nor{u}{L^\infty}^2 \nor{(\zeta_A u)'}{L^2}^2,\\
\int_\R \zeta_A^2(u_3')^2dx&\lesssim
\int_\R \zeta_A^2|u|^2|u'|^2dx =\int_\R |u|^2 |(\zeta_A u)' -\zeta_A' u|^2dx\\
&\lesssim \int_\R |u|^2 |(\zeta_A u)'|^2dx + \int_\R \frac 1{A^2} \zeta_A^2|u|^4dx\\
&\lesssim C \nor{u}{L^\infty}^2 \nor{(\zeta_A u)'}{L^2}^2.
\end{align*}
By integration by parts, we also have
\begin{align*}
\left|\int_\R \vph_A' \frac{u_3u_1u_1''}{2+u_3} dx\right|
&= \left|-\int_\R \left(\zeta_A^2 \frac{u_3u_1}{2+u_3}\right)' u_1' dx\right|\\
&= \left|\int_\R 
-2\zeta_A\zeta_A' \frac{u_3u_1u_1'}{2+u_3} - \zeta_A^2 \frac{(u_3'u_1+u_3u_1')u_1'}{2+u_3}
+ \zeta_A^2 \frac{u_3u_3'u_1u_1'}{(2+u_3)^2}  dx\right| \\
&\lesssim  \int_\R\zeta_A^2 (\frac 1A |u|^3u' +|u|^2|u'|^2+|u|^4|u'|^2)dx. 
\end{align*}
Here, note that
\begin{equation}\label{eq:Aug17}
\zeta_A^2|u'|^2 =|(\zeta_Au)' -\zeta_A'u|^2\le 2|(\zeta_A u)'|^2 + \frac{2}{A^2}|\zeta_A u|^2.
\end{equation}
Hence, by noting $\nor{u}{L^\infty}\ll 1$, we have
\begin{align*}
& \int_\R\zeta_A^2 (\frac 1A|u|^3u' +|u|^2|u'|^2+|u|^4|u'|^2)dx\\
&\lesssim \int_\R \frac{1}{A} \zeta_A|u|^3 |(\zeta_Au)'| + \frac{1}{A^2} \zeta_A^2 |u|^4 + |u|^2 |(\zeta_Au)'|^2
\\
&\lesssim 
\frac 1A\left[\int_\R \zeta_A^2|u|^6 dx\right]^{1/2} \nor{(\zeta_A u)'}{L^2}^2 + \frac{1}{A^2} \int_\R \zeta_A^2 |u|^4dx + 
\nor{u}{L^\infty}^2 \nor{(\zeta_A u)'}{L^2}^2
\\
&\lesssim \nor{u}{L^\infty}^2 \nor{(\zeta_A u)'}{L^2}^2.
%+ \nor{u}{L^\infty}^2\nor{(\zeta_A u)'}{L^2}
% + \frac{1}{A^2}\nor{u}{L^\infty}^2 \nor{\zeta_A u}{L^2}^2.
\end{align*}
The same is true for $\int_\R \vph_A' \frac{u_3u_2u_2''}{2+u_3}dx$.
Consequently, combining all the above estimates concludes \eqref{eq:Aug31}.
\subsection{Conclusive step}
Let us conclude \eqref{E4.1Apr}. 
By \eqref{eq:Aug11}, \eqref{eq:Aug32}, and
\eqref{eq:Aug31}, 
\begin{align*}
\nor{(\zeta_Au)'}{L^2}^2 \lesssim -\frac{d}{dt} \boI(u) + \nor{b_{10}u}{L^2}^2 + C(A\nor{u}{L^\infty} + A^2\nor{u}{L^\infty}^2) \nor{(\zeta_Au)'}{L^2}^2
\end{align*}
by supposing $\nor{u}{L^\infty}\lesssim 1$ and $A\gtrsim 1$. By assumption, we have $A\nor{u}{L^\infty} \lesssim A\nor{u}{H^1}\lesssim A\del$ for $t\in [0,T]$. 
Therefore, if $A\del$ is sufficiently small, we have
\[
\nor{(\zeta_Au)'}{L^2}^2 \lesssim -\frac{d}{dt} \boI(u) + \nor{b_{10}u}{L^2}^2,\qquad t\in [0,T].
\]
Integrating over $[0,T]$ yields
\begin{align*}
\int_0^T  \nor{(\zeta_Au)'}{L^2}^2 dt
\lesssim \boI(u(0)) - \boI(u(T)) + \int_0^T \nor{b_{10}u}{L^2}^2dt.
\end{align*}
Now we have
\[
|\boI(t)| \le C\nor{\vph_A}{L^\infty} \nor{u(t)}{H^1}^2\le CA \del^2 ,\qquad t\in[0,T].
\]
By \cite{MR4404795}*{Lemma 4}, given $w\in H^1(\R)$, we have
\[
\int_\R b_A|w|^2 dx \le CA \nor{b_{10}^2w}{L^2}^2 + CA^2 \nor{w'}{L^2}^2 .
% \le  
% CA \nor{b_{10}^2u}{L^2}^2 + CA^2 \del^2.
\]
By this and \eqref{eq:Aug14}, we have
\begin{align*}
\frac {1}{A^2}\int_\R b_A^2 | u|^2 dx = \frac{1}{A^2}\int_\R b_A\zeta_{A}^{-2} \cdot b_A|\zeta_Au|^2dx
&\le \frac{C}{A^2} \int_\R b_A|\zeta_Au|^2 dx\\
&\le C\nor{(\zeta_A u)'}{L^2}^2 + \frac{C}{A} \nor{b_{10}^2 (\zeta_A u)}{L^2}^2. 
\end{align*}
On the other hand, \eqref{eq:Aug14} and \eqref{eq:Aug17} imply
\begin{align*}
\int_\R b_A^2 |u'|^2 dx 
&\le C\int_\R b_A \zeta_A^2 |u'|^2 dx
\le C \int_\R b_A |(\zeta_A u)'|^2dx + \frac{C}{A^2} \int_\R |\zeta_A u|^2dx\\
&\le C \int_\R |(\zeta_A u)'|^2dx + \frac{C}{A^2} \int_\R b_A^2|u|^2dx.
\end{align*}
Therefore, we obtain
\begin{align*}
\int_0^T\nor{b_Au'}{L^2}^2 + \frac{1}{A^2}\nor{b_Au}{L^2}^2 dt
&\le C\int_0^T \nor{(\zeta_A u)'}{L^2}^2 + \nor{b_{10}^2 u}{L^2}^2 dt\\
&\le C\sup_{t\in [0,T]}|\boI(t)| + C \int_0^T\nor{b_{10}^2 u}{L^2}^2 dt\\
&\le CA\del^2+ C \int_0^T\nor{b_{10}^2 u}{L^2}^2 dt,
\end{align*}
thereby completing the proof.

\section{Proof of asymptotic stability}
\label{S7}

\subsection{Bootstrapping}
In this section, we prove Theorem \ref{Tx1.2Jun} using the bootstrap framework introduced in \cite{MR4864876}. The central lemma is stated as follows:
% \dlt{Suppose $\eps>0$. 
% The key bootstrapping lemma that we shall prove is the following:}
%We shall prove Theorem \ref{Tx1.2Jun} by following the framework of Cuccagna and Maeda \cite{MR4864876}. The estimate proceeds by bootstrap argument:
%
%
\begin{prop}(Bootstrapping lemma)\label{P6.2Aug}
Suppose $\eps>0$. 
There exists $\nu_0>0$ such that for $0<\nu<\nu_0$, 
%Let $\nu$ be a given small parameter. 
there exist $\del(\nu)>0$, $A(\nu)>0$ satisfying the following property: 
Let $u$ be a solution to \eqref{E2.2Apr} with $u_1,u_2$ odd, $u_3$ even, and 
satisfying $\nor{u}{L^\infty_tH^1_x}\le \del$. 
Assume that
\begin{equation}\label{E5.1Apr}
\nor{b_{A}u'}{L^2_t(0,T\colon L^2_x)}
+\frac{1}{A}\nor{b_{A}u}{L^2_t(0,T\colon L^2_x)}
+\nor{b_{10}u}{L^2_t(0,T\colon L^2_x)} \le \nu
\end{equation}
for some $T>0$. Then we have
%claim that if $0<\del<\del_0(\nu)$, we have
\begin{equation}\label{eq:Aug22}
\nor{b_{A}u'}{L^2_t(0,T\colon L^2_x)}
+\frac{1}{A}\nor{b_{A}u}{L^2_t(0,T\colon L^2_x)}
+\nor{b_{10}u}{L^2_t(0,T\colon L^2_x)} \le o_{\nu\to 0}(\nu).
\end{equation}
\end{prop}
As verified later, it suffices to choose $\del=\del(\nu)$, $A=A(\nu)$ satisfying
\begin{align}\label{eq:Aug21}
\frac 1A = o_{\nu\to 0} (1),&&
A^{1/2}\del =  o_{\nu\to0}(\nu),
&&A^2\del = o_{\nu\to0}(1).
\end{align}
Assuming Proposition \ref{P6.2Aug} for now, we first complete the proof of Theorem \ref{Tx1.2Jun}.
\begin{proof}[Proof of Theorem \ref{Tx1.2Jun}]
%[Proof of Theorem \ref{Tx1.2Jun}]
Let $\nu>0$ be a sufficiently small number such that the RHS of \eqref{eq:Aug22} can be bounded by $\frac 12\nu$. 
Let $\del_0>0$ be a number fixed later. 
Let $A=A(\nu)$, 
$\del=\del(\nu)$ be numbers as in \eqref{eq:Aug21}. 
% At the same time, we choose $\del_0(\nu)$ sufficiently small such that the argument in Section 1 is also true. By Section 1, 
By Theorem \ref{T2.1}, if $\nor{\bom(0,\cdot)-\bom_0}{H^1_x}\le \del_0$, then 
the solution $\bom$ to \eqref{E5.1Mar} is global-in-time, and satisfies $\nor{\bom(t)-\bom_0}{L^\infty_tH^1_x}\le C\del_0$. 
In particular, by choosing $\del_0\ll 1$, it holds that $\nor{u}{L^\infty_tH^1_x} <\del(\nu)$. 
Therefore, the bootstrap argument works, and we obtain
\[
\nor{b_{A}u'}{L^2_t(0,\infty\colon L^2_x)}
+\frac{1}{A}\nor{b_{A}u}{L^2_t(0,\infty\colon L^2_x)}
+\nor{b_{10}u}{L^2_t(0,\infty\colon L^2_x)} \le \nu.
\]
Hence, we have
\begin{align*}
&\int_0^\infty \nor{b_A\rd_x(\bom(t)-\bom_0)}{L^2}^2 + \nor{b_A(\bom(t)-\bom_0)}{L^2_x}^2 dt\\
&\lesssim \int_0^\infty \nor{b_Au_1'}{L^2}^2 + \nor{b_A(u_2'+bu_3)}{L^2}^2
 + \nor{b_A(u_3'-bu_2)}{L^2}^2 
 + \nor{b_A u}{L^2}^2 + \nor{b_Au_3^2}{L^2_x}^2 dt\\
 &\lesssim \int_0^\infty \nor{b_A u'}{L^2_x}^2 + \nor{b_Au}{L^2_x}^2 dt \lesssim 1
\end{align*}
as desired.
\end{proof}

\subsection{Proof of Proposition \ref{P6.2Aug}}
\subsubsection{Reduction}
We begin by reducing our goal \eqref{eq:Aug22} in several steps. 
If $\nu\ll 1$, then Proposition \ref{P4:4.1} can be applied, yielding
\[
\nor{b_{A}u'}{L^2_t(0,T\colon L^2_x)}
+\frac{1}{A}\nor{b_{A}u}{L^2_t(0,T\colon L^2_x)} \le C A^{1/2}\del + C\nor{b_{10}u}{L^2_t(0,T\colon L^2_x)},
\]
and thus it suffices to show
\begin{align}\label{E5.2Apr}
\nor{b_{10}u}{L^2_t(0,T\colon L^2)} \le o(\nu).
\end{align}
Now, we introduce a new parameter $B=A^{1/3}$. Note that $A \gg B  \gg 1.$
Let $\chi$ be an even cut-off function with
$\chi\equiv 1$ on $[-1,1]$ and $\chi=0$ on $\R\setminus (-2,2)$. 
Define $\chi_B = \chi(B^{-1}\cdot)$. 
To show \eqref{E5.2Apr}, it suffices to show
\begin{align}\label{E5.3Apr}
\nor{\chi_Bu}{L^2_t(0,T\colon L^{2,-s})} \le o(\nu),\qquad s>\frac 32,
\end{align}
since
\begin{align*}
\nor{b_{10}u}{L^2_t(0,T\colon L^2)} &\le \nor{\chi_B \jb{x}^s b_{10}u}{L^2_tL^{2,-s}_x} + \nor{(1-\chi_B) b_{10}u}{L^2_tL^{2}_x}\\
&\lesssim \nor{\chi_B u}{L^2_tL^{2,-s}_x} + e^{-(\frac{1}{10}-\frac{1}{A})B}\nor{b_Au}{L^2_tL^2_x}
%(\because\ \jb{x}^{s}b_{10} \lesssim 1)
\\
&\le  o(\nu) + Ae^{-(\frac{1}{10}-\frac{1}{A})B} \nu
\le o(\nu)
\end{align*}
since 
$Ae^{-(\frac{1}{10}-\frac{1}{A})B}\le Ae^{-\frac{B}{20}}\le \frac{1}{A}= o(1)$.
\subsubsection{Stereographic projection and localization}

We shall obtain \eqref{E5.3Apr} by applying linear estimates to the $(u_1,u_2)$ system. Within the framework of \cite{MR4864876}, second-order derivative nonlinearities pose significant challenges, rendering the original system \eqref{E2.2Apr} unsuitable for direct estimation. To circumvent this, we again employ the semilinearized equation \eqref{eq:Aug3}, which is free of these second-order terms.\par
Let $v \coloneqq (v_1,v_2)= {}^t (\frac{u_1}{2+u_3}, \frac{u_2}{2+u_3})$. Then $v$ satisfies
\begin{align}\label{E5.5Apr2}
\rd_t v = \boA_\eps v +\tilde{N}(u),\qquad \tilde{N}={}^t(\tilde{N_1},\tilde{N}_2),
\end{align}
where $\tilde{N}_j$, $j=1,2$ is as defined in \eqref{E5.5Apr}. Multiplying \eqref{E5.5Apr2} by $\chi_B$, we obtain
\begin{align}\label{E5.6Apr}
\begin{aligned}
\rd_t (\chi_B v) 
&= \boA_\eps (\chi_Bv) + [\chi_B,\boA_\eps] v + \chi_B \tilde{N}(u)\\
&= \boA_\eps (\chi_B v) + J(\chi_B'' v +2\chi_B' \rd_xv)  + \chi_B \tilde{N}(u)
\end{aligned}
\end{align}
where 
$J=\begin{pmatrix}
    0 & -1\\
    1 & 0
\end{pmatrix}$. 
Here, 
by virtue of the spatial symmetry of $u$, 
both sides of \eqref{E5.6Apr} are odd, and in particular included in the range of $P_c(\boA_\eps)$. 
Hence we can write
%\begin{screen}
\begin{align*}
\chi_Bv = e^{t\boA_\eps} P_c \chi_B v(0) 
&+ 
\int_0^t e^{(t-s)\boA_\eps} P_c (J (\chi_B'' v +2\chi_B' \rd_xv)) ds\\
&+ \int_0^t e^{(t-s)\boA_\eps} P_c \chi_B \tilde{N}(u) ds.
\end{align*}
%\end{screen}
%We estimate the $L^2_tL^{2,-s}_x$-norm:
Since $\nor{u}{L^\infty_TL^\infty_x}\le C\del\ll 1$, $(2+u_3)^{-1}$  admits a uniform bound by a constant. 
Hence, \eqref{eq:Aug5} implies
\begin{align*}
\nor{\chi_Bu}{L^2_tL^{2,-s}_x} 
&\lesssim
\nor{\chi_Bv}{L^2_tL^{2,-s}_x} \\
&\lesssim \nor{e^{t\boA_\eps} P_c \chi_B v(0)}{L^2_tL^{2,-s}_x} 
+ 
\nor{\int_0^t e^{(t-s)\boA_\eps} P_c (\si_3 (\chi_B'' v +2\chi_B' \rd_xv)) ds}{L^2_tL^{2,-s}_x}\\
&\qquad + \nor{\int_0^t e^{(t-s)\boA_\eps} P_c \chi_B \tilde{N}(u) ds}{L^2_tL^{2,-s}_x}.
\end{align*} 
Among these, we begin by observing the first two terms. 
%the first term can be estimated from 
By Proposition \ref{P4:3.3}, we have
\[
\nor{e^{t\boA_\eps} P_c \chi_B v(0)}{L^2_tL^{2,-s}_x} \le C \nor{\chi_Bv(0)}{L^2} \le C\del \lesssim o(\nu) ,\qquad \text{if } \del(\nu) \lesssim o(\nu).
\]
For the second term, we borrow the estimate in \cite{MR4864876}*{Lemma 7.2}. Let us take $\tau\in (\frac 12,1)$. 
Then Proposition \ref{P4:3.2} implies
\begin{align*}
&\nor{\int_0^t e^{(t-s)\boA_\eps} P_c (J(\chi_B'' v +2\chi_B' \rd_xv)) ds}{L^2_tL^{2,-s}_x} \\
&\lesssim 
\nor{\chi_B'' v +2\chi_B' \rd_x v}{L^2_tL^{2,\tau}_x}
\lesssim 
\nor{\chi_B'' u +2\chi_B' \rd_x u}{L^2_tL^{2,\tau}_x}
\\ &\lesssim B^{\tau-1} (\nor{b_A u'}{L^2_{t,x}}+ \frac 1A \nor{b_A u}{L^2_{t,x}} + \nor{b_{10} u}{L^2_{t,x}}) \lesssim B^{\tau-1} \nu\lesssim o(\nu).
\end{align*}
%the remaining task is to obtain the bound the following bound
%\begin{itembox}[l]{\textbf{Desired estimate}}
%\end{itembox}
%
\subsubsection{Nonlinear estimate}
It remains to show
\begin{equation}\label{E5.4Apr}
\nor{\int_0^t e^{(t-t')\boA_\eps} P_c \chi_B \tilde{N}(u) dt'}{L^2_tL^{2,-s}_x} \lesssim o(\nu),\qquad \text{for some }s>\frac 32.
\end{equation}
%We shall prove \eqref{E5.4Apr}. 
We begin by classifying the nonlinearity. 
Noting that $\nor{u}{L^\infty_{t,x}}\le C\del \ll 1$ and \eqref{eq:Aug5}, 
% The constraint $u_1^2+u_2^2+(1+u_3)^2=1$ implies that
% \[
% u_3 = -\frac{u_1^2+u_2^2}{\sqrt{1-u_1^2-u_2^2}},\qquad u_3' \sim u  u'.
% \]
% Now let $u= {}^t (u_1,u_2)$. 
% Noting that $\nor{u}{L^\infty_{t,x}}\le C\del \ll 1$, we may estimate the order as
% \[
% u_3 \sim O(u^2),\qquad u_3' \sim O(u u').
% \]
% Accordingly, we may ignore $2+u_3$ as it is bounded by constant. 
%Therefore, 
the nonlinearity essentially consists of
\[
O(|u|^2),\qquad O(|u| |u'|),\qquad 
O(|u||u'|^2).
\]
For the first two cases, we apply the inhomogeneous Kato smoothing (Proposition \ref{P4:3.2}), yielding
\begin{align*}
&\nor{\int_0^t e^{(t-t')\boA_\eps} P_c \chi_B O(|u(t')|^2)dt'}{L^2_tL^{2,-s}_x}\\
&\lesssim \nor{\chi_B u^2}{L^2_{t}L^{2,\tau}_x}\le 
\nor{u}{L^\infty_tL^\infty_x} \nor{\chi_B \jb{x}^{\tau}u}{L^2_{t,x}}
\lesssim \del B^\tau e^{\frac{2B}{A}} \nor{b_A u}{L^2_{t,x}} \le \del AB^\tau \nu,\\
&\nor{\int_0^t e^{(t-t')\boA_\eps} P_c \chi_B O(|u(t')||u'(t')|)dt'}{L^2_tL^{2,-s}_x}\\
&\lesssim \nor{\chi_B uu'}{L^2_{t}L^{2,\tau}_x}\le 
\nor{u}{L^\infty_tL^\infty_x} \nor{\chi_B \jb{x}^{\tau}u'}{L^2_{t,x}}
\lesssim \del B^\tau e^{\frac{2B}{A}} \nor{b_A u'}{L^2_{t,x}} \le \del B^\tau \nu.
\end{align*}
Here, \eqref{eq:Aug21} implies $\del AB^\tau \lesssim \del A^2\lesssim o_{\nu\to 0}(1)$, and thus we obtain
\begin{align*}
&\nor{\int_0^t e^{(t-t')\boA_\eps} P_c \chi_B O(|u(t')|^2)dt'}{L^2_tL^{2,-s}_x}\\
&+
\nor{\int_0^t e^{(t-t')\boA_\eps} P_c \chi_B O(|u(t')||u'(t')|)dt'}{L^2_tL^{2,-s}_x}
\lesssim o_{\nu\to 0} (\nu).
\end{align*}
For the last nonlinearity, we apply dual Strichartz estimate (Proposition \ref{P4:3.5}), yielding
\begin{align}\label{eq:Sep10}
\begin{aligned}
\nor{
\int_0^t e^{i(t-t')\boA_\eps}\chi_B O(|u(t')||u'(t')|^2] dt'}{L^2_tL^{2,-s}_x} 
&\le 
C \nor{\chi_B (u')^2u}{L^{\frac 43}_tL^1_x} \\
&\le C\nor{u}{L^\infty_TL^\infty_x} \nor{\chi_{2B}u'}{L^{\frac 83}_tL^2_x}^2 \\
&\le C\nor{u}{L^\infty_TL^\infty_x} (\nor{\chi_{2B}u'}{L^{2}_tL^2_x}^{\frac 34} 
\nor{\chi_{2B}u'}{L^{\infty}_tL^2_x}^{\frac 14})^2\\
&\le C \del^{\frac 32} \nu^{\frac 32}= o(\nu).
\end{aligned}
\end{align}
Therefore, we conclude \eqref{E5.4Apr}, and hence \eqref{E5.3Apr}.

\appendix
\titleformat{\section}[block]
  {\normalfont\Large\bfseries}
  {Appendix \thesection}
  {1em}
  {}
\section{Spectral analysis: Proof of Proposition \ref{Tx1}}\label{SA}
In this section, we prove Proposition \ref{Tx1}. Recall
\[
\boA_\eps = 
\begin{pmatrix}
    0 & -L_0\\
    L_\eps & 0
\end{pmatrix}
,\qquad 
L_\eps = -\rd_x^2 +1-2b^2+2\eps b.
\]
$\boA_\eps$ is structurally analogeous to the linearized operator around standing waves for NLS. 
The study of the spectral stability for such operators is well established, and Proposition \ref{Tx1} can be proven by combining known techniques in this context.\par 
\subsection{Basic properties}
We first review the basic properties of $L_0$, $L_\eps$ and $\boA_\eps$. 
The terminologies in the operator theory are based on the style of \cite{MR1335452}. 
We begin with the properties of the Schr\"odinger operators $L_\eps$, $\eps\in\R$. 
$L_\eps$ 
%is the Schr\"odinger operators; in particular 
is self-adjoint on $L^2(\R)$ with domain $\boD(L_\eps)=H^2(\R)$. 
More specifically, the potential of $L_0$ coincides with the P\"oschl--Teller potential, arising in the cubic nonlinear Schr\"odinger equations. Note that $L_0= F^*F$ with $F\coloneqq b \rd_x b^{-1}$, and thus $L_0$ is nonnegative. 
In more detail, we have
\begin{align}\label{eq:m1}
\si(L_0)= \{0\}\cup [1,\infty),\qquad \Ker L_0 = \Span \{b\}.
\end{align}
We also have $L_0p =p$, and hence $L_0$ has a threshold resonance at $1$. For $L_\eps$, Weyl's essential spectrum theorem (see Theorem IV.5.26 in \cite{MR1335452}) implies that $\si_{\ess} (L_\eps) = [1,\infty)$. 
Since $\eps\mapsto L_\eps$ is an analytic perturbation, the simple eigenvalue $0$ for $L_0$ bifurcates as $\eps$ is perturbed. In detail, we have the following:
%to $\la_\eps$ for $L_\eps$ near $\eps\sim0$, which is increasing in $\eps$ since so is $L_\eps$. However, the 
\begin{lem}\label{L:AugA5}
%We have $\si_\ess(L_\eps)=[1,\infty)$ for all $\eps\in\R$. Moreover, 
There exists $\eps_0>0$ such that
\begin{align*}
\si(L_\eps) = \{\la_\eps\}\cup [1,\infty) \quad\text{if } 0\le\eps<\eps_0,&&
\si(L_\eps) = [1,\infty) \quad\text{if } \eps\ge\eps_0
\end{align*}
with $\la_\eps \in [0,1)$. 
Furthermore, $\la_\eps$ is simple and strictly increasing in $\eps$ as long as $0\le \eps <\eps_0$.
\end{lem}
\begin{proof}
%\begin{myleftbar}
Since $L_{\eps} < L_{\eps'}$ for $\eps<\eps'$, 
the min-max principle, together with \eqref{eq:m1}, implies that $L_\eps$ has at most one point spectrum in $(0,1)$ if $\eps>0$. 
% On the other hand, 
% the standard perturbation theory implies that if $L_{\eps_1}$ has a simple eigenvalue $\la_{\eps_1}$ in $(-\infty,1)$ for some $\eps_1$, then so does $L_\eps$ for $\eps$ close to $\eps_1$, and the corresponding eigenvalue $\la_\eps$ is strictly increasing in $\eps$. 
% there exist $\del>0$ and a continuous function $\la(\eps)$ on $\{|\eps-\eps_1|<\del_1\}$ such that $\la(\eps_1)=\la_{\eps_1}$ and $\la(\eps)$ is a simple eigenvalue of $L_\eps$ for $|\eps-\eps_1|<\del$. 
% Moreover, since $L_{\eps} < L_{\eps'}$ if $\eps<\eps'$, $\la(\eps)$ is strictly increasing.\par
In light of this, define
\[
\la_\eps \coloneqq \inf_{f\in H^1,\ \nor{f}{L^2}=1} \inp{L_\eps f}{f}_{L^2},\qquad 
% \]
% \[
\eps_0 \coloneqq \inf \{\eps\in\R \mid  \la_\eps \ge1\}.
\]
Obviously, $\la_\eps$ is non-decreasing in $\eps$. 
If $0<\eps<\eps_0$, then $\la_\eps$ is the unique point spectrum of $L_\eps$ by the min-max principle, 
and strictly increasing in light of analytic perturbation theory. 
Since $L_1 = -\rd_x^2 +1 -2b^2 + 2b \ge 1$, it follows that $\eps_0\le 1$. We also have $\la_{\eps_0}=1$; otherwise, $\la_{\eps_0}<1$ would contradicts the definition of $\eps_0$ and the analytic perturbation theory. 
%\end{myleftbar}
\end{proof}
Next, we observe the basic properties of $\boA_\eps$. This is a closed operator on $L^2(\R\colon \C)^2$ with domain $\boD(\boA_\eps)=H^2(\R)^2$. Since $\boA_\eps$ is a relatively compact perturbation of $\boA_*\coloneqq 
\begin{pmatrix}
0 & -(-\rd_x^2 +1) \\
-\rd_x^2 +1 & 0
\end{pmatrix}$, we have 
%\[
$\si_{\ess} (\boA_\eps) = i(-\infty,-1] \cup i[1,\infty)$. 
%\]
Moreover, for $\la\not \in \si_\ess (\boA_\eps)$, 
%$\la\in \si(\boA_\eps) \setminus \si_{\ess \boA_\eps}$,
the stability theorem \cite{MR1335452}*{Theorem IV.5.26} implies that 
%\[
$\ind (\boA_\eps- \la) = \ind (\boA_*-\la) =0$. 
%\]
Hence, every $\la\in \si(\boA_\eps) \setminus \si_{\ess } (\boA_\eps)$ is an eigenvalue of $\boA_\eps$. Next,
$\boA_\eps$ can be rewritten as 
%\[
$\boA_\eps = 
\begin{pmatrix}
0 & -1\\
1 & 0
\end{pmatrix}
\begin{pmatrix}
L_\eps & 0\\
0 & L_0
\end{pmatrix}
$, 
%\]
which takes the form considered in the Grillakis--Shatah--Strauss theory \cite{MR1081647}. 
By \cite{MR1081647}*{Lemma 5.6}, $\si(\boA_\eps)$ is symmetric with respect to both real and imaginary axes. Furthermore, in light of \eqref{eq:m1}, it follows from \cite{MR1081647}*{Lemma 5.8} that 
\begin{align}
&\dim 
\Span\{ \phi\in (H^2)^2 \mid \phi\colon \text{eigenfunction for } \boA_\eps \text{ with eigenvalue } \la \text{ satisfying } \Re \la>0  \}\notag\\
&\le \dim 
\Span\{ f\in H^2 \mid f\colon \text{eigenfunction for } L_\eps \text{ with negative eigenvalue} \}\label{eq:m2}.
\end{align}

\subsection{Criterion}
Following \cite{MR2368894}, we define
\begin{align}\label{eq:m3}
\mu_\eps := 
\inf_{\substack{f\in H^1(\R)\\ f\perp b,\ \nor{f}{L^2} =1}} 
\inp{L_\eps f}{f}_{L^2(\R)}.
\end{align}
Then, we claim that $\mu_\eps$ gives the criterion for the linear stability.
%\begin{screen}
\begin{prop}\label{P1.1}
The following are true:\\
(i) If $\mu_\eps \ge 0$, then $\si(\boA_\eps) \subset i\R$.\\
(ii) If $\mu_\eps <0$, then there exist eigenvalues $\pm \la$ with $\la>0$.
\end{prop}
\begin{proof}
(i) Suppose $\mu_\eps\ge 0$. We first claim that the dimension of negative eigenspace for $L_\eps$ is at most $1$. 
Suppose that there is two linearly-independent negative eigenfunction $\phi_1,\phi_2$ of $L_\eps$, 
with eigenvalues $\la_1,\la_2<0$, respectively. 
Since $L_\eps$ is self-adjoint, we may suppose that $\phi_1\perp \phi_2$ in $L^2(\R^2)$, and $\nor{\phi_1}{L^2}=\nor{\phi_2}{L^2}=1$. 
If $\inp{\phi_j}{b}_{L^2}=0$ for either $j=1,2$, then 
\[
\inp{L_\eps \phi_j}{\phi_j}_{L^2} = \la_j \nor{\phi_j}{L^2}^2 <0,
\]
which contradicts $\mu_\eps \ge 0$. Hence $\inp{\phi_j}{b}_{L^2}\neq 0$ for $j=1,2$. 
Here, define $\phi := -\inp{b}{\phi_2} \phi_1 + \inp{b}{\phi_1} \phi_2$.
Then $\phi\neq 0$ and satisfies $\inp{b}{\phi}=0$. Moreover,
\begin{align*}
\inp{L_\eps \phi}{\phi} 
&= 
\inp{-\la_1 \inp{b}{\phi_2} \phi_1 + \la_2 \inp{b}{\phi_1} \phi_2}{
-\inp{b}{\phi_2} \phi_1 + \inp{b}{\phi_1} \phi_2
} \\
&= \la_1 |\inp{b}{\phi_2}|^2 + \la_2 |\inp{b}{\phi_1}|^2 < 0,
\end{align*}
which contradicts $\mu_\eps\ge 0$.\par
Assuming that there exists $\la_0\in \si(\boA_\eps) \setminus i\R$, we will prove contradiction. 
By the symmetry of $\si(\boA_\eps)$, we may assume $\Re \la_0>0$. 
Then by \eqref{eq:m2}, the number of eigenvalues with positive real part is at most $1$, and thus it should hold that $\la_0>0$. 
Then, there exists ${}^t(f,g)\in H^2(\R)^2\setminus \{0\}$ such that 
\[
\boA_\eps 
\begin{bmatrix}
f \\ g
\end{bmatrix}
= \la_0 
\begin{bmatrix}
f \\ g
\end{bmatrix}
.
\]
Note that $f\neq 0$; otherwise $\la_0g=L_0f=0$, which contradicts ${}^t(f,g)\neq0$. 
Then, we have
\[
\boA_\eps^2 
\begin{bmatrix}
f \\ g
\end{bmatrix}
=
\begin{bmatrix}
- L_0 L_\eps f \\ 
- L_\eps L_0 g
\end{bmatrix}
=
\la_0^2
\begin{bmatrix}
f \\ g
\end{bmatrix}
.
\]
In particular, we have $f\in R(L_0) = \ovl{\ker L_0}^\perp 
% = \ovl{\Span \{b\}}^\perp
= 
(\Span \{b\})^\perp$. 
Here, we note that $L_0\ge 1$ on $(\Span \{b\})^\perp$. 
If we let $h:= (L_0)^{-1/2} f$, we have 
$(L_0)^{1/2} L_\eps (L_0)^{1/2} h = -\la_0^2 h$, which gives $\inp{L_\eps f}{f}_{L^2} =-\la_0^2 \nor{h}{L^2}^2 <0$. 
% \[
% (L_0)^{1/2} L_\eps (L_0)^{1/2} h
% = (L_0)^{1/2} L_\eps f
% = -\la_0^2 (L_0)^{-1/2} f = -\la_0^2 h,
% \]
% which gives
% \[
% \inp{L_\eps f}{f}_{L^2} =
% \inp{L_\eps (L_0)^{1/2} h}{ (L_0)^{1/2}h}_{L^2} =
% -\la_0^2 \nor{h}{L^2}^2 <0.
% \]
Since $f\perp b$, it follows that $\mu_\eps \le \inp{L_\eps f}{ f}_{L^2} <0$, which is a contradiction.\par
(ii) The proof is essentially the same as \cite{MR2368894}, and thus we only sketch the outline. See also \cite{MR4899905}*{Section 2.1.2} for more explicit argument. Set
\begin{equation}\label{xa1.2}
\nu_\eps =
\inf_{\substack{f\in H^1\setminus \{0\} \\ \inp{f}{b}_{L^2} =0 }}
\frac{
\inp{L_\eps f}{f}_{L^2}
}{
\inp{(L_0)^{-1} f}{f}_{L^2}
}.
\end{equation}
Then, thanks to the assumption $\mu_\eps<0$, 
the direct method shows that \eqref{xa1.2} has a minimizer $f\in H^1\setminus \{0\}$ with $f\perp b$. 
By the Lagrange multiplier theorem, $f$ satisfies 
%\[
$L_\eps f = \al (L_0)^{-1} f + \be b$.
%\]
for some $\al,\be\in\R$. 
Taking inner product with $f$, we have $\al= \nu_\eps$. 
Letting $e_\eps:= \sqrt{-\nu_\eps} >0$, we have
\[
\begin{pmatrix}
0 & -L_0 \\
L_\eps & 0
\end{pmatrix}
\begin{bmatrix}
f\\
-e_\eps (L_0)^{-1} f + e_\eps^{-1} \be b
\end{bmatrix}
=
e_\eps 
\begin{bmatrix}
f\\
-e_\eps (L_0)^{-1} f + e_\eps^{-1} \be b
\end{bmatrix}
.
\]
Thus $e_\eps$ is an eigenvalue of $\boA_\eps$. By the symmetry of $\si(\boA_\eps)$, so is $-e_\eps$.
\end{proof}

%\end{screen}
%
%
\subsection{Sign of $\mu_\eps$}
Theorem \ref{Tx1} is an immediate consequence of the following claim:
%\begin{screen}
\begin{prop}\label{P1.2}
$\mu_\eps$ is finite-valued. Moreover, the map $\eps\mapsto \mu_\eps$ is \\
(i) non-decreasing;\\
(ii) continuous; \\
(iii) If $\eps\ge -\frac 12$, we have $\mu_\eps\ge 0$;\\
(iv) $\lim_{\eps\to - \infty} \mu_\eps =- \infty$.\\ 
In particular, 
by setting $\eps_* \coloneqq \inf\{\eps \in\R \mid \mu_\eps\ge 0\}$, 
%there is a unique $\eps=\eps_*$ such that 
then $\mu_\eps\ge 0$ if $\eps \ge \eps_*$, and $\mu_\eps<0$ if $\eps<\eps_*$.
\end{prop}
%\end{screen}

%Although 
% It is straightforward that 
% Proposition \ref{P1.2} is enough to conclude Proposition \ref{Tx1}. 
For the sake of the use in the proof of orbital stability, 
we further observe the following property. 
%, which will be used later.
%the positivity of $\mu_\eps$
\begin{prop}\label{P1.3}
If $\eps>\eps_*$, then $\mu_\eps >0$.
\end{prop}

%\subsection{Proof of Proposition \ref{P1.2}}
\begin{proof}[Proof of Proposition \ref{P1.2}]
We only give a proof of (iii), while the others are straightforward. 
Let $f\in H^1(\R)$ with $\inp{f}{b}_{L^2} =0$. 
Since $L_0 \ge 1$ on $(\Span \{ b\})^\perp$, we have
\[
\inp{L_\eps f}{f} = \inp{L_0 f}{f} + 2\eps \inp{b f}{f} 
\ge \nor{f}{L^2}^2 + 2\eps \inp{b f}{f} 
= \int_{\R} (1+2\eps b) |f|^2 dx\ge 0 
\]
by noting $1 + 2\eps b \ge 0$. 
% Since $0 < b \le 1$, 
% when $\eps \ge -\frac 12$ we have $1 + 2\eps b \ge 0$. Hence %the claim follows.
% \[
% \inp{L_\eps f}{f}_{L^2} \ge 0.
% \]
By taking infimum, we have $\mu_\eps \ge 0$.
\end{proof}

%\subsection{Proof of Proposition \ref{P1.3}}\label{SA.5}
%
\begin{proof}[Proof of Proposition \ref{P1.3}]
First note that the continuity of $\eps\mapsto \mu_\eps$ implies that $\mu_{\eps_*}=0$. Suppose that $\mu_{\eps_0}=0$ for some $\eps_0>\eps_*$. By a standard argument, it follows that the infimum \eqref{eq:m3} is attained for $\eps=\eps_0$. Letting $f_*$ be a minimizer, we have
\[
0\le \inp{L_{\eps_*} f_*}{f_*} = \inp{L_{\eps_0} f_*}{f_*} + 2(\eps_*-\eps_0) \int_\R b |f_*|^2 dx 
< 0,
\]
which is a contradiction.
\end{proof}

\subsection{Numerics}\label{SA.6}

\begin{figure}
    \includegraphics[scale=0.5]{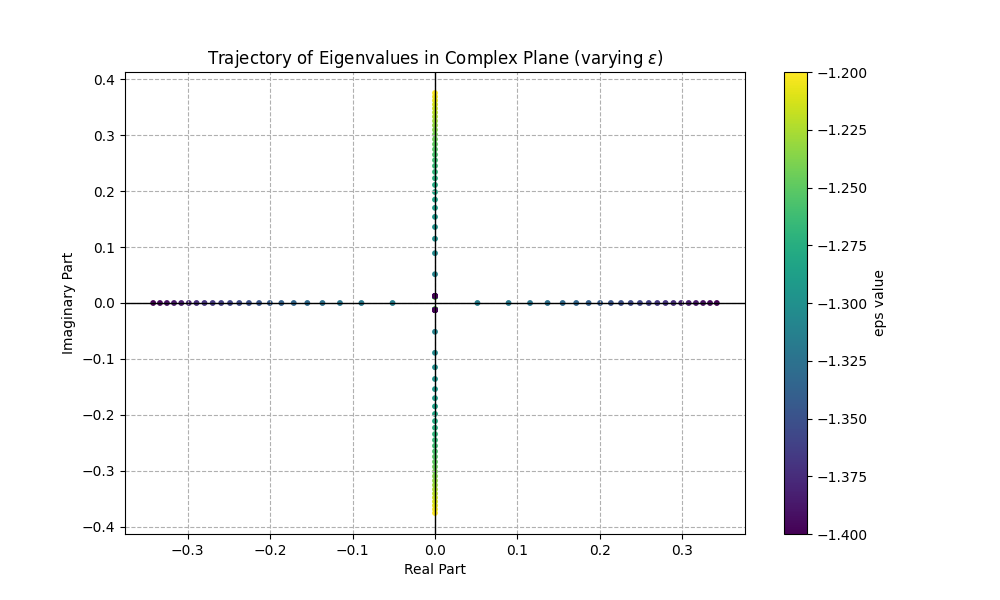}
    \caption{The numerical plot of the eigenvalues of $\boA_\eps$ near the origin.}
    \centering
    \label{F1}
\end{figure}

Figure \ref{F1} displays the four eigenvalues of $\boA_\eps$ closest to the origin, computed via a standard discretization of the eigenvalue problem using the scheme from \cite{MR2368894}*{Appendix} and the ARPACK library in \texttt{SciPy.sparse}. The two eigenvalues at the origin correspond to the kernel generated by translation invariance. As $\eps$ decreases, a pair of imaginary eigenvalues collides at the origin at a critical parameter $\eps=\eps_*$ before bifurcating onto the real axis. This critical value corresponds to the spectral stability threshold in Proposition \ref{Tx1}. From Figure \ref{F1}, we estimate $\eps_* \approx -1.3$.

\subsubsection*{Numerical method}

\begin{itemize}
\item The domain: $I=[-T_{\max},T_{\max}]$ with $T_{\max}=10$
\item Mesh width: $\del=0.04$,\qquad $N=\frac{T_{\max}}{\del}=250$
\item For $u\colon I\to \C$, define $(u_j)_{-N\le j\le N}$ with $u_j \coloneq u(\del j)$
\item Discretize $\rd_{xx} u$ as
\[
\rd_{xx} u \mapsto \frac{u_{j+1} -2u_j +u_{j-1}}{\del^2},\qquad 
u_{-N-1}=u_{N+1} =0\quad (\text{Dirichlet B.C.}).
\]
Hence the operator $\rd_{xx}$ turns to be the matrix $A$ with
\[
A = 
\frac{1}{\del^2}
\begin{bmatrix}
-2 & 1 & 0 & \cdots &0&0\\
1 & -2 & 1 & \cdots &0&0\\
0 & 1 & 2 & \cdots & 0&0\\
\vdots & & & \ddots & 0&0\\
0 & 0 & 0 &\cdots &2&1\\
0 & 0 & 0 &\cdots &1 &2
\end{bmatrix}
\]
\item Let $f^{\eps} (x) = 1-2b^2(x) + 2\eps b(x)$, and set $[f^\eps] \coloneq \mathrm{diag} (f^\eps (\del j))$. 
\item Consequently, the eigenvalue problem for $\boA_\eps$ turns out into
\[
\begin{bmatrix}
O & A - [f^0 ]\\
-A +[f^\eps] & O 
\end{bmatrix}
\begin{bmatrix}
u_{-N} \\ \vdots \\ u_{N}
\end{bmatrix}
=
\la 
\begin{bmatrix}
u_{-N} \\ \vdots \\ u_{N}
\end{bmatrix}
\]
\end{itemize}

\section{Spectral properties II: Internal modes, threshold resonances, embedded eigenvalues}\label{SB}
%
%
%\subsection{Background}
%
%
In this section, we provide a proof of Proposition \ref{P3.2Jun}. 
Before the proof, we introduce the notions related to the statement. 
%As observed in Section \ref{SA}, we have $\si_\ess(\boA_\eps) = i(-\infty,-1]\cup i[1,\infty)$. 
%of internal modes, threshold resonances, and embedded eigenvalues. 
%In the study of linear Schro\"odinger theory, these are known as an obstacle to proving the smoothing estimates, such as the Strichartz estimate. 
% \par 
% To begin with, we recall that the essential spectrum of $\boA_\eps$ is $i(-\infty,-1]\cup i[1,\infty)$. These notions are define via the position relative to the essential spectrum.
%
%
\begin{definition}
(i) An \textit{internal mode} is a point spectrum $i\la$ of $\boA_\eps$ with $\la\in [-1,1]\setminus \{0\}$.\\
(ii) $\boA_\eps$ has a \textit{threshold resonance} at $i$ ($-i$, resp.) 
% is a \textit{threshold resonance} 
if and only if there exists $f\in L^\infty(\R)^2\setminus L^2(\R)^2$ such that $\boA_\eps f=if$ ($\boA_\eps f=-if$, resp.).\\
(iii) $\la$ is an \textit{embedded eigenvalue} of $\boA_\eps$ if and only if $\la\in i(-\infty,-1)\cup i(1,\infty)$ and there exists $f\in D(\boA_\eps)$ with $\boA_\eps f=\la f$.
\end{definition}
\begin{remark}
At the linear level, an internal mode $i\la$ yields a bound state of the linearized equation. Nonlinearly, however, this mode radiates its energy into the continuous spectrum and decays via dispersive effects, unless nonlinear interactions perfectly balance (e.g., \cites{MR1218303,MR4228957}). Establishing asymptotic stability therefore requires a complete understanding of the internal modes of $\boA_\eps$. Conversely, the absence of threshold resonances and embedded eigenvalues is sufficient for the smoothing estimates to hold (see \cites{MR2219305,MR4864876}). 
%Propositions \ref{P4:3.2} and \ref{P4:3.4}. 
Importantly, for $\eps=0$, $\boA_0$ exhibits threshold resonances at $\pm i$ with $f= \begin{pmatrix}
    p \\ \mp p
\end{pmatrix}$, explaining why this case is omitted from Proposition \ref{Tx1.2Jun}.

\end{remark}
%
% \begin{remark}
% We discuss the 
% \end{remark}

% \noindent{}\textbf{\underline{Background.}} 
% (i) At linear level, a internal mode $i\la$ produces a bound state of linearized equation, oscillating with period $\la$. However, at nonlinear level, this mode will radiate its energy to continuous spectrum via nonlinear interaction and end up decaying by dispersive effect of continuous spectrum, except quite rare case when nonlinear interaction perfectly balances. Accordingly, for the asymptotic stability, we need to know complete information of internal modes for $\boA_\eps$ to follow this process precisely. \par
% (ii) The threshold resonance is known as a technical obstacle to prove a dispersive estimates for $\boA_\eps$. 
% %for Schr\"odinger operators. 
% This notion also arises in the study of standard Schr\"odinger operators $H=-\Del +V$, and causes the same problem. 
% On the other hand, this is not an essential factor to make dispersive estimates fail; indeed, in the case of $H=-\Del+V$, 
% a stronger assumption would imply a dispersive estimate whether or not there is a threshold resonance. See for example \textcolor{red}{[Goldberg--Schlag 2004]}. \par 
% (iii) The embedded eigenvalues are also a technical obstacle for dispersive estimates. This notion does not arise in the case of standard Schr\"odinger operators $H=-\Del +V$, since the possibility of such a type of eigenvalues is a priori excluded. 
%the free Schr\"odinger operator $\Del$ 
%This process is called the Fermi golden rule.
%
%
\subsection{Jost solutions
%Preliminaries
}\label{SB1}
%
%
%For the proof of Proposition xxxx, we give some preparations.
%
We begin by introducing some preliminaries. 
Recall the definition of the operator:
\[
\boA_\eps = 
\begin{pmatrix}
0 & -L_0\\
L_\eps & 0
\end{pmatrix}
,\qquad 
L_\eps =-\rd_x^2 +1-2b^2 +2\eps b.
\]
We observe that
\begin{equation}\label{Eq:B.1Jun}
-i\boH_\eps = U^{-1}\boA_\eps U,\qquad U= \frac 1{\sqrt2}
\begin{pmatrix}
1 & 1\\
i & -i
\end{pmatrix}
,
\end{equation}
where
%Considering $\boA_\eps$ is equivalent to $-i\boH_\eps$:
\begin{gather*}
\boH_\eps = \si_3 H_\eps,\qquad 
\si_3= 
\begin{pmatrix}
1 & 0\\
0 & -1
\end{pmatrix}
,\qquad 
H_\eps \coloneqq (-\rd_x^2+1) I - 
\begin{pmatrix}
2b^2-\eps b & -\eps b\\
-\eps b & 2b^2-\eps b
\end{pmatrix}
.\\
\text{i.e.}\qquad 
\boH_\eps = 
\begin{pmatrix}
-\rd_x^2+1 + 2b^2-\eps b & -\eps b\\
-\eps b & \rd_x^2 -1 +2b^2 -\eps b
\end{pmatrix}
.
\end{gather*}
We set the potential matrix of $H_\eps$ as
\[
V \equiv V_\eps \coloneqq 
\begin{pmatrix}
2b^2 -\eps b & -\eps b\\
-\eps b & 2b^2-\eps b
\end{pmatrix}
.
\]
%\eqref{Eq:B.1Jun} offers an equivalent 
In view of the unitary equivalence \eqref{Eq:B.1Jun}, we will freely switch between the spectral problems for $\boA_\eps$ and $\boH_\eps$.\par 
%
%
%\subsubsection*{2. Jost solutions}
%
%
% Since $\boH$ satisfies all the assumptions in \cite{MR2219305}*{Definition 5.1}, we may employ the results in \cite{MR2219305}*{Section 5}. 
Next, we recall
%For $\la\in \R$, we recall 
the properties of the eigen-equation
\begin{equation}\label{Eq:B.2Jun}
\boH_\eps f=(\la^2+1)f,\qquad \la\in\R
\end{equation}
established in \cite{MR2219305}*{Section 5}. 
This knowledge will play important roles not only in the analysis of threshold resonance and embedded eigenvalues, but also the proof of the resolvent estimate for $\boH_\eps$ in Section \ref{SC}. 
In light of this, we record all the information that will be needed in the present paper.\par
In what follows, let $\gam\in (0,1)$ be a fixed number. (Any number in $(0,1)$ is fine.) 
We often abbreviate $\boH_\eps$ as $\boH$.
%
%
%\begin{itemize}
For $\la\in \R \setminus \{0\}$, 
the space of complex-valued solutions to \eqref{Eq:B.2Jun} is spanned by
%we consider the linear equation
% \[
% \boH f = (\la^2+1)f.
% \]
% The solution space is of dimension $4$. This space is spanned by
\begin{align}\label{Eq:B.3Jun}
f_1(x,\la),\quad f_2(x,\la),\quad f_3(x,\la),\quad f_4(x,\la)
\end{align}
with the following properties.
\begin{itemize}
\item 
There exists $x_0>0$ such that for $x\ge x_0$, we have
\begin{align*}
f_1(x,\la) = e^{i\la x} m_1(x,\la),\qquad 
|\rd_x^k\rd_\la^l(m_1(x,\la)-
\begin{pmatrix}
    1 \\0
\end{pmatrix}
)| \le C_{k,l}\jb{\la}^{-1} e^{-\gam x},\\
f_2(x,\la) = e^{-i\la x} m_2(x,\la),\qquad 
|\rd_x^k\rd_\la^l(m_2(x,\la)-
\begin{pmatrix}
    1 \\0
\end{pmatrix}
)| \le C_{k,l}\jb{\la}^{-1} e^{-\gam x},\\
f_3(x,\la) = e^{-\sqrt{2+\la^2}x} m_3(x,\la),\qquad 
|\rd_x^k\rd_\la^l(m_3(x,\la)-
\begin{pmatrix}
    0 \\1
\end{pmatrix}
)| \le C_{k,l}\jb{\la}^{-1} e^{-\gam x},\\
f_4(x,\la) = e^{\sqrt{2+\la^2}x} m_4(x,\la),\qquad 
|\rd_x^k\rd_\la^l(\tilde{m}_4(x,\la)-
\begin{pmatrix}
    0 \\1
\end{pmatrix}
)| \le C_{k,l}\jb{\la}^{-1} e^{-\gam x}
\end{align*}
for all $k,l\in\Z_{\ge 0}$. 
%%%
%%% Intuitional descriptions
%%%
% \textcolor{blue}{The functions} $f_j$ are called \textit{the Jost solutions}. 
% More intuitionally, \dlt{$f_j$'s} {\color{blue}they} have asymptotics like
% \begin{align*}
% f_1(x,\la) &= 
% \begin{pmatrix}
%     e^{ix\la} \\ 0
% \end{pmatrix}
% + O(e^{-\gam x}),\qquad x\to\infty,
% %+
% %O(\frac{e^{-\gam x}}{\sqrt{\la^2+2}}).
% \\
% f_2(x,\la) &= \ovl{f_1(x,\la)}
% =
% \begin{pmatrix}
%     e^{-ix\la} \\ 0
% \end{pmatrix}
% + O(e^{-\gam x}),\qquad x\to\infty,\\
% f_3(x,\la) &= 
% \begin{pmatrix}
% 0 \\ e^{-\sqrt{\la^2+2} x}
% \end{pmatrix}
% + O(e^{-(\sqrt{\la^2+2}+\gam )x}),\qquad x\to\infty,\\
% f_4(x,\la) &= 
% \begin{pmatrix}
% 0 \\ e^{\sqrt{\la^2+2} x}
% \end{pmatrix}
% + O(e^{(\sqrt{\la^2+2}-\gam )x}),\qquad x\to\infty.
% \end{align*}
% Note that the derivatives of the remainders also \dlt{exponentially decays} \textcolor{blue}{decay exponentially} as $x\to\infty$. 
%%%
%%% Intuitional descriptions end
%%%
%These solutions further satisfies
Furthermore, they satisfy
\begin{align*}
f_1(x,-\la) = \ovl{f_1(x,\la)} = f_2(x,\la),&&
f_2(x,-\la) = \ovl{f_2(x,\la)} = f_1(x,\la),\\
f_3(x,-\la) = \ovl{f_3(x,\la)} = f_3(x,\la),&&
f_4(x,-\la) = \ovl{f_4(x,\la)} = f_4(x,\la).
\end{align*} 
\item For $f,g:\R\to\C^2$, we define the Wronskian
\[
W[f,g] = {}^t(f')g-{}^tfg'.
\]
Then we have
%Regarding the Wronskian, we have
\begin{gather*}
W[f_1,f_2]=2i\la,\qquad
W[f_3,f_4]=-2\sqrt{2+\la^2},\\
W[f_1,f_3]=W[f_1,f_4]=W[f_2,f_3]=W[f_2,f_4]=0.
\end{gather*}
\item Since $V_\eps$ is even, 
the Jost solutions for $x\to-\infty$ can be simply obtained by reflecting $f_j$. Namely, we define 
%\[
$g_j (x,\la) \coloneqq f_j(-x,\la)$. 
%\]
Then $g_1,g_2,g_3,g_4$ forms a basis of the solution space to \eqref{Eq:B.2Jun}. 
\end{itemize}
For $\la=0$, the solution to $\boH_\eps f=f$ is spanned by
\begin{align}\label{eq:m4}
f_1(x,0),\quad \rd_\la f_1(x,0),\quad f_3(x,0),\quad f_4(x,0).
\end{align}
%%%
%%% Asymptotics of \rd_\la f_1
%%%
% $\rd_\la f_1$ has the asymptotics
% \[
% \rd_\la f_1 (x,0) = 
% \begin{pmatrix}
% ix \\ 0
% \end{pmatrix}
% + 
% O(xe^{-\gam x}).
% \]
%%%
%%% Asymptotics of \rd_\la f_1 end
%%%
%\end{itemize}
%
%
\subsection{Generalized kernel}
%
%
%\begin{screen}
\begin{prop}
If $\eps>0$, then the generalized kernel for $\boA_\eps$ is
\[
\boN_g (A_\eps) = \Span \left\{
\begin{bmatrix}
    0 \\ b
\end{bmatrix}
,
\begin{bmatrix}
    L_\eps^{-1} b \\ 0
\end{bmatrix}
.
\right\}
\]
% \[
% \si(\boA_\eps) \cap i[-1,1] = \emptyset.
% \]
\end{prop}
%\end{screen}
%
%
\begin{proof}
Suppose $\eps>0$. Recall that $L_\eps$ is invertible 
since $L_\eps\ge \la_\eps>0$ for some $\la_\eps>0$. 
Hence, we have
$$
\boA_\eps 
\begin{bmatrix}
f \\ g
\end{bmatrix}
= 
0 \quad \Longleftrightarrow \quad
\left\{
\begin{aligned}
L_0 g &=0,\\
-L_\eps f & =0.
\end{aligned}
\right.
\quad \Longleftrightarrow \quad
\left\{
\begin{aligned}
g &= c\, b,\\
f & =0.
\end{aligned}
\right. \quad
( c\in\R).
$$
$$
\begin{aligned}
\boA_\eps^2 
\begin{bmatrix}
f \\ g
\end{bmatrix}
= 
0 \quad \Longleftrightarrow \quad
\left\{
\begin{aligned}
-L_0L_\eps f &=0,\\
-L_\eps L_0 g & =0.
\end{aligned}
\right.
&
\quad \Longleftrightarrow \quad
\left\{
\begin{aligned}
L_\eps f &= c\, b,\\
L_0 g & =0.
\end{aligned}
\right. \quad
( c\in\R)
\\
&
\quad \Longleftrightarrow \quad
\left\{
\begin{aligned}
f &= c\, L_\eps^{-1}b,\\
g & = c'\, b.
\end{aligned}
\right. \quad
( c,c'\in\R).
\end{aligned}
$$
Furthermore, 
%if $\inp{R_\eps}{b}_{L^2(\R)}\neq 0$, 
$$
\begin{aligned}
\boA_\eps^3 
\begin{bmatrix}
f \\ g
\end{bmatrix}
= 
0 
\quad \Longleftrightarrow \quad
\boA_\eps^2 (\boA_\eps \begin{bmatrix}
f \\ g
\end{bmatrix}) =0
% &\quad \Longleftrightarrow \quad
% \left\{
% \begin{aligned}
% -L_0L_\eps L_0 g & =0,\\
% L_\eps L_0 L_\eps f &=0.\\
% \end{aligned}
% \right.
% \\
&
\quad \Longleftrightarrow \quad
\left\{
\begin{aligned}
L_0 g &= c\, L_\eps^{-1}b,\\
L_\eps f & = c' b.
\end{aligned}
\right. \quad
( c,c'\in\R).
\end{aligned}
$$
From the second equation, we have $f = c' L_\eps^{-1} b$. On the other hand, taking inner product of the first equation with $b$, we have
$$
0 = \inp{L_0 g}{b}_{L^2} = c \inp{L_\eps^{-1}b}{b}_{L^2} 
\quad\Longleftrightarrow \quad 
c=0.
$$
Hence $g\in \Ker L_0 = \Span b$, and thus $\Ker A_\eps^3 = \Ker A_\eps^2$ follows.
\end{proof}

\subsection{Nonexistence of internal modes and threshold resonances}
\begin{prop}\label{Psub4.1}
Suppose $\eps>0$. Then $\boA_\eps$ has no internal mode or threshold resonance.
% the following are true:\\
% (i) 
% %If $\eps>0$, then 
% $\boA_\eps$ has no internal mode; i.e. there is no eigenvalue of $\boA_\eps$ in $i[-1,1]\setminus \{0\}$.\\
% (ii) $\boA_\eps$ has no threshold resonance at $\pm i$; i.e, there is no 
% $u\in H^2_{\loc}(\R)^2 \cap  L^\infty(\R)^2\setminus L^2(\R)^2$ satisfying either of $\boA_\eps u =\pm iu$.
\end{prop}
%}
%
% Following \cite{MR2219305}*{Proposition 9.2}, we split the proof into two regimes for $\eps$. 
We follow the strategy of \cite{MR2219305}*{Proposition 9.2}. However, in our setting, the spectral properties of $L_\eps$ depend on $\eps$, requiring us to divide the proof into two regimes. 
By Lemma \ref{L:AugA5}, for $0 < \eps < \eps_0$, $L_\eps$ has point spectrum below $1$, requiring an extra argument to restrict to $(\Span b)^\perp$. For $\eps \ge \eps_0$, the bound $L_\eps \ge 1$ significantly simplifies the argument.\par 
%}
%
%
\subsubsection{A sufficient condition}\label{Sop2.2}
First, for $0<\eps<\eps_0$, 
we give a sufficient condition for excluding the internal modes.
% We begin with the case $0<\eps<\eps_0$. 
% %
% We first give a sufficient condition to obtain the conclusion. 
Define
\[
h(\eps,\la) \coloneqq \inp{(L_\eps -\la)^{-1}b}{b}_{L^2},\qquad 
\la\in (\la_\eps ,1),
\]
%$h(\la)$ is well-defined on $\la\in (\la_\eps,1)$,
where $\la_\eps>0$ is the least eigenvalue of $L_\eps$.
$h(\eps,\la)$ is strictly increasing in $\lambda$, since
\[
\frac{\rd h}{\rd \la}(\la) = \inp{(L_\eps-\la)^{-2}b}{b} >0.
\]
%Since $(L_\eps-\la)^{-1}$ is singular at $\la=\la_\eps$, we have $\lim_{\la \to \la_\eps +0} h(\la) = -\infty$. 
%
%
%{\color{red}
% \begin{prop}
% $h(\la)<0$ for $\la\in (\la_\eps,1)$.
% \end{prop}
%}
%
%
\begin{lem}\label{Lop2.2}
Let $0<\eps<\eps_0$. 
If $h(\eps,\la)<0$ for all $\la\in (\la_\eps,1)$, then Proposition \ref{Psub4.1} is true.
\end{lem}
%
%
% We first give \textcolor{blue}{a} proof of Lemma \ref{Lop2.2}, then verify the hypothesis later. 
The proof of Lemma \ref{Lop2.2} is almost parallel to that in \cite{MR2219305}*{Proposition 9.2}; in our case, the argument gets much simpler since there is no non-positive spectrum for $L_\eps$. 
For clarity, we reproduce the proof in the sequel.\par
%
%\subsubsection{Proof of Lemma \ref{Lop2.2}}
%\label{Sop2.2}
%\subsubsection*{\fbox{Step 0.}}
We denote by $\phi_\eps$ the positive $L^2$-normalized eigenfunction associated with $\la_\eps$. Let $P_0$ be the orthogonal projection onto $\Span b$, and set $\ovl{P_0}=1-P_0$. 
Setting $L^2_\perp = \ovl{P_0}(L^2(\R))$ and 
$H^s_\perp = \ovl{P_0}(H^s(\R))$ for $s>0$, we introduce 
\[
L_{\eps,0} \coloneqq \ovl{P_0} L_\eps \ovl{P_0} ,\qquad 
L_{0,0} \coloneqq \ovl{P_0} L_0 \ovl{P_0}
%\colon H^2_\perp \subset L^2_\perp\to L^2_\perp.
\]
as an operator from $H^2_\perp \subset L^2_\perp$ to $L^2_\perp$. Note that $L_{\eps,0}$, $L_{0,0}$ is still self-adjoint on $L^2_\perp$.

%\subsubsection*{\fbox{Step 1.}}
% Let $P_0$ be the orthogonal projection into $\Span b$. Set $\ovl{P_0}=1-P_0$.
\begin{lem}\label{L4.2}
$L_{\eps,0}$ has no eigenvalue in $(-\infty,1)$.
\end{lem}
\begin{proof}
Suppose there exists $\eta\in H^2_\perp\setminus \{0\}$, $\la<1$ such that  $L_{\eps,0} \eta =\la \eta$. 
Then $\eta\in L^2_\perp =\Ran (\ovl{P_0})$, and
\[
\inp{(L_\eps -\la)\eta}{\eta} 
= \inp{(L_{\eps,0} -\la)\eta}{\eta} =0.
\]
In particular,
\[
\la_\eps \nor{\eta}{L^2}^2 \le 
\inp{L_\eps \eta}{\eta} = \la \nor{\eta}{L^2}^2.
\]
Since $\eta\neq 0$, we have $\la\ge \la_\eps$. 
If $\la=\la_\eps$, then $\eta$ is a minimizer of $\inf_{f\in H^2} \frac{\inp{L_\eps f}{f}}{\nor{f}{L^2}^2}$, which implies $\eta= c \phi_\eps$ for some $c\in\C\setminus \{0\}$. However, 
$\inp{\phi_\eps}{b}>0$ since both $\phi_\eps$ and $b$ are strictly positive, which contradicts $\eta\in \Ran (\ovl{P_0})$. Hence $\la \in (\la_\eps,1)$. Then we claim that $h(\eps,\la)=0$, contradicting the hypothesis. 
Indeed, since $\ovl{P_0}(L_{\eps}-\la) \eta =0$, we have $(L_{\eps}-\la) \eta =  c b$ for some $c\in\C$.  
Noting $\la\in (\la_\eps,1)$, we can write 
$\eta = c(L_\eps-\la)^{-1} b$. 
Clearly $c\neq 0$. 
Therefore, 
\[
h(\eps,\la) = \inp{(L_\eps -\la)^{-1}b}{b} = \inp{%\frac 1c
c^{-1}\eta}{b}=0.
\]
\end{proof}
\begin{lem}\label{L4.3}
The following are true:\\
(i) $\inf_{\nor{f}{L^2}=1}\inp{L_{\eps,0}f}{f}\ge 1$.\\
(ii) $\inp{L_{0,0}^{-1}f}{f} < \nor{f}{L^2}^2$ for all $f \in L^2_\perp\setminus \{0\}$.
\end{lem}
\begin{proof}
(i) 
Set $\la_* \coloneqq \inf_{f\in H^1_\perp,\ \nor{f}{L^2}=1}\inp{L_{\eps,0}f}{f}$. 
Note that $\la_*\ge \la_\eps>-\infty$. 
Suppose to the contrary that $\la_*<1$. 
Then, since $\si_{\ess} (L_{\eps,0}) = [1,\infty)$, $\la_*$ is a point spectrum of $L_{\eps,0}$, which contradicts Lemma \ref{L4.2}.\par
(ii) Since $\nor{L_{0,0}^{-1}}{}\le 1$, we have $\inp{L_{0,0}^{-1}f}{f} \le \nor{f}{L^2}^2$. 
Suppose the equality holds at some $f=f_*\in L^2_\perp\setminus \{0\}$. Then $f_*$ is a maximizer of
\[
\sup_{f\in L^2_\perp,\  f\neq 0} \frac{\inp{L_{0,0}^{-1}f}{f}}{\nor{f}{L^2}^2}.
\]
By variational principle, we have $L_{0,0}^{-1}f_* = f_*$, and hence $L_0 f_* =f_*$. This contradicts $\ker (L_0-1)=\{0\}$ in  $L^2$.
\end{proof}

% %\subsubsection*{\fbox{Step 2.}}
% \begin{lem}\label{L4.4}
% $\inp{L_{0,0}^{-1}f}{f} < \nor{f}{L^2}^2$ for all $f \in L^2_\perp\setminus \{0\}$.
% \end{lem}
% %
% \begin{proof}
% Since $\nor{L_{0,0}^{-1}}{}\le 1$, we have $\inp{L_{0,0}^{-1}f}{f} \le \nor{f}{L^2}^2$. 
% Suppose the equality holds at some $f=f_*\in L^2_\perp\setminus \{0\}$. Then $f_*$ is a maximizer of
% \[
% \sup_{f\in L^2_\perp,\  f\neq 0} \frac{\inp{L_{0,0}^{-1}f}{f}}{\nor{f}{L^2}^2}.
% \]
% By variational principle, we have $L_{0,0}^{-1}f_* = f_*$, and hence $L_0 f_* =f_*$. \textcolor{blue}{This contradicts $\ker (L_0-1)=\{0\}$ in  $L^2$.}
% \end{proof}
%
%
\begin{proof}[Proof of Lemma \ref{Lop2.2}]
Suppose there exists an eigenvalue $iE$ with $E\in [-1,1]\setminus \{0\}$. Let $u={}^t(u_1,u_2)$ be an associated eigenfunction. Then, $f_1\in \Ran \ovl{P_0}\setminus \{0\}$, and satisfies $E^2f_1=L_0L_\eps f_1  = L_{0,0} L_{\eps,0} f_1$. Therefore, Lemma \ref{L4.3} implies
\begin{align*}
\nor{f_1}{L^2}^2
>
\inp{L_{0,0}^{-1}f_1}{f_1} 
% =
% \frac{1}{E^2} \inp{\ovl{P_0} L_\eps\psi}{\psi}
= \frac{1}{E^2} \inp{L_{\eps,0}f_1}{f_1} 
\ge
% \overset{\text{Lemma \ref{L4.3}}}{\ge}
\nor{f_1}{L^2}^2,
\end{align*}
which is a contradiction.
\end{proof}
%
% %\subsubsection*{\fbox{Step 3.}}
% \begin{lem}\label{L4.5}
% Suppose there exists an eigenvalue $iE$ with $E\in [-1,1]\setminus \{0\}$. 
% \textcolor{blue}{Then,} there exists $\psi\in H^4_\perp\setminus \{0\}$ such that $L_0L_\eps \psi =E^2 \psi$. 
% \end{lem}
% %
% \begin{proof}
% Let 
% $f= {}^t(f_1,f_2)
% % \begin{pmatrix}    f_1\\f_2
% % \end{pmatrix}
% $ 
% be an associated eigenfunction. 
% Then we have $L_0L_\eps f_1 =E^2f_1$ and $L_\eps L_0 f_2 = E^2 f_2$. 
% % \[
% % \boA_\eps^2 f = -E^2 f\quad \Longleftrightarrow\quad 
% % \begin{cases}
% %     L_0L_\eps f_1 =E^2f_1,\\
% %     L_\eps L_0 f_2 = E^2 f_2.
% % \end{cases}
% % \]
% If $f_1=0$, then $iE f_2= L_\eps f_1 = 0$ and thus $f_2=0$, which contradicts $f\neq 0$. 
% Thus $f_1\neq 0$. 
% Moreover, we have $f_1 = E^{-2}L_0L_\eps f_1\in \Ran (L_0) = \Ran \ovl{P_0}$. 
% Hence, it suffices to set $\psi=f_1$. 
% \end{proof}
% %
% %
% %\textit{Conclusive step.} 
% \begin{proof}[Proof of Lemma \ref{Lop2.2}]
% Suppose there exists an eigenvalue $iE$ with $E\in [-1,1]\setminus \{0\}$. Let $\psi$ be as in Lemma \ref{L4.5}. Note that $E^2\psi= L_0L_\eps\psi = L_{0,0} \ovl{P_0} L_\eps\psi$. 
% Thus, we have%since $E^2\le 1$, we have
% \begin{align*}
% \nor{\psi}{L^2}^2
% \overset{\text{Lemma \ref{L4.4}}}{>}
% \inp{L_{0,0}^{-1}\psi}{\psi} 
% =
% \frac{1}{E^2} \inp{\ovl{P_0} L_\eps\psi}{\psi}
% = \frac{1}{E^2} \inp{L_{\eps,0}\psi}{\psi} 
% \overset{\text{Lemma \ref{L4.3}}}{\ge}\nor{\psi}{L^2}^2,
% \end{align*}
% which is a contradiction.
% \end{proof}
%
%
%
%
%
% \subsubsection{Verification of the \textcolor{blue}{hypothesis}}
\subsubsection{Proof of Proposition \ref{Psub4.1}}
We begin with the proof of the nonexistence of internal modes. 
%We first consider the case 
Let $0<\eps<\eps_0$. 
%Suppose $0<\eps<\eps_0$. 
By Lemma \ref{Lop2.2}, it suffices to verify that $h(\eps,\la)<0$. To this end, we exploit the feature of $L_\eps$. 
% The most nontrivial part is the verification of the hypothesis $h(\eps,\la)<0$. 
% We verify this by \textcolor{blue}{exploiting the feature} of $L_\eps$. 
%
%
% \begin{lem}\label{P4.7}
% If $\eps>0$, then 
% $h(\eps,\la)<0$ for $\la\in (\la_\eps,1)$. 
% \end{lem}
% %
% %
% \begin{proof}
Set $f_\eps = (L_\eps-\la)^{-1}b$. 
Then $h(\eps,\la) = \inp{f_\eps}{b}$. Moreover, we have 
%\[
$f_0 = (L_0-\la)^{-1}b = -\frac{1}{\la} b$. 
%\]
Now, differentiating $(L_\eps-\la)f_\eps =b$ by $\eps$ yields 
$\rd_\eps f_\eps = -2(L_\eps - \la)^{-1} (bf_\eps)$. 
% \[
% (L_\eps -\la) \rd_\eps f_\eps + 2bf_\eps =0\quad
% \Longleftrightarrow\quad 
% \rd_\eps f_\eps = -2(L_\eps - \la)^{-1} (bf_\eps).
% \]
Therefore,
\begin{align*}
\frac{\rd h}{\rd \eps}(\eps,\la)
= \inp{\rd_\eps f_\eps}{b}
&= -2\inp{(L_\eps - \la)^{-1} (bf_\eps)}{b}\\
% &= -2 \inp{ bf_\eps}{(L_\eps - \la)^{-1}b} \\
&= -2 \inp{bf_\eps}{f_\eps} 
= -2\int_\R b(x) |f_\eps(x)|^2 dx <0,\\
h(0,\la) &= \inp{f_0}{b} = 
-\frac 1\la \nor{b}{L^2}^2 = -\frac{2}{\la}<0,
\end{align*}
by noting $\la>\la_\eps >0$ for $\eps>0$. 
Since $\eps\mapsto \la_\eps$ is monotonically increasing, 
if $\la>\la_\eps$, then $ h(\cdot,\la)$ is smooth on $[0,\eps]$. Therefore, we conclude that $h(\eps,\la)< h(0,\la) <0$.\par
%\end{proof}
%
%
% \subsubsection{Case $\eps\ge \eps_0$}
%
Next, let $\eps \ge \eps_0$. 
Since $L_\eps\ge 1$, 
the same assertion as Lemma \ref{L4.3} holds without extra conditions. Hence, the argument from Section \ref{Sop2.2} concludes the proof.\par
Finally, we consider the threshold resonances. 
Since the proof is almost identical to that in \cite{MR2219305}*{Proposition 9.2}, we only sketch the proof. Since $\pm$ can be interchanged by complex conjugation, we only discuss the $-i$ case. 
Suppose to the contrary that there exists $u\in (L^\infty)^2\setminus (L^2)^2$ such that $\boA_\eps u =-iu$. Let $\chi\in C_0^\infty(\R)$ be a cut-off function with 
$\chi(x) =1$ on $x\in[-1,1]$ and 
$\chi(x)=0$ on $|x|\ge 2$, and define $\psi_\del \coloneqq \ovl{P_0} (\psi\chi(\del\cdot))$ for $\del>0$. Then by the same argument as in \cite{MR2219305}, it follows that $\inp{L_\eps \psi_\del}{\psi_\del} < \nor{\psi_\del}{L^2}^2$ for sufficiently small $\del$, contradicting Lemma \ref{L4.3}.

\subsection{Embedded eigenvalues}

\begin{prop}
For $\eps\in\R$, $\boA_\eps$ has no embedded eigenvalue. Namely, there exists no pair of $E>1$ and $u\in H^2(\R\colon \C^2)\setminus \{0\}$ satisfying either of $\boA_\eps u = \pm iEu$.
%satisfying
% \[
% {}^\exists u\in H^2(\R\colon \C^2)\setminus \{0\}\quad \text{s.t.}\quad 
% \boA_\eps u = \pm iEu.
% \]
\end{prop}
To rule out embedded eigenvalues, we adapt a recent approach by Collot, Germain, and Pacherie \cite{MR5093257}, which represents the eigenfunction using the kernels of the off-diagonal operators in $\boA_\eps$. While $L_\eps$ lacks a kernel in our setting, substituting the Jost solution to $L_\eps f=0$ allows the proof to proceed largely as in \cite{MR5093257}.

%their idea perfectly fits our case by replacing .
%\subsubsection*{Positivity of Jost function}
\begin{lem}\label{L:em3}
Let $f=f_\eps$ be the Jost solution to $L_\eps f= 0$ with the asymptotics
\[
f(x) = e^{-x} + O(e^{-(1+\gam)x}),\qquad x\to+\infty.
\]
Then $f>0$ on $\R$.  
\end{lem}
\begin{proof}
Suppose the contrary. Let
\[
x_* =  \inf\{x\in \R \mid f(t) >0 \text{ for } t\in(x,\infty)\}.
\]
%(The largest zero of $f$.) 
Let us consider the Wronskian $W(x) = f'b-fb'$. 
Then $W(\infty)=0$ and 
$W'(x) = 2\eps fb >0$ on $(x_*,\infty)$.
Hence, we have $W(x_*)<0$. On the other hand, since $f'(x_*)>0$, we have $W(x_*) = f'(x_*)b(x_*) >0$, a contradiction. 
\end{proof}

%\subsubsection*{Inversion formula}
This is the key factor of the proof of \cite{MR5093257}:
\begin{lem}\label{L:em1}
(i) Consider the equation $L_0 u=v$. Assume that
\begin{align}\label{eq:em1}
\sup_{x\in(0,\infty)}(|u(x)|+|v(x)|)e^{(1+\gam)x} <\infty.
\end{align}
Then we have
\begin{align}\label{eq:em2}
u(x) =  - b(x) \int_x^\infty \frac{1}{b(t)^2} \int_t^\infty b(s)v(s)ds.
\end{align}
(ii) Consider the equation $L_\eps u=v$. Suppose that $u,v$ satisfies \eqref{eq:em1}. Then
\begin{align}\label{eq:em4}
u(x) =  - f(x) \int_x^\infty \frac{1}{f(t)^2} \int_t^\infty f(s)v(s)ds.
\end{align}
\end{lem}
\begin{proof}
Direct computation yields
\[
L_\eps = G^*G,\qquad G=f\rd_x \frac{1}{f} = \rd_x - \frac{f'}{f}.
\]
With this propery, the proof proceeds almost similarly to \cite{MR5093257}*{Lemma 2.2}, and hence we omit the details.
% \dlt{
% Therefore, the proof is the same as in} \cite{MR5093257}*{Lemma 2.2}.
\end{proof}
For completeness, we reproduce how to derive a contradiction. 
Suppose to the contrary that $\boH_\eps$ has a $L^2$-solution of
\begin{align}\label{eq:em5}
\boH_\eps 
\begin{pmatrix}
    f \\g
\end{pmatrix}
= E
\begin{pmatrix}
    f \\g
\end{pmatrix},
\qquad E>1.
\end{align}
By setting $u=g-f$, $v=g+f$, we have 
$L_0u=-Ev$ and $L_\eps v= -Eu$. 
% \[
% \begin{cases}
%     L_0 u= -E v,\\
%     L_\eps v = - E u.
% \end{cases}
% \]
%
%As observed in the Jost solutions, 
%As a consequence of 
Here, the results summarized in Section \ref{SB1} imply that 
$L^2$-solutions to \eqref{eq:em5} 
% \begin{equation}\label{eq:em5}
% \boH_\eps f = (1+\la^2)f
% \end{equation}
behave as 
$C_\pm e^{-\mu |x|} 
\begin{pmatrix}
    0 \\1
\end{pmatrix}
+ O(e^{-(\mu+\gam) |x|})$ as $x\to \pm\infty$ 
% \[
% \begin{cases}
% C_+ e^{-\mu |x|} 
% \begin{pmatrix}
%     0 \\1
% \end{pmatrix}
% + O(e^{-(\mu+\gam) |x|}),\qquad & x\to\infty,\\
% C_-
% e^{-\mu |x|} 
% \begin{pmatrix}
%     0 \\1
% \end{pmatrix}
% + O(e^{-(\mu+\gam) |x|}),\qquad & x\to-\infty
% \end{cases}
% \]
with $\mu= \sqrt{2+E}$, $0<\gam<1$, and $C_\pm\in \C\setminus \{0\}$. 
It holds that the eigenspace for \eqref{eq:em5} is one dimensional, and hence $f$ is either even or odd. 
We may suppose $C_+=1$, then $f$ is real-valued. By this, we have
\begin{align}\label{eq:em6}
\begin{pmatrix}
    u\\v
\end{pmatrix}
= 
e^{-\mu x}
\begin{pmatrix}
    1\\1
\end{pmatrix}
+ O(e^{-(\mu+\gam)x}),\qquad x\to+\infty.
\end{align}
Therefore, $u$, $v$ satisfies \eqref{eq:em1}, and thus Lemma \ref{L:em1} yields
\begin{align}
u(x) =   E b(x) \int_x^\infty \frac{1}{b(t)^2} \int_t^\infty b(s)v(s)ds,\label{eq:em7}\\
v(x) =   Ef(x) \int_x^\infty \frac{1}{f(t)^2} \int_t^\infty f(s)u(s)ds.\label{eq:em8}
\end{align}
Here, \eqref{eq:em6} implies $u(x),v(x)>0$ for $x\gg 1$. 
Since $b(x)$, $f(x)>0$ by Lemma \ref{L:em3}, \eqref{eq:em7} and \eqref{eq:em8} imply that $u(x),v(x)>0$ for $x>0$. 
% We claim that $u(x)>0$, $v(x)>0$ for $x\in (0,\infty)$. In fact, \eqref{eq:em6} implies that this is true for $x\gg 1$. Since $b(x)$, $f(x)>0$ by Lemma \ref{L:em3}, \eqref{eq:em7} and \eqref{eq:em8} lead to the claim. 
In particular, this yields 
$u(0)>0$ and $u'(0)<0$. 
However, since $u$ is either even or odd, we must have
$u(0)=0$ or $u'(0)=0$, which leads to a contradiction.

\section{Resolvent estimate}\label{SC}
Verifying Propositions \ref{P4:3.3} and \ref{P4:3.2} requires a resolvent estimate akin to \cite{MR4864876}*{Proposition 8.6}. Since their proof relies on proximity to the integrable equation, we develop an alternative approach for our setting.

For $\boH=\boH_\eps$ and real $|E|\ge1$, define
\[
R_\boH^\pm(E) \coloneqq 
(\boH_\eps-(E\pm i0))^{-1} \equiv \lim_{\mu\to +0} (\boH_\eps-(E\pm i\mu))^{-1},
\]
understood in the usual limiting sense \cites{MR2219305,MR4864876}. For $\la\in\R$, let $\boR_\boH^\pm (x,y,\la)$ be the integral kernel of $R_{\boH}^\pm(1+\la^2)$. We claim the following:

\begin{prop}[Resolvent bound]\label{PC.1}
Suppose $\eps>0$. Then, for all $\la\in\R$, we have
\begin{align}\label{eq:m9}
|R^\pm_{\boH} (x,y,\la)| \le C 
\begin{cases}
1+x^- + y^+\quad & \text{if } x\ge y,\\
1+x^+ + y^-\quad & \text{if } x\le y,
\end{cases}
\end{align}
where $x^\pm = \max\{0, \pm x\}$. 
\end{prop}
The proof of the counterpart to Proposition \ref{PC.1} in \cite{MR4864876} uses the representation formula \eqref{eq:m8}, 
which involves 
%and 
the Jost matrix $D(\la)$ from \cite{MR2219305}. For $|\la|\gtrsim 1$, the estimate follows via the reflection and transmission coefficients $A(\la)$ and $B(\la)$, as the diagonality of $A(\la)D^{-1}(\la)$ cancels the exponentially growing factors in \eqref{eq:m8}. 
The regime $|\la|\ll 1$ is more delicate because $A(\la)D^{-1}(\la)$ is singular at $\la = 0$. In \cite{MR4864876}, this is handled by defining $D(\la)$ via a Jost solution $f_1(\cdot ,\la)$ with at most linear growth as $x\to -\infty$. 
% constructing a Jost solution $f_1(\cdot,\la)$ (unique up to adding $cf_4$, $c\in\R$) with at most linear growth as $x\to -\infty$. 
However, %this 
the construction of such $f_1$ 
requires the exponential growth of $f_3(\cdot,0)$ as $x\to -\infty$ (\cite{MR4864876}*{Claim 8.12}), which relies on explicit integrable formulas and fails outside small perturbations of the integrable regime. Instead, we introduce modified coefficients $A^\#(\la)$ and $B^\#(\la)$ via a new basis for which $A^\#(\la)D^{-1}(\la)$ remains regular at $\la = 0$, providing a robust method for more general settings.

\subsection{Preliminaries}
%First of all
For \eqref{eq:m9}, it suffices to consider only the $+$ case, since
\[
R_\boH^+ (1+\la^2) = \ovl{R_\boH^- (1+\la^2)},\qquad \text{and thus }\qquad R_\boH^+(x,y,\la) = \ovl{R^-_\boH (x,y,\la)}.
\]
Moreover, it suffices to consider the case $x\ge y$. In fact, we have
\[
R_\boH^+ (x,y,\la) = \si_3 R_\boH^+ (y,x,\la)^t \si_3.
\]
%
%
%\subsubsection*{\fbox{Step 2.}}
Here, we recall 
% \dlt{some identities for the Jost solutions and} 
the representation formula of $\boR_{\boH}^\pm (x,y,\la)$ 
in \cite{MR2219305}*{Sections 5 and 6}. 
By using the Jost solutions \eqref{Eq:B.3Jun}, we define
\begin{align*}
F_1(x,\la) = (f_1(x,\la), f_3(x,\la)),
&&
F_2(x,\la) = (f_2(x,\la), f_4(x,\la)),\\
G_1(x,\la) = (g_2(x,\la), g_4(x,\la)),
&&
G_2(x,\la) = (g_1(x,\la), g_3(x,\la))
\end{align*}
for $\la\in\R$. 
For $F,G:\R\to M_2(\C)$, we define the matrix Wronskian as
\[
W[F,G] = (F')^t G - F^t G'.
\]
Note that if $F,G$ are solutions to \eqref{Eq:B.2Jun}, then $W[F,G]$ is independent of $x$. Note also that 
if we write $F=(f_1,f_2)$, $G=(g_1,g_2)$, then
\[
W[F,G] = 
\begin{pmatrix}
W[f_1,g_1] & W[f_1,g_2]\\
W[f_2,g_1] & W[f_2,g_2]
\end{pmatrix}
.
\]
In particular, $W[F,F]$ is not necessarily $0$, as the off-diagonal part can be non-zero. Now we define
\[
D(\la) \coloneqq W[F_1(\cdot,\la),G_2(\cdot,\la)],
\]
%$D(\la)$ is
which we call the \textit{Jost matrix} by following \cite{MR2598115}. Then the following properties are known:
\begin{itemize}
\item $D(\la)^t=D(\la) = \ovl{D(-\la)}$.
\item $D(\la)$ is smooth in $\la$.
\item If $\boH_\eps$ has no threshold resonances or embedded eigenvalues, then $D(\la)$ is invertible for all $\la\in\R$. In particular, this is true for $\eps>0$.
%It is known that if $\eps>0$, then $\det D(\la)$ does not vanish for all $\la\in\R$. (This follows from the absence of threshold resonance nor embedded eigenvalues.)
% If $\boA_\eps$ does not have threshold resonance nor embedded eigenvalues, then $\det D(\la)\neq 0$ on $\R$. Therefore, by Proposition \ref{P3.2Jun}, this condition holds if $\eps>0$.
\end{itemize}
%
%
%
% Finally, we recall the identities related to the reflection and transmission coefficients. 
% For $\la>0$, 
% there exist matrices $A(\la)$, $B(\la)$ such that
% \begin{align}\label{eq:m7}
% F_1(x,\la) = G_1(x,\la)A(\la) + G_2(x,\la) B(\la),
% \end{align}
% with the following properties:
% \begin{itemize}
% \item $G_2(x,\la) = F_1(x,\la) B(\la) + F_2(x,\la) A(\la)$,
% \item $D(\la) = W[F_1,G_2]= A(\la)^t 
% \begin{pmatrix}
% 2i\la & 0\\
% 0 & -2\sqrt{2+\la^2}
% \end{pmatrix}$
% ,
% \item $ W[F_1,G_1]= -B(\la)^t 
% \begin{pmatrix}
% 2i\la & 0\\
% 0 & -2\sqrt{2+\la^2}
% \end{pmatrix}$.
% \end{itemize}
%
%
For $\la\ge 0$, we have the following representation formula of the resolvent:
%\begin{screen}
\begin{align}\label{eq:m8}
R_\boH^+ (x,y,\la) = 
\begin{cases}
- F_1(x,\la) D^{-1}(\la) G_2(y,\la)^t \si_3\quad & x\ge y,\\
- G_2(x,\la) D^{-1}(\la) F_1(y,\la)^t  \si_3 \quad &x\le y.
\end{cases}
\end{align}
%\end{screen}
% (The two formulas in the right hand side coincide at $x=y$, as verified in \cite{MR2219305}.)
%
%
\subsection{Modified transmission and reflection coefficients}
%
%
%To deal with the case $|\la|\lesssim 1$, 
We introduce a new version of transmission and reflection coefficients, modified from those in \cite{MR2219305}. 
% Jost matrices and observe their coefficients. 
First, based on the Jost functions in \eqref{Eq:B.3Jun}, we define
\[
\tilde{f}_2(x,\la) \coloneqq \frac{f_1(x,\la)-f_1(x,-\la)}{2\la} \quad \text{if } \la\neq 0,\qquad 
\rd_\la f_1(x,0) \quad \text{if } \la=0.
\]
\begin{lem}
$\tilde{f}_2$ is smooth both in $x$ and $\la$. Moreover, for $x\ge 0$ and $|\la|\lesssim 1$, we have
\begin{align}\label{eq:Aug61}
\jb{x}^{-1}|\tilde{f}_2 (x,\la)|\le C.
\end{align}
\end{lem}
\begin{proof}
We observe
\begin{align*}
\tilde{f}_2(x,\la) 
= \frac 12 \int_{-1}^1 e^{it\la x}
[\rd_\la m_1] (x,t\la) dt
+ \frac 12 \int_{-1}^1 ix e^{it\la x} m_1(x,t\la) dt,
\end{align*}
from which it can be easily seen that $\tilde{f}_2$ is smooth both in $x$ and $\la$. Moreover, we have
\begin{align*}
|e^{it\la x}
[\rd_\la m_1+ixm_1] (x,t\la)|\lesssim \jb{t\la}^{-1} (1+x) e^{-\gam x} \le C,
\end{align*}
and hence \eqref{eq:Aug61} follows.
\end{proof}
We set $\tilde{g}_2(x,\la) = \tilde{f}_2(-x,\la)$. Then, we introduce
\begin{align*}
F_1^\# (\la) \coloneqq F_1(\la) = (f_1, f_3),&& F_2^\# (\la) \coloneqq  (\tilde{f}_2,f_4),\\
G_1^\# (\la) \coloneqq (\tilde{g}_2, g_4),&& G_2^\# (\la) \coloneqq G_2(\la) = (g_1,g_3).
\end{align*}
\begin{lem}\label{LC.3}
For all $\la\in\R$, 
$f_1$, $\tilde{f}_2$, $f_3$, $f_4$ is a basis of the solution space for $\boH f=(1+\la^2)f$. Moreover, if we write
\begin{align}\label{Eq:Jun3}
F_1(x,\la) = G_1^\# (x,\la) A^\#(\la) + G_2^\# (x,\la) B^\#(\la).
\end{align}
the following hold:\\
(i) $A^\#(\la)$, $B^\#(\la)$ are smooth for all $\la\in\R$;\\
(ii) $A^\#(\la)D(\la)^{-1} = 
\begin{pmatrix}
i & 0\\
0 & -\frac{1}{2\sqrt{2+\la^2}}
\end{pmatrix}.$\\
(iii) $B^\#(\la) = 
\begin{pmatrix}
-i & 0\\
0 & \frac{1}{2\sqrt{2+\la^2}}
\end{pmatrix}
W[F_1,G_1^\#]^t$.
\end{lem}
%
%\begin{myleftbar}
\begin{proof}
\begin{align*}
W[F_1,G_2^\#] &= ((G_1^\#)' A^\# + (G_2^\#)' B^\#)^t G_2^\# - (G_1^\#A^\# + G_2^\# B^\#)^t (G_2^\#)'\\
&= (A^\#)^t W[G_1^\#,G_2^\#] + (B^\#)^t W[G_2^\#,G_2^\#]\\
&=-(A^\#)^t W[F_2^\#,F_1^\#] - (B^\#)^t W[F_1^\#,F_1^\#]\\
&=-(A^\#)^t 
\begin{pmatrix}
W[\tilde{f}_2,f_1] & W[\tilde{f}_2,f_3]\\
W[f_4,f_1] &
W[f_4,f_3]
\end{pmatrix}
- (B^\#)^t
\begin{pmatrix}
0 & W[f_1,f_3]\\
W[f_3,f_1] &
0
\end{pmatrix}
\\
&=-(A^\#)^t 
\begin{pmatrix}
i & 0\\
0 & 2\sqrt{2+\la^2}
\end{pmatrix}
,
\end{align*}
where
\[
W[f_1,\tilde{f}_2] = W[f_1 , \frac{f_1-f_2}{2\la}] = -\frac{1}{2\la} W[f_1,f_2]= -i.
\]
We also have
\[
W[F_1,G_2^\#] = W[F_1,G_2]= D(\la).
\]
Since $W[F_1,G_2^\#]$ is smooth on $\R$, the smoothness of $A^\#$ follows. 
Using the non-resonant condition, we obtain (ii). On the other hand,
\begin{align*}
W[F_1,G_1^\#] &= ((G_1^\#)' A^\# + (G_2^\#)' B^\#)^t G_1^\# - (G_1^\#A^\# + G_2^\# B^\#)^t (G_1^\#)'\\
&= (A^\#)^t W[G_1^\#,G_1^\#] + (B^\#)^t W[G_2^\#,G_1^\#]\\
&=-(A^\#)^t W[F_2^\#,F_2^\#] - (B^\#)^t W[F_1^\#,F_2^\#]\\
&=-(A^\#)^t 
\begin{pmatrix}
0 & W[\tilde{f}_2,f_4]\\
W[f_4,\tilde{f}_2] &
0
\end{pmatrix}
- (B^\#)^t
\begin{pmatrix}
W[f_1,\tilde{f}_2] & W[f_1,f_4]\\
W[f_3,\tilde{f}_2] & W[f_3,f_4]
\end{pmatrix}
\\
&=-(B^\#)^t 
\begin{pmatrix}
-i & 0\\
0 & -2\sqrt{2+\la^2}
\end{pmatrix}
,
\end{align*}
which concludes the proof.
\end{proof}
%\end{myleftbar}

\subsection{Proof of Proposition \ref{PC.1}}
Now we prove Proposition \ref{PC.1} by dividing into several steps.\par
\textit{Step 1.}
%\subsubsection*{\fbox{Step 3.}}
%For $\la\gg 1$,
For $|\la|\gg 1$, we claim a stronger bound 
\begin{align}\label{eq:m6}
\sup_{x\ge y, |\la| \gg 1} |R_{\boH}^+(x,y,\la)| \le C_0.
\end{align}
For this, it suffices to show
%This bound indirectly follows by showing
\[
\nor{R_\boH^+(\la^2+1)}{L^1\to L^\infty} \le C.
\]
%which has been shown in \cite{MR2219305}. 
By writing
\[
\boH = \boH_0+V,\qquad 
\boH_0 \coloneqq
\begin{pmatrix}
-\rd_x^2 +1 & 0\\
0 & \rd_x^2-1
\end{pmatrix}
,\qquad V\coloneqq 
\begin{pmatrix}
 -2b^2+\eps b & \eps b\\
 -\eps b & 2b^2-\eps b
\end{pmatrix}
,
\]
the resolvent can be written as
%We use the Born series
\[
%(\boH-(\la^2+1+i0))^{-1} 
R_\boH^+(\la^2+1)
= 
\sum_{n=0}^\infty (-1)^n (\boH_0-(\la^2+1+i0))^{-1} (V(\boH_0-(\la^2+1+i0))^{-1})^n,
\]
known as the Born series. 
%If $n=0$, the well-known formula for free resolvent gives
%Here, since $(\boH_0-(\la^2+1+i0))^{-1}$ consists of the free resolvent, we have
Since
\[
\nor{(\boH_0-(\la^2+1+i0))^{-1}}{L^1\to L^\infty} \le C|\la|^{-1}
\]
(see \cite{MR2219305}*{Equation (7.1)}), we have
%This yields
\[
\nor{V(\boH_0-(\la^2+1+i0))^{-1}}{L^1\to L^1} \le \nor{V}{L^1} \nor{(\boH_0-(\la^2+1+i0))^{-1}}{L^1\to L^\infty} \le C|\la|^{-1}.
\]
Therefore, for sufficiently large $\la$, this series converges in $B(L^1, L^\infty)$.\par
%
%
%
%
%\subsubsection*{\fbox{Step 5.}}
\textit{Step 2.} Next, we observe the case $|\la|\lesssim 1$. 
It suffices to show
%, we prove
\[
\sup_{x\ge y, |\la|\lesssim 1} 
|R^\pm_{\boH} (x,y,\la)| \le C (1+x^-+y^+).
\]
By definition, we may restrict $\la\ge0$. 
% \begin{lem}
% For all $|\la|\lesssim 1$, we have
% \[
% |R^\pm_{\boH} (x,y,\la)| \le C 
% \begin{cases}
% 1+x^- + y^+\quad & \text{if } x\ge y,\\
% 1+x^+ + y^-\quad & \text{if } x\le y.
% \end{cases}
% \]
% \end{lem}
%
%
%\begin{proof}
%Again, we may suppose $x\ge y$. 
%We obtain the desired bound from the formula
Let $x\ge y$ and $0\le\la\lesssim 1$. Recall the formula
\[
R_H^+(x,y,\la) = - F_1(x,\la) D^{-1}(\la) G_2(y,\la)^t \si_3.
\]
We divide the cases as follows.\par
\underline{Case 2-1. $x\ge 0\ge y$.} We write
\begin{align*}
&F_1(x,\la) D^{-1}(\la) G_2(y,\la)^t \si_3 \\
&= [f_1(x,\la)\ 0] D^{-1}(\la) [f_1(-y,\la)\ 0]^t \si_3 
+[f_1(x,\la)\ 0] D^{-1}(\la) [0\ f_3(-y,\la)]^t \si_3 \\
&\qquad +[0\ f_3(x,\la)] D^{-1}(\la) [f_1(-y,\la)\ 0]^t \si_3 
+[0\ f_3(x,\la)] D^{-1}(\la) [0\ f_3(-y,\la)]^t \si_3.
\end{align*}
For $x\gg 1$, we have 
%We have
%\[
$|f_1(x,\la)| = |m_1(x,\la)| \le C$ uniformly in $\la\in \R$.
%\]
Since $\la$ is restricted on a compact interval, we have
\begin{equation}\label{Eq:Jun2}
\sup_{|\la|\lesssim 1, x\ge 0} |f_1(x,\la)| <\infty.
\end{equation}
For $y'\ge 0$, we also have 
%\[
$|f_3(y',\la)| = e^{-\sqrt{\la^2+2}y'} |m_3(y',\la)| \le C$.
%\]
Since $|D(\la)^{-1}|$ is uniformly bounded for $0\le \la\lesssim 1$, we conclude that
\[
\sup_{x\ge0\ge y, 0\le\la\lesssim 1}
|F_1(x,\la) D^{-1}(\la) G_2(y,\la)^t \si_3| \le C.
\]
\underline{Case 2-2. $0\ge x\ge y$.} We apply \eqref{Eq:Jun3}:
\begin{align*}
&F_1(x,\la) D^{-1}(\la) G_2(y,\la)^t \si_3 \\
&= G_1^\#(x,\la) A^\#(\la )D^{-1}(\la) G_2(y,\la)^t \si_3
+ G_2^\#(x,\la) B^\#(\la )D^{-1}(\la) G_2(y,\la)^t \si_3.
\end{align*}
For the first term, we have
\begin{align*}
&G_1^\#(x,\la) A^\#(\la ) D^{-1}(\la) G_2(y,\la)^t \si_3 
\\
&= [\tilde{f}_2(-x,\la)\ 0] A^\#(\la )D^{-1}(\la) [f_1(-y,\la)\ 0]^t \si_3 
+[\tilde{f}_2(-x,\la)\ 0] A^\#(\la )D^{-1}(\la) [0\ f_3(-y,\la)]^t \si_3 \\
&\qquad +[0\ f_4(-x,\la)] A^\#(\la )D^{-1}(\la) [f_1(-y,\la)\ 0]^t \si_3 
+[0\ f_4(-x,\la)] A^\#(\la )D^{-1}(\la) [0\ f_3(-y,\la)]^t \si_3.
\end{align*}
% Note that %$A(\la)D(\la)^{-1}$ is 
% \[
% A^\#(\la)D(\la)^{-1} = 
% \begin{pmatrix}
% i &0\\
% 0 & -\frac{1}{2\sqrt{\la^2+2}}
% \end{pmatrix}
% .
% \]
By Lemma \ref{LC.3} (ii), $A^\#(\la)D(\la)^{-1}$ is diagonal, and thus the second and third terms vanish. 
Since $A^\#(\la)D(\la)^{-1}$ is uniformly bounded in $0\le\la\lesssim 1$, 
we bound the first term by% because of \eqref{Eq:May1}. 
\begin{align*}
|[\tilde{f}_2(-x,\la)\ 0] A^\#(\la )D^{-1}(\la) [f_1(-y,\la)\ 0]^t \si_3 |\lesssim |\tilde{f}_2(-x,\la)| |f_1(-y,\la)|\lesssim 1+x^- .
\end{align*}
For the fourth term, we have
\begin{align*}
&|[0\ f_4(-x,\la)] A^\#(\la )D^{-1}(\la) [0\ f_3(-y,\la)]^t \si_3| \\
&\lesssim |f_4(-x,\la)||f_3(-y,\la)|
= e^{-\sqrt{\la^2+2} (x-y)} |m_4(-x,\la)|
|m_3(-y,\la)| \le C
\end{align*}
uniformly in $0\ge x\ge y$ and $0\le\la\lesssim 1$. The $B^\#(\la)$ part consists only of $f_1$ and $f_3$, and thus every component is uniformly bounded in light of Lemma \ref{LC.3} (iii). Hence, we obtain
\[
\sup_{0\ge x\ge y, 0\le\la\lesssim 1}
(1+x^-)^{-1}
|F_1(x,\la) D^{-1}(\la) G_2(y,\la)^t \si_3| \le C.
\]
%\par
\underline{Case 2-3. $x\ge y \ge 0$.} We observe that
\begin{align*}
F_1(x,\la) D^{-1}(\la) G_2(y,\la)^t \si_3 
&= G_2(-x,\la) D^{-1}(\la)  F_1(-y,\la)^t \si_3\\
&= \si_3 \left[ F_1(-y,\la) D^{-1}(\la) G_2(-x,\la)^t \si_3
\right]^t \si_3
\end{align*}
with $0\ge -y\ge -x$. Hence the estimate in Case 2-2 yields
\[
|F_1(x,\la)D^{-1}(\la)G_2(y,\la)^t \si_3| \le C(1+(-y)^-) \le C (1+y^+)
\]
uniformly in $x\ge y\ge 0$ and $|\la|\lesssim 1$. Consequently, the proof of Proposition \ref{PC.1} is complete.

\section{Local well-posedness}\label{SD}

We first observe that \eqref{E5.1Mar} is locally well-posed in $H^s$ with $s\ge 2$. After that, we claim that the solution does not blow up as long as $\nor{\bom-\bom_0}{H^1}$ remains finite. 
\subsection{High regularity well-posedness}
\begin{prop}[Local well-posedness]\label{P:wp1}
Let $s\ge 2$ be an integer. Then for $\bom|_{t=0}$ with 
$\bom|_{t=0}-\bom_0\in H^s(\R)$, 
\eqref{E5.1Mar} has a unique solution $\bom$ with $\bom-\bom_0 \in C((-T_{\max}^-,T_{\max}^+)\colon H^s(\R))$ with $\bom(0,x)=\bom|_{t=0}$, 
where 
$T_{\max}^\pm = T_{\max}^\pm (\bom|_{t=0})$ is the maximal existence time in this class. 
Moreover, %the map $H^s\ni \bom_0 \mapsto \bom \in C_t H^s$ is continuous. 
the solution map $\bom|_{t=0} \mapsto \bom$ is continuous in the topology $H^{s-1}$ to $C_tH^{s-1}$. 
Furthermore, if $T_{\max}^+ (\bom|_{t=0})<\infty$ ($T_{\max}^- (\bom|_{t=0})<\infty$, resp.), then we have
\begin{equation}\label{eq:WP1}
\lim_{t\to T_{\max}^+-0} \nor{\bom(t)-\bom_0}{H^s(\R)} = \infty,\quad 
(\lim_{t\to T_{\max}^-+0} \nor{\bom(t)-\bom_0}{H^s(\R)} = \infty, \text{ resp.})
\end{equation}
\end{prop}
We only give an outline of the proof. 
First notice that 
%Note that 
%We use \eqref{E2.4Aug}. 
%hen 
the well-posedness of \eqref{E5.1Mar} in $\bom_0+H^s$ is equivalent to that for \eqref{E2.4Aug} in $\hat\boe_3+H^s$. 
Then the classical energy method for Schr\"odinger maps is available, yielding the conclusions (see \cites{MR2304153,MR4482034,MR1877231}.) Although the term $-2(b-\eps) w_1\bow'$ from DMI is a derivative nonlinearity, the energy method still works since the derivative in $\bow'$ can be moved to $w_1$ by applying integration by parts. See \cite{MR4482034} in detail.
%, thank to the structure of this term.
\par
We remark that the above strategy fails in $H^1$-regularity. 
%Let us observe how 
%the above classical energy method in $H^1$ fails due to DMI. 
%We compute
In fact, a direct computation yields
\begin{align*}%\label{eq:Sep31}
\begin{aligned}
\frac 12 \rd_t \nor{\bm{w}'}{L^2}^2 &= 
%\inp{\bm{w}'\times \bbH \bm{w}}{\bm{w}'}+ 
\int_\R \bm{w}\times \bm{w}'''\cdot \bm{w}'dx
-2 \int_\R (b-\eps) w_1 \bm{w}'' \cdot \bm{w}'dx + \text{Remainder}
\\
&= -\int_\R (b-\eps) w_1' |\bow'|^2dx + \text{Remainder}
%\label{eq:Sep32}
%
%\inp{\bm{w}\times (\bbH \bm{w})'}{\bm{w}'} 
%- \inp{2(b-\eps) w_1\bm{w}'}{}
\end{aligned}
\end{align*}
by using $|\bow|\equiv 1$, 
where we exploited the cancellation arising from integration by parts. Then, one observes that the first term on RHS cannot be controlled with respect to $\nor{\bow-\hat\boe_3}{L^\infty_tH^1_x}$. 
% the first term in the RHS in the first line \textcolor{blue}{\eqref{eq:Sep31}} vanishes by integration by part. 
% Then, one can observe that
% the first term in RHS in \textcolor{blue}{\eqref{eq:Sep32}} can not be \textcolor{blue}{controlled} \textcolor{blue}{with respect to} $\nor{\bow-\hat\boe_3}{L^\infty_tH^1_x}$.
% will vanish by integration by part. However, the second term can only be bounded by
% \begin{align*}
% \int_\R 2(b-\eps) w_1 \bm{w}'' \cdot \bm{w}'dx &= -\int_\R (b-\eps)' w_1 |\bm{w}'|^2dx - \int_\R (b-\eps) w_1' |\bm{w}'|^2dx\\
% &\le C \nor{\bm{w}'}{L^2}^2 \nor{\bm{w}'}{L^\infty} + C \nor{\bm{w}-\hat\boe_3}{H^1}^3.
% \end{align*}
% % \begin{align*}
% % \int_\R \bm{w}\times \bm{w}'''\cdot \bm{w}'dx = - 
% % \end{align*}
% Hence, this term cannot be handled in $H^1$.
%
%
\subsection{Refined blow-up criterion}
%
%
%When $s=1$, the energy method in the proof of Proposition \ref{P:wp1} does not work, since \eqref{eq:WP5} is no longer true. 
Nevertheless, 
%if we restrict ourselves to small solutions, then 
we are able to improve the blow-up criterion for $H^2$-solution 
%can be improved 
to the boundedness of $H^1$-norm.
\begin{prop}\label{P:wp2}
%There exists $\del_0$ such that the following holds: 
If a solution $\bom\in C((T_{\max}^-,T_{\max}^+)\colon H^2(\R))$ satisfies $\nor{\bom-\bom_0}{L^\infty_tH^1_x}<\infty$, then $T_{\max}^\pm=\pm\infty$.
%Let $\bom|_{t=0}-\bom_0 \in H^2(\R)$
\end{prop}
We use the strategy of modified energy as in \cite{MR4881654}. The goal is to find a functional $\boI$ satisfying
\[
\boI(\bow) \ge C \nor{\bow''}{L^2}^2,\qquad 
|\rd_t \boI(\bow) | \le C (1+P(\nor{\bow-\hat\boe_3}{H^1}) \nor{\bow-\hat\boe_3}{H^2}^2,
\]
which easily leads to Proposition \ref{P:wp2} by a standard argument. 
Toward the construction of $\boI$, 
the starting point is the $H^2$-energy for Schr\"odinger map, explicitly defined as
\begin{align}\label{eq:WP11}
\boE_2 (\bow) \coloneqq \nor{\bow''}{L^2}^2 - \frac 54 \nor{\bow'}{L^4}^4.
\end{align}
For the Schr\"odinger map equation, $\boE_2$ is conserved. If we apply this for \eqref{E2.4Aug}, then $\rd_t \boE_2$ contains the terms whose integrands are of order $|\bow'| |\bow''|^2$, which cannot be bounded as in \eqref{eq:WP11}. To cancel these terms, we modify $\boE_2$ by adding suitable correction terms. 
The main difficulty lies in finding such corrections that precisely eliminate these bad terms.
% \dlt{
% The most nontrivial part is to find such correction terms, which exactly cancel the bad terms.}
%
%
\subsection{Proof of Proposition \ref{P:wp2}}

We first introduce a notation. 
Given a functional $I(\bow)$, $J(\bow)$ of $\bow$,  
we write $I(\bow)\sim J(\bow)$ if we have 
\begin{align*}
|I(\bow) - J(\bow) |\le C (1+P(\nor{\bow-\boe_3}{H^1(\R)})) \nor{\bow-\boe_3}{H^2}^2
\end{align*}
for some $C>0$ and a polynomial $P$. 
Note that $\sim$ is an equivalent relation. Typical examples are
\begin{align}
\int_\R P(w_1,w_2,w_3-1) |\bow''|^j dx \sim 0\qquad \text{for } j=0,1,2,\\
\int_\R P(w_1,w_2,w_3-1) |\bow'|^k |\bow''|dx\sim 0\qquad \text{for } k=0,1,2,3. \label{eq:WPX2}\\
\int_\R P(w_1,w_2,w_3-1) |\bow'|^k dx\sim 0\qquad \text{for } k=0,1,...,6. \label{eq:WPY1}
\end{align}
where each $P$ is a polynomial such that the integrand is higher then or equal to quadratic order.
In fact, \eqref{eq:WPX2} for $k=3$ and \eqref{eq:WPY1} for $k=6$ follows from
\begin{align*}
\int_\R P(w_1,w_2,w_3-1) |\bow'|^3 |\bow''|dx&\le C\nor{ P(w_1,w_2,w_3-1)}{L^\infty}
\nor{\bow'}{L^6}^3 \nor{\bow''}{L^2} \\
&\le C (1+\nor{\bow-\boe_3}{H^1}^K) \nor{\bow'}{L^2}^2 \nor{\bow''}{L^2}^2,\\
\int_\R P(w_1,w_2,w_3-1) |\bow'|^6dx&\le C\nor{ P(w_1,w_2,w_3-1)}{L^\infty}
\nor{\bow'}{L^6}^6 \\
&\le C (1+\nor{\bow-\boe_3}{H^1}^K) \nor{\bow'}{L^2}^4 \nor{\bow''}{L^2}^2
\end{align*}
for some $K\ge 0$ by the Gagliardo--Nirenberg inequality.\par
We rewrite the equation as
\[
\rd_t \bow = \bow\times  \bow'' - 2(b-\eps) w_1 \bow' + F(\bow-\boe_3)
\]
where $F(\bow-\boe_3)$ involves no derivatives of $\bow$. By noting this, we have
%
%
% \begin{lem}
% \end{lem}
\begin{align}
&\rd_t \left(
\nor{\bow''}{L^2}^2 -\frac 54 \nor{\bow'}{L^4}^4
\right)\notag\\
%&= 2\inp{\bow_t''}{\bow''} - 5 \inp{|\bow'|^2\bow}{\bow_t'}\\
&\sim -4 \int \bow''\cdot(b-\eps)(w_1\bow')''dx +10 \int |\bow'|^2\bow'\cdot (b-\eps)(w_1\bow')' dx\notag\\
&\qquad +2 \int \bow''\cdot \bF'' dx- 5 \int |\bow'|^2\bow'\cdot \bF'' dx\notag\\
&\sim \int-4 (b-\eps) (w_1''\bow'\cdot\bow'' + 2w_1'|\bow''|^2 +w_1\bow'''\cdot \bow''
) \notag\\
% &\qquad +\int10(b-\eps)|\bow'|^2\bow' (w_1'\bow' + w_1\bow'')\\
&\sim \int-4 (b-\eps) w_1''\bow'\cdot\bow'' -6 (b-\eps) w_1' |\bow''|^2 \notag\\
&\eqqcolon -4B_1(\bow) - 6B_2(\bow). \label{eq:WP14}
\end{align}
Our aim is to cancel $B_1$ and $B_2$ by adding appropriate correction terms. %First, 
We introduce
\begin{align}
\boI (\bow) &\coloneqq \int_\R (b-\eps) \bow''\cdot \boe_1\times \bow' dx,\label{eq:Sep6}\\
\boJ (\bow) &\coloneqq \int_\R (b-\eps) |\bow'|^2\bow\cdot \boe_1\times \bow'dx.\label{eq:Sep7}
\end{align}
%The first one is the following:
\begin{lem}
\begin{align}\label{eq:WP12}
\rd_t \boI(\bow) &\sim 6B_1(\bow) -B_2(\bow),\\
\rd_t \boJ (\bow) &\sim -2B_1(\bow) -2 B_3(\bow).\label{eq:WP13}
\end{align}
where
\[
B_3(\bow) \coloneqq \int (b-\eps) [\bow \cdot \bow'\times \bow''] [\boe_1\cdot \bow\times \bow''].
\]
\end{lem}
\begin{proof}
For \eqref{eq:WP12}, we have
\begin{align*}
\rd_t \boI &\sim \int (b-\eps) (\bow\times \bow'')''\cdot \boe_1\times \bow' +
\int (b-\eps) \bow''\cdot \boe_1\times  (\bow\times \bow'')'\\
&\eqqcolon I_1+I_2.
\end{align*}
By integration by part, we obtain
\begin{align*}
I_1 &\sim \int -(b-\eps) (\bow\times\bow'')' \cdot\boe_1\times \bow''\\
&\sim \int -(b-\eps) (\bow'\times \bow'' + \bow\times \bow''')\cdot\boe_1\times \bow''\\
&\sim \int -(b-\eps)(w_1'|\bow''|^2 -w_1''\bow'\cdot\bow'' + w_1\bow''\cdot \bow'''-w_1'''\bow\cdot \bow'')\\
% &\sim \int -(b-\eps)(w_1'|\bow''|^2 -w_1''\bow'\cdot\bow'' -\frac 12  w_1' |\bow''|^2
% +w_1''' |\bow'|^2)\\
&\sim \int - \frac 12(b-\eps ) w_1'|\bow'|^2 +3 (b-\eps) w_1'' \bow'\cdot \bow'',
\end{align*}
where we used
\[
\bow\cdot \bow''= -|\bow'|^2,\qquad \bow\cdot \bow'''= -3 \bow'\cdot \bow'',
\]
which follow from $|\bow|\equiv 1$. Similarly, we have
\begin{align*}
I_2 &= \int (b-\eps) \bow''\cdot \boe_1 \times  (\bow'\times \bow'' + \bow\times \bow''')\\
&\sim \int (b-\eps) \bow'' \cdot (w_1''\bow' - w_1'\bow'' + w_1'''\bow - w_1\bow''')\\
% &\sim \int (b-\eps) (w_1'' \bow'\cdot\bow'' - w_1' |\bow''|^2 - w_1''' |\bow'|^2 + \frac 12 w_1'|\bow''|^2)\\
&\sim \int (b-\eps) (3w_1''\bow'\cdot \bow'' - \frac 12 w_1' |\bow''|^2).
\end{align*}
Therefore, combining these yields \eqref{eq:WP12}. Next, 
% \begin{align*}
% \rd_t \boI_1 \sim \int (b-\eps) (-w_1'|\bow''|^2 + 6w_1''\bow'\cdot \bow'')
% \end{align*}
% \begin{align*}
% \boI_2 = \int(b-\eps) |\bow'|^2\bow\cdot \boe_1\times \bow'
% \end{align*}
\begin{align*}
\rd_t\boJ &= \int (b-\eps) 2 (\bow_t' \cdot \bow') \bow\cdot \boe_1\times \bow' 
+ \int (b-\eps) |\bow'|^2\bow_t \cdot \boe_1\times \bow' \\
&\qquad + \int (b-\eps) |\bow'|^2\bow\cdot \boe_1\times \bow_t'\\
&\eqqcolon J_1+J_2+J_3.
\end{align*}
Each term can be computed as
\begin{align*}
I_1 &\sim \int 2(b-\eps) 
[(\bow\times \bow'')' \cdot \bow'] [\bow\cdot \boe_1\times \bow'] 
%\\&
\sim \int -2(b-\eps) [\bow\times \bow'' \cdot \bow'] [\bow\cdot \boe_1\times \bow''],\\
% \end{align*}
% \begin{align*}
I_2 &\sim \int (b-\eps) |\bow'|^2 (\bow\times \bow'') \cdot \boe_1 \times \bow' \sim 0. \\
% \end{align*}
% \begin{align*}
I_3&\sim \int (b-\eps) |\bow'|^2\bow \cdot \boe_1\times (\bow\times \bow'')'\\
&= \int (b-\eps) |\bow'|^2\bow \cdot \boe_1\times (\bow'\times \bow'' + \bow\times \bow''')\\
&\sim \int (b-\eps) |\bow'|^2\bow \cdot (w_1'''\bow - w_1 \bow''')\\
&\sim \int (b-\eps) (-2w_1'' \bow'\cdot \bow'' + 3 w_1|\bow'|^2 \bow'\cdot \bow'' )\sim -2 \int (b-\eps) w_1'' \bow'\cdot \bow''.
\end{align*}
Therefore, we obtain \eqref{eq:WP13}.
% \begin{align*}
% \rd_t\boI_2 \sim \int 2(b-\eps) w_1'' \bow'\cdot \bow''.
% \end{align*}
\end{proof}
To cancel $B_3$, we further introduce
\begin{align}\label{eq:Sep8}
\boK(\bow) = \int(b-\eps) w_1\bow\times\bow'' \cdot \bow'.
\end{align}
\begin{lem}
\[
\rd_t \boK(\bow) \sim B_1(\bow)+B_2(\bow) -B_3(\bow)
\]
\end{lem}
\begin{proof}
\begin{align*}
\rd_t\boI 
&= \int w_{1t} \bow\times\bow'' \cdot \bow' 
+ w_1\bow_t \times \bow''\cdot \bow'
+(b-\eps) w_1 \bow\times \bow_t'' \cdot \bow'
+(b-\eps) w_1 \bow\times \bow'' \cdot \bow_t'\\
&= I_1+I_2+I_3+I_4.
\end{align*}
Now we substitute \eqref{E2.4Aug}. 
It is straightforward to see that the terms where $\rd_t\bow$ is replaced by $-2(b-\eps)w_1\bow'+\bF$ are good terms. Hence
\begin{align*}
I_1 &\sim \int (b-\eps) [(\bow\times \bow'')\cdot \boe_1] [\bow\times \bow''\cdot \bow']\\
I_2 &\sim \int (b-\eps) w_1(\bow\times \bow'')\times \bow'' \cdot \bow'
% \sim \int (b-\eps) w_1 (\bow\cdot \bow'')\bow'' \cdot \bow'\\
% &
= \int -(b-\eps) w_1|\bow'|^2\bow''\cdot \bow' \sim 0,
\\
% &\sim \int \frac 14(b-\eps) w_1' |\bow'|^4
%\sim 0
% I_2 &\sim \int (b-\eps) w_1(\bow\times \bow'')\times \bow'' \cdot \bow'\\
% &\sim \int (b-\eps) w_1(-|\bow''|^2\bow + (\bow\cdot \bow'')\bow'') \cdot \bow'\\
% &\sim \int -(b-\eps) w_1|\bow'|^2\bow''\cdot \bow' \\
% &\sim \int \frac 14(b-\eps) w_1' |\bow'|^4\sim 0
% \end{align*}
% \begin{align*}
I_3 &\sim \int (b-\eps) w_1\bow\times (\bow\times \bow'')'' \cdot \bow'\\
&\sim \int -(b-\eps) \bow\times (\bow\times \bow'')' \cdot (w_1\bow')'\\
&\sim \int -(b-\eps) (\bow\times (\bow'\times \bow'') - \bow\times (\bow\times \bow''')) \cdot (w_1'\bow'+w_1\bow'')\\
&\sim \int -(b-\eps) (-|\bow'|^2\bow' - 3(\bow'\cdot\bow'')\bow + \bow''')\cdot (w_1'\bow'+w_1\bow'')\\
&\sim \int -(b-\eps) \bow'''\cdot (w_1'\bow'+w_1\bow'')\\
&\sim \int (b-\eps) w_1'' \bow'\cdot \bow'' + \frac 32 (b-\eps) w_1' |\bow''|^2 \\
% \end{align*}
% \begin{align*}
I_4 &\sim \int (b-\eps) w_1\bow\times \bow'' \cdot (\bow\times \bow'')'\\
&\sim \int -\frac12 (b-\eps) w_1' (\bow\times \bow'')\cdot (\bow\times \bow'')\\
&\sim \int -\frac12 (b-\eps) w_1' (|\bow''|^2 -(\bow\cdot \bow'')^2)
%\\
% &= \int -\frac 12(b-\eps) w_1' |\bow''|^2
% + \frac 12(b-\eps) w_1' |\bow'|^4\\
\sim \int-\frac 12 (b-\eps) w_1'|\bow''|^2.
\end{align*}
%Hence, 
Combining the above implies the conclusion.
% \begin{align*}
% \rd_t \boI_1 &\sim \int (b-\eps) [(\bow\times \bow'')\cdot \boe_1] [\bow\times \bow''\cdot \bow']\\
% &\qquad +\int (b-\eps) (w_1''\bow'\cdot \bow'' + w_1'|\bow''|^2)
% \end{align*}
\end{proof}

Combining \eqref{eq:WP14}, \eqref{eq:WP12}, and \eqref{eq:WP13}, we have obtained the following:
\begin{prop}
Let
% \[
% \boE_2^* (\bow) \coloneqq \nor{\bow''}{L^2}^2 - \frac 54\nor{\bow'}{L^4}^4 - 6\boI(\bow)
% %\int (b-\eps )\bow''\cdot \boe_1\times \bow' 
% +20 \boJ(\bow).
% %\int(b-\eps) |\bow'|^2 \bow\cdot \boe_1\times \bow'.
% \]
\begin{align}\label{eq:Sep5}
\boE_2^* (\bow) \coloneqq \nor{\bow''}{L^2}^2 - \frac 54\nor{\bow'}{L^4}^4 + \boI(\bow)+ \frac 52 \boJ (\bow) -5 \boK(\bow)
\end{align}
Then we have $\boE_2^*(\bow)\sim 0$.
% \begin{align}
% \rd_t \left(
% \nor{\bow''}{L^2}^2 - \frac 54\nor{\bow'}{L^4}^4 - 6\int (b-\eps )\bow''\cdot \boe_1\times \bow' +20 \int(b-\eps) |\bow'|^2 \bow\cdot \boe_1\times \bow'
% \right) \sim 0.
% \end{align}
In particular, it holds that
%we write the LHS as $\mathcal{I}$, then 
\begin{align}\label{eq:WP15}
\mathcal{I}(\bow(t)) \le  \mathcal{I}(\bow(s)) + \int_s^t (1+P(\nor{\bow(t')-\boe_3}{H^1})) \nor{\bow(t')-\boe_3}{H^2}^2 dt'
\end{align}
for $T_{\max}^-<s<t<T_{\max}^+$.
\end{prop}
\begin{proof}[Proof of Propositoin \ref{P:wp2}]
We only show $T_{\max}^+=\infty$. 
Suppose to the contrary that $T_{\max}^+<\infty$. Let $T\in (0,T_{\max}^+)$ be a number fixed later. For $t\in (T,T_{\max}^+)$, 
Young's inequality yields
\begin{align*}
    \nor{\bow''(t)}{L^2}
    &\le \nor{\bow''(t)}{L^2}-\frac 54\nor{\bow'(t)}{L^4}^4\\
    &\le \boE_2^*(\bow(t)) + 6 \boI(\bow(t)) -20 \boJ(\bow(t))\\
    &\le \mathcal{E}_2^*(\bow(t)) + \frac 12 \nor{\bow''(t)}{L^2} + C P(\nor{\bow(t)-\boe_3}{H^1}),
\end{align*}
which implies
%for arbitrary $\nu>0$. 
% Hence, by letting $M\coloneqq \nor{\bow-\hat\boe_3}{L^\infty_tH^1_x}$
% , we have
\[
\nor{\bow''(t)}{L^2} \le C_0 \mathcal{E}_2^*(\bow(t))+ C_1
\]
with $C_0,C_1>0$ depending only on $\nor{\bow-\hat\boe_3}{L^\infty_tH^1_x}$. Therefore, \eqref{eq:WP15} gives
\[
\nor{\bow''(t)}{L^2} \le \boE_2^*(\bow(T)) + C_1 + C_2 (t-T) \nor{\bow''}{L^\infty((T,t)\colon L^2_x)}
\]
for $T<t<T_{\max}^+$. By choosing $T$ such that $C_2(T_{\max}^+-T)<\frac 12$, we obtain
\[
\nor{\bow''}{L^\infty((T,t)\colon L^2_x)} \le C (\boE_2^*(\bow(T)) + 1),\qquad t\in (T,T_{\max}^+),
\]
with $C>0$ independent of $t$. This contradicts \eqref{eq:WP1}.
\end{proof}

\subsection*{Acknowledgements}
S.I. acknowledges financial support from the Natural Sciences and Engineering Research Council of Canada (NSERC) under Grant No. 371637-2025. Part of this work was discussed at the Research Institute for Mathematical Sciences-Kyoto,  during Fall 2025. S.I and I.S would like to thank the staff of RIMS for their hospitality and for providing a stimulating environment for collaboration. 
C.M. is indebted to Stavros Komineas for valuable discussions on the subject matter.
I.S. is supported by JSPS KAKENHI Grant Number JP26K17016.

\bibliography{BW}

\bigskip
{\small
\noindent
\textsc{Slim Ibrahim}\\
Department of Mathematics and Statistics, University of Victoria, BC-Canada\\
\textit{Email:} ibrahims@uvic.ca

\medskip

\noindent
\textsc{Christof Melcher}\\
Chair of Applied Analysis, RWTH Aachen
University, Aachen, 52026, Germany
\\
\textit{Email:} melcher@rwth-aachen.de

\medskip

\noindent
\textsc{Ikkei Shimizu}\\
Department of Mathematics, Graduate School of Science, Kyoto University, Kyoto, 606-8502, Japan\\
\textit{Email:} shimizu.ikkei.8s@kyoto-u.ac.jp
}

\end{document}